\documentclass[10pt]{article}

\usepackage[all,cmtip]{xy}

\usepackage{setspace}
\usepackage{bbm,youngtab}
\usepackage{amscd,mathrsfs}
\usepackage{amssymb,amsmath,dsfont,amsfonts,amsthm}
\usepackage{enumerate}
\usepackage[mathscr]{eucal}
\usepackage{stmaryrd} 
\usepackage[cbgreek]{textgreek} 
\usepackage{booktabs}
\usepackage{multirow}
\usepackage{hhline} 

\usepackage{rotating}
\usepackage{ftnxtra}
\usepackage{fnpos}
\usepackage[titletoc,toc,title]{appendix}
\usepackage{graphicx}
\usepackage{epsfig}
\usepackage{epstopdf}
\DeclareGraphicsExtensions{.pdf,.png,.jpg,.eps}
\usepackage{upgreek}
\usepackage{tikz-cd}

\usepackage{authblk}
\usepackage{bm}
\usepackage{caption}
\usepackage{mathtools}
\usepackage{soul}            
\usepackage{xcolor}
\usepackage[colorlinks,citecolor=blue,linkcolor=blue]{hyperref}

\usepackage{natbib}

\numberwithin{equation}{section}

\newtheorem{lemma}{Lemma}
\newtheorem{corollary}{Corollary}[lemma] 

\theoremstyle{definition}
\newtheorem{remark}{Remark}

\theoremstyle{definition}
\newtheorem{example}{Example}

\title{
\vspace{-3cm}
\textbf{Second-Order Compatible-Strain Mixed Finite Elements for $2$D Compressible Nonlinear Elasticity}}

\author[1,2]{Mohsen Jahanshahi\thanks{Corresponding author, e-mail: jahanshahi@sharif.edu}}

\author[1]{Damiano Pasini}
\affil[1] {\small \textit{Department of Mechanical Engineering, McGill University, Montreal, Quebec H3A 0C3, Canada}}
\affil[2]{\small \textit{Department of Civil Engineering, School of Science and Engineering, Sharif University of Technology, International Campus, Kish Island, P.O. Box 76417-76655, Iran}}

\author[3,4]{Arash Yavari}
\affil[3]{\small \textit{School of Civil and Environmental Engineering, Georgia Institute of Technology, Atlanta, GA 30332, USA}}
\affil[4]{\small \textit{The George W. Woodruff School of Mechanical Engineering, Georgia Institute of Technology, Atlanta, GA 30332, USA}}

\begin{document}

\maketitle

\begin{abstract}
In recent years, a new class of mixed finite elements---compatible-strain mixed finite elements (CSMFEs)---has emerged that uses the differential complex of nonlinear elasticity. Their excellent performance in benchmark problems, such as numerical stability for modeling large deformations in near-incompressible solids, makes them a promising choice for solving engineering problems. Explicit forms exist for various shape functions of first-order CSMFEs. In contrast, existing second-order CSMFEs evaluate shape functions using numerical integration. In this paper, we formulate second-order CSMFEs with explicit shape functions for the displacement gradient and stress tensor. Concepts of vector calculus that stem from exterior calculus are presented and used to provide efficient forms for shape functions in the natural coordinate system. Covariant and contravariant Piola transformations are then applied to transform the shape functions to the physical space. Mid-nodes and pseudo-nodes are used to enforce the continuity constraints for the displacement gradient and stress tensor over the boundaries of elements. The formulation of the proposed second-order CSMFEs and technical aspects regarding their implementation are discussed in detail. Several benchmark problems are solved to compare the performance of CSMFEs with first-order CSMFEs and other second-order elements that rely on numerical integration. It is shown that the proposed CSMFEs are numerically stable for modeling near-incompressible solids in the finite strain regime.
\end{abstract}

\begin{description}
\item[Keywords:] Mixed finite elements, finite element exterior calculus, compatible-strain finite elements, nonlinear elasticity.
\end{description}

\tableofcontents

\section{Introduction}
\label{sec:1}

The conventional theory of the Finite Element Method (FEM) is well-established \citep{BAT96,ZIE05,HUG00,BEL05,WRI08}. A well-known feature of displacement-based finite elements is the discontinuity of the displacement gradient and stress tensor across the element boundaries. Conventional displacement-based finite element methods often lead to inaccurate results when used for modeling incompressible or nearly incompressible solids. This is because the space of the displacement gradient derived from the displacement field is not sufficiently large to encompass the solution space \citep{GLO84,GLO88,SIM91}. Additionally, they are unable to accurately solve plasticity problems due to the isochoric nature of plastic deformations. To address these issues, mixed finite elements have been proposed as an alternative to displacement-based FEMs.

The theoretical background for linear mixed FEMs is discussed in detail in \citep{WAS75,BRA07,BRE07}. The first hybrid elements for linear elasticity were introduced by \citet{PIA64}. In his formulation, he defined certain functions in terms of the generalized displacements to satisfy the equilibrium equations and maintain the compatibility between adjacent elements. \citet{PIA84} then developed an efficient and accurate finite element based on the Hellinger-Reissner variational principle. In some of the early works on modeling finite strain problems, the deformation gradient was first multiplicatively decomposed into a volumetric and an isochoric part, and then the displacement, pressure and dilation were treated as independent fields \citep{SIM85,SIM88a,SIM88b,SIM91}. Enhanced strain methods or the methods of incompatible modes were designed and developed based on mixed formulation \citep{SIM90,SIM92,SIM93}. Various stabilization techniques were later proposed to address the deficiencies attributed to these methods and to improve their performance \citep{KOR96,WRI96,GLA97,REE00}. For quasi-incompressible finite elasticity, \citet{SHR11} developed a new mixed finite element method, which is based on the approximations of the minors of the deformation gradient. \citet{BON15} extended the formulation by considering the deformation gradient, its adjoint and its determinant as independent kinematic variables. \citet{NEU21} introduced three mixed formulations for nonlinear elasticity with the following independent fields: (i) displacement, deformation gradient, and the first Piola–Kirchhoff stress; (ii) displacement, the right Cauchy–Green strain, and the second Piola–Kirchhoff stress; and (iii) displacement, deformation gradient, the right Cauchy–Green strain, and both the first and second Piola–Kirchhoff stresses. \citet{SIM88a,SIM88b} developed a framework for finite strain elastoplasticity based on a mixed formulation and the multiplicative decomposition of the deformation gradient into elastic and plastic parts. More recently, \citet{JAH15} introduced an integration algorithm for $J_{2}$ plasticity that uses a three-field Hu-Washizu principle and a similar decomposition for the deformation gradient. \citet{KHO16} employed a mixed formulation using the reversed multiplicative decomposition of the deformation gradient into plastic and elastic parts to extend the classical theory of plasticity to nano-structures.

Over the years, mixed finite element methods have been successfully used to solve a wide range of engineering problems. However, their performance is not always satisfactory. For example, undesired hourglass modes appear when low order elements are used for the analysis of bending dominated problems in the near-incompressible regime and for modeling isochoric plastic deformations in metals. Hourglass or zero-energy modes are observed when the eigenvalues corresponding to the bending modes of an element vanish during the course of loading \citep{GLA97,REE02,REE05}. Numerous stabilization techniques have been proposed to address the issue of hourglassing. \citet{JAQ86} introduced a method to filter hourglass modes from the global solution. Their approach involved obtaining the solution using underintegration and then eliminating the instabilities present in the solution through a post-processing operation. Their method is, however, limited to elastic problems. An alternative approach resorts to a displacement interpolation which encompasses a linear and orthogonal component. The latter is used to derive a stabilization matrix \citep{KOS78,BEL84}. Although it is possible to compute the stiffness parameters for the stabilization matrix from the equivalence of mixed methods and reduced integration procedure \citep{KOS78,MAL78}, it is not straightforward to apply the technique to nonlinear problems \citep{BEL93,REE98,WRI08}. Numerical stability issues in the near-incompressible regime have been addressed using enhanced strain or enhanced assumed strain elements, originally proposed by \citet{SIM90}. Their formulation was later extended to large deformations \citep{SIM92,SIM93}. The enhanced strain elements typically show good performance. However, it has been shown that they can become unstable under compression \citep{WRI96}. This issue was partially resolved via modified shape functions \citep{KOR96,GLA97,REE00}. An in-depth study of enhanced strain-based elements shows that the deformation gradient in these elements is enhanced through internal degrees of freedom, which however are not compatible with the displacement field approximated by shape functions.

Recently, a new approach for formulating mixed finite elements for linear elasticity, which is based on the differential complexes of linear ealsticity, has gained popularity \citep{ARN02,ARN06,ARN18}. \citet{ARN06} discretized the linear elasticity complex in such a way that the discrete complex preserves all the topological structures of the linear elasticity complex. They generalized the concept and developed mixed methods that are numerically stable for linear operators applicable to certain classes of differential complexes. \citet{ARN06,ARN10,ARN18} proposed the Finite Element Exterior Calculus (FEEC). FEEC utilizes the tools from geometry and topology to develop numerical techniques for a class of Partial Differential Equations (PDEs) that are generally hard to solve using conventional finite element methods. In the context of nonlinear elasticity, \citet{ANG15,ANG16} introduced a differential complex that is suitable for describing the kinematics and kinetics of large deformations. They expressed the nonlinear elasticity complex in terms of displacement, displacement gradient and the first Piola-Kirchhoff stress tensor. They also introduced a corresponding Hilbert complex of nonlinear elasticity. The nonlinear elasticity Hilbert complex is then discretized so that the discrete complex inherits the topological properties of the original complex. Similar discrete complexes can also be obtained using edge and face finite elements for vector fields that were introduced by \citet{RAV77}, \citet{BRE85} and \citet{NED80,NED86}. \citet{ANG17} proposed a mixed formulation for solving two-dimensional compressible elasticity problems in the finite strain regime by considering the displacement, the displacement gradient and the first Piola-Kirchhoff stress tensor as the independent fields of a Hu-Washizu functional. Since the displacement gradient in their approach satisfies the classical Hadamard jump condition for the compatibility of non-smooth strain fields, this method was named Compatible-Strain Mixed Finite Element (CSMFE). The performance of CSMFE for bending dominated and near-incompressible problems was demonstrated through several numerical examples. Recently, \citet{DHA22a} proposed a mixed variational approach in nonlinear elasticity using Cartan's moving frames and implemented it via the application of FEEC. They extended their formulation to 3D nonlinear elasticity \citep{DHA22b}. \citet{SHO18} utilized a four-field Hu-Washizu functional to extend CSMFE to the solution of incompressible nonlinear elasticity problems. In their work, the displacement, the displacement gradient, the first Piola-Kirchhoff stress tensor and pressure were considered as independent fields. \citet{SHO19} employed three and four-field functionals to develop new compatible-strain finite element methods applicable to the solution of 3D compressible and incompressible nonlinear elasticity problems. In these and many similar works, a comprehensive study of the performance of CSMFEs in comparison with other classical formulations is missing. \citet{JAH22} developed a compatible mixed finite element method using the idea of the mid-nodes for the imposition of the continuity constraints over the boundaries of elements. It was shown that CSMFEs have excellent performance both in terms of the computational efficiency and the ability to address the deficiencies that typically emerge in the conventional first-order finite elements. Moreover, the performance of the proposed CSMFE was compared to that of enhanced strain-based elements through several numerical examples.

In this paper, a second-order mixed finite element method is introduced for the analysis of two-dimensional compressible solids in the finite strain regime. A three-field Hu-Washizu functional, with displacement, displacement gradient and stress tensor considered as independent fields, is used. A novel aspect of the proposed formulation is that the shape functions that are used to interpolate the displacement gradient and stress tensor are derived explicitly. This is in contrast to the existing formulations of the second-order compatible strain elements where the shape functions are evaluated using numerical integration \citep{ANG16}. The explicit forms help to study the variation of shape functions over the domain of the elements and decide on the form of the global shape functions.
Unlike previous formulations, which are based on exterior calculus \citep{JAH22}, the concepts in vector calculus are first developed using their equivalents in exterior calculus and are then applied to formulate the shape functions for displacement gradient and stress tensor. The shape functions are derived in the natural coordinate system based on the specific requirements that they need to satisfy in the physical space. They are then transformed to the physical space using the covariant and contravariant Piola transformations as discussed in detail in \citep{ROG10, AZN22}. Certain constraints are imposed on the displacement gradient and stress tensor to satisfy the required continuity conditions across the boundaries of the elements \citep{ANG16,JAH22}. In order to facilitate the application of these constraints, the degrees of freedom that are associated to the edges are assigned to the mid-nodes of these edges. For the degrees of freedom associated to the element itself, pseudo-nodes are defined at the centers of the elements. The concept of mid-nodes and pseudo-nodes, which helps to provide an elegant implementation for the formulation, is another novel aspect of the proposed formulation and the earlier work of \citet{JAH22}.

This paper is organized as follows. The space of polynomial vector fields that provide the basis for deriving shape functions is discussed in \S\ref{sec:2}. The order of polynomials and the independent variables are specifically chosen for two-dimensional problems. The derivation of second-order shape functions is described in \S\ref{sec:3}. The explicit forms of shape functions are provided in Appendix \ref{app:A}. The formulation of the second-order CSMFEs is discussed in \S\ref{sec:4}.  Technical aspects of the mixed formulation, such as assigning different degrees of freedom to nodes and pseudo-nodes, efficient methods for implementing the continuity constraints over the edges and minimum dimensions of sub-matrices to preserve the rank of the overall stiffness matrix, are also discussed. Several numerical examples are presented in \S\ref{sec:5}. The results of the proposed formulation are compared with those of the existing mixed formulations. Finally, concluding remarks are given in \S\ref{sec:6}.

\section{Polynomial Vector Spaces}
\label{sec:2}

The shape functions that approximate the displacement gradient and the stress tensor belong to the spaces $\mathcal{P}_{r}\left(\otimes^{2}T\mathbb{R}^{2}\right)$, $\mathcal{P}^{-}_{r}\left(\otimes^{2}T\mathbb{R}^{2}\right)$ and $\mathcal{P}^{\ominus}_{r}\left(\otimes^{2}T\mathbb{R}^{2}\right)$ of polynomial tensors. These spaces can be constructed using the polynomial vector spaces $\mathcal{P}_{r}\left(T\mathbb{R}^{2}\right)$, $\mathcal{P}^{-}_{r}\left(T\mathbb{R}^{2}\right)$ and $\mathcal{P}^{\ominus}_{r}\left(T\mathbb{R}^{2}\right)$. The spaces $\mathcal{P}^{-}_{r}\left(T\mathbb{R}^{2}\right)$ and $\mathcal{P}^{\ominus}_{r}\left(T\mathbb{R}^{2}\right)$ are intermediate between $\mathcal{P}_{r-1}\left(T\mathbb{R}^{2}\right)$ and $\mathcal{P}_{r}\left(T\mathbb{R}^{2}\right)$, meaning that $\mathcal{P}_{r-1}\left(T\mathbb{R}^{2}\right)\subset\mathcal{P}^{-}_{r}\left(T\mathbb{R}^{2}\right)\subset\mathcal{P}_{r}\left(T\mathbb{R}^{2}\right)$ and $\mathcal{P}_{r-1}\left(T\mathbb{R}^{2}\right)\subset\mathcal{P}^{\ominus}_{r}\left(T\mathbb{R}^{2}\right)\subset\mathcal{P}_{r}\left(T\mathbb{R}^{2}\right)$. Given a degree $r$, the rows of the displacement gradient can be chosen from either $\mathcal{P}_{r}\left(T\mathbb{R}^{2}\right)$ or $\mathcal{P}^{-}_{r}\left(T\mathbb{R}^{2}\right)$, while the rows of the stress tensor belong to either $\mathcal{P}_{r}\left(T\mathbb{R}^{2}\right)$ or $\mathcal{P}^{\ominus}_{r}\left(T\mathbb{R}^{2}\right)$. The details of constructing the spaces $\mathcal{P}_{r}\left(T\mathbb{R}^{2}\right)$, $\mathcal{P}^{-}_{r}\left(T\mathbb{R}^{2}\right)$ and $\mathcal{P}^{\ominus}_{r}\left(T\mathbb{R}^{2}\right)$, and then restricting them to a given simplex $\mathscr{T}$ are discussed in this section.

We denote the space of homogeneous polynomials in $n$ variables of degree $r$ by $\mathcal{H}_{r}\left(\mathbb{R}^{n}\right)$, which has the dimension \citep{ARN06}:
\begin{equation} \label{eqn:2.1}
\text{dim}\,\mathcal{H}_{r}\left(\mathbb{R}^{n}\right)=\binom{n+r-1}{n-1}\,.
\end{equation}
The space $\mathcal{P}_{r}\left(\mathbb{R}^{n}\right)$ of polynomials in $n$ variables of degree at most $r$ can be constructed by the direct sum of the spaces $\mathcal{H}_{k}\left(\mathbb{R}^{n}\right),k=1,\ldots,r$ as:
\begin{equation}\label{eqn:2.2}
\mathcal{P}_{r}\left(\mathbb{R}^{n}\right)=\overset{r}{\underset{k=0}{\bigoplus}}\,\mathcal{H}_{k}\left(\mathbb{R}^{n}\right)\,.
\end{equation}
Its dimension is given by \citep{ARN06}:
\begin{equation}\label{eqn:2.3}
\text{dim}\,\mathcal{P}_{r}\left(\mathbb{R}^{n}\right)=\binom{n+r}{n}\,.
\end{equation}
The space $\mathcal{H}_{r}\left(T\mathbb{R}^{n}\right)$ of homogeneous polynomial vector fields in $n$ variables of degree $r$ and the space $\mathcal{P}_{r}\left(T\mathbb{R}^{n}\right)$ of polynomial vector fields in $n$ variables of degree at most $r$ can be defined as extensions of $\mathcal{H}_{r}\left(\mathbb{R}^{n}\right)$ and $\mathcal{P}_{r}\left(\mathbb{R}^{n}\right)$, which have, respectively, the dimensions \citep{ARN06}:
\begin{equation}\label{eqn:2.4}
\text{dim}\,\mathcal{H}_{r}\left(T\mathbb{R}^{n}\right)=n\binom{n+r-1}{n-1}\,,\qquad
\text{dim}\,\mathcal{P}_{r}\left(T\mathbb{R}^{n}\right)=n\binom{n+r}{n}\,.
\end{equation}
An extension of \eqref{eqn:2.2} to vector spaces shows that the space $\mathcal{P}_{r}\left(T\mathbb{R}^{n}\right)$ can be constructed as:
\begin{equation}\label{eqn:2.5}
\mathcal{P}_{r}\left(T\mathbb{R}^{n}\right)=\mathcal{P}_{r-1}\left(T\mathbb{R}^{n}\right)+\mathcal{H}_{r}\left(T\mathbb{R}^{n}\right)\,.
\end{equation}
We can now define the polynomial vector spaces $\mathcal{P}^{-}_{r}\left(T\mathbb{R}^{n}\right)$ and $\mathcal{P}^{\ominus}_{r}\left(T\mathbb{R}^{n}\right)$ that are intermediate between the spaces $\mathcal{P}_{r-1}\left(T\mathbb{R}^{n}\right)$ and $\mathcal{P}_{r}\left(T\mathbb{R}^{n}\right)$.

In the following, we restrict the discussion to $\mathbb{R}^{2}$. As a result, the dimensions of $\mathcal{H}_{r}\left(T\mathbb{R}^{2}\right)$ and $\mathcal{P}_{r}\left(T\mathbb{R}^{2}\right)$ are $2\left(r+1\right)$ and $\left(r+1\right)\left(r+2\right)$, respectively. If $\mathbf{x}$ is the position vector of a given point in $\mathbb{R}^{2}$ and $p\left(x_{1}, x_{2}\right)\in\mathcal{P}_{r}\left(\mathbb{R}^{2}\right)$ is a polynomial at that point, then the operator $\mathbf{K}:\mathcal{P}_{r}\left(\mathbb{R}^{2}\right)\rightarrow\mathcal{P}_{r+1}\left(T\mathbb{R}^{2}\right)$ is defined as \citep{JAH22}:
\begin{equation}\label{eqn:2.6}
\mathbf{K}\left(p\right)=\left(p{\epsilon_{j}}^{i}x^{j}\right)\mathbf{e}_{i}\,,
\end{equation}
where, the permutation symbol ${\epsilon_{i}}^{j}$ takes the value of either 1 for $i=1,j=2$, $-1$ for $i=2,j=1$, and zero otherwise. In addition, for the vector $\mathbf{v}\left(x_{1},x_{2}\right)\in\mathcal{P}_{r}\left(T\mathbb{R}^{2}\right)$, the operator $\mathbf{k}:\mathcal{P}_{r}\left(T\mathbb{R}^{2}\right)\rightarrow\mathcal{P}_{r+1}\left(\mathbb{R}^{2}\right)$ is defined as:
\begin{equation}\label{eqn:2.7}
\mathbf{k}\left(\mathbf{v}\right)=\delta_{ij}v^{i}x^{j}\,,
\end{equation}
where $\delta_{ij}$ is the Kronecker delta. Evidently, at point $\mathbf{x}$, the operator $\mathbf{k}$ applied to $\mathbf{v}$ leads to the inner product of $\mathbf{v}$ and $\mathbf{x}$, i.e. $\mathbf{k}\left(\mathbf{v}\right)=\mathbf{v}\cdot\mathbf{x}$. It is straightforward to verify that the following property holds for the operators $\mathbf{K}$ and $\mathbf{k}$: $\mathbf{k}\circ\mathbf{K}\left(p\right)=p\epsilon_{ij}x^{i}x^{j}=0$. Using the definition of the operator $\mathbf{K}$ in \eqref{eqn:2.6}, it was shown in \citep{JAH22} that the space $\mathcal{P}_{r}\left(T\mathbb{R}^{2}\right)$ can be decomposed as:
\begin{equation}\label{eqn:2.8}
    \mathcal{P}_{r}\left(T\mathbb{R}^{2}\right)=
    \mathbf{K}\left(\mathcal{P}_{r-1}\left(\mathbb{R}^{2}\right)\right)\oplus
    \mathbf{grad}\left(\mathcal{P}_{r+1}\left(\mathbb{R}^{2}\right)\right)\,.
\end{equation}
Similar to $\mathbf{K}$, the operator $\widetilde{\mathbf{K}}:\mathcal{P}_{r}\left(\mathbb{R}^{2}\right)\rightarrow\mathcal{P}_{r+1}\left(T\mathbb{R}^{2}\right)$ gives the following vector when operated on $p\left(x_{1},x_{2}\right)\in\mathcal{P}_{r}\left(\mathbb{R}^{2}\right)$:
\begin{equation}\label{eqn:2.9}
\widetilde{\mathbf{K}}\left(p\right)=\left(px^{i}\right)\mathbf{e}_{i}\,.
\end{equation}
However, the operator $\tilde{\mathbf{k}}:\mathcal{P}_{r}\left(T\mathbb{R}^{2}\right)\rightarrow\mathcal{P}_{r+1}\left(\mathbb{R}^{2}\right)$ acting on the vector $\mathbf{v}\left(x_{1},x_{2}\right)\in\mathcal{P}_{r}\left(T\mathbb{R}^{2}\right)$ results in the following scalar:
\begin{equation}\label{eqn:2.10}
\tilde{\mathbf{k}}\left(\mathbf{v}\right)=\epsilon_{ij}v^{i}x^{j}\,.
\end{equation}
It is straightforward to verify the property $\tilde{\mathbf{k}}\circ\widetilde{\mathbf{K}}\left(p\right)=0$. Using \eqref{eqn:2.9}, it can be shown that the following decomposition also holds for the space $\mathcal{P}_{r}\left(T\mathbb{R}^{2}\right)$ \citep{JAH22}:
\begin{equation}\label{eqn:2.11}
\mathcal{P}_{r}\left(T\mathbb{R}^{2}\right)=
\widetilde{\mathbf{K}}\left(\mathcal{P}_{r-1}\left(\mathbb{R}^{2}\right)\right)\oplus
\mathbf{curl}\left(\mathcal{P}_{r+1}\left(\mathbb{R}^{2}\right)\right)\,,
\end{equation}
where the $\mathbf{curl}$ operator for $p\left(x_{1},x_{2}\right)\in\mathcal{P}_{r}\left(\mathbb{R}^{2}\right)$ is defined as $\mathbf{curl}\,p=\epsilon^{ji}\frac{\partial p}{\partial x_{i}}\mathbf{e}_{j}$. From \eqref{eqn:2.5} and noting the definitions of the operators $\mathbf{K}$ and $\mathbf{\widetilde{K}}$ in \eqref{eqn:2.6} and \eqref{eqn:2.9}, the following spaces can be defined as subspaces of $\mathcal{P}_{r}\left(T\mathbb{R}^{2}\right)$:
\begin{equation}\label{eqn:2.12}
\mathcal{P}^{-}_{r}\left(T\mathbb{R}^{2}\right)=
\mathcal{P}_{r-1}\left(T\mathbb{R}^{2}\right)+
\mathbf{K}\left(\mathcal{H}_{r-1}\left(\mathbb{R}^{2}\right)\right)\,,
\end{equation}
and
\begin{equation}\label{eqn:2.13}
\mathcal{P}^{\ominus}_{r}\left(T\mathbb{R}^{2}\right)=
\mathcal{P}_{r-1}\left(T\mathbb{R}^{2}\right)+
\mathbf{\widetilde{K}}\left(\mathcal{H}_{r-1}\left(\mathbb{R}^{2}\right)\right)\,.
\end{equation}
It is important to note that the subspaces in \eqref{eqn:2.12} and \eqref{eqn:2.13} are defined in the context of vector calculus; similar decompositions are derived in \citep{ANG17}. The subspaces in \eqref{eqn:2.12} and \eqref{eqn:2.13} and those given in \citep{ANG17} are based on similar decompositions in exterior calculus as proposed by \citet{ARN06}. Eqs.~\eqref{eqn:2.8}, \eqref{eqn:2.11}, \eqref{eqn:2.12} and \eqref{eqn:2.13} are used in the following subsections to construct the polynomial vector spaces that are used to approximate the displacement gradient and the stress tensor.

\subsection{Polynomial vector spaces on a simplex}
\label{sec:2.1}

Let us first define a space of shape functions for each simplex $\mathscr{T}$. This space can be chosen among $\mathcal{P}_{r}\left(T\mathscr{T}\right)$, $\mathcal{P}^{-}_{r}\left(T\mathscr{T}\right)$ and $\mathcal{P}^{\ominus}_{r}\left(T\mathscr{T}\right)$, which are obtained, respectively, by restricting the spaces $\mathcal{P}_{r}\left(T\mathbb{R}^{2}\right)$, $\mathcal{P}^{-}_{r}\left(T\mathbb{R}^{2}\right)$ and $\mathcal{P}^{\ominus}_{r}\left(T\mathbb{R}^{2}\right)$ to the simplex $\mathscr{T}$. It is also necessary to specify how the simplices are connected to obtain a global discretization. To achieve this goal, a basis dual to the space of shape functions is defined. Instead of associating the degrees of freedom directly to shape functions (coefficients of shape functions), they are defined on the dual space. The degrees of freedom in the dual space are associated to the subsimplices of the triangulation. When a subsimplex is shared by more than one simplex, the degree of freedom associated to that subsimplex must be single-valued. Relating  the degrees of freedom to subsimplices enables one to decompose the dual space of shape functions on $\mathscr{T}$ into a direct sum of subspaces that are indexed by the subsimplices of $\mathscr{T}$. This geometric decomposition of the dual space determines the interelement continuity rather than the specific choice of the degrees of freedom \citep{ARN06}.

\begin{figure}
\centering
\includegraphics[width=0.3\textwidth]{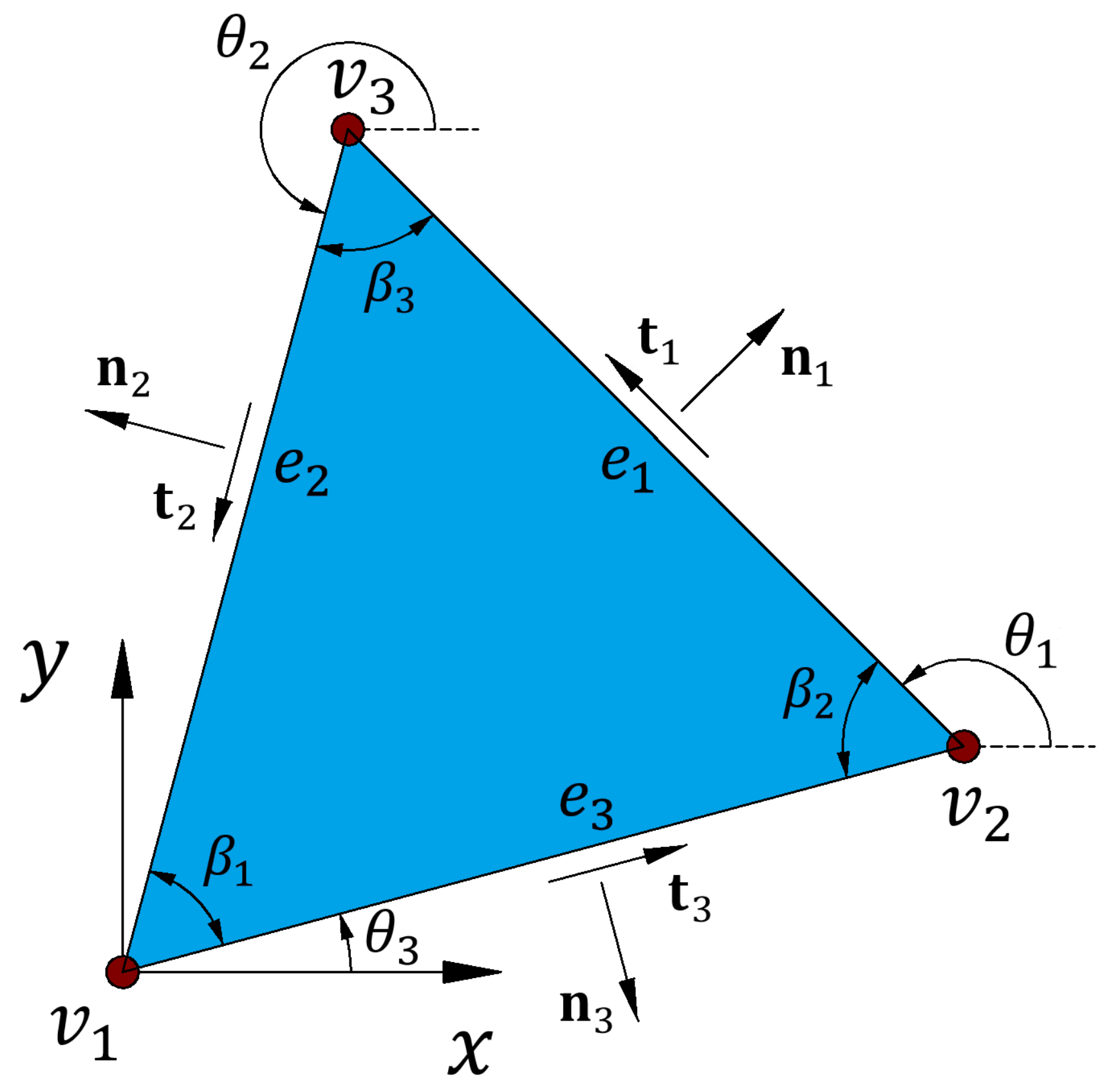}
\vskip 0.1in
\caption{Two-simplex $\mathscr{T}$ with vertices $v_{i}$ and edges $e_{i},i=1,2,3$. Tangent and normal vectors for each edge are denoted by $\mathbf{t}_{i}$ and $\mathbf{n}_{i}$, respectively. The orientation of edges is induced by the way the vertices are labeled.}
\label{fig:21}
\end{figure}

A two-simplex $\mathscr{T}$ is shown in Figure \ref{fig:21}. The vertices $v_{i},i=1,2,3$ are labeled in the counterclockwise direction and the origin of the coordinate system is chosen at $v_{1}$. The edge in front of the vertex $v_{i}$ is labeled as $e_{i}$. Thus, the edges are oriented according to the same scheme as the vertex labeling. For each vertex $v_{i}$, the internal angle between the edges incident to that vertex is denoted by $\beta_{i}$. The angle that the edge $e_{i}$ makes with the positive direction of the $x$ axis is identified with $\theta_{i}$. The vectors $\mathbf{t}_{i}$ and $\mathbf{n}_{i}$ are, respectively, the tangent and normal vectors for the edge $e_{i},i=1,2,3$. The following lemmas prove useful for defining the dual spaces, $\mathcal{P}_{r}\left(\mathscr{T}\right)^{*}$, $\mathcal{P}_{r}\left(T\mathscr{T}\right)^{*}$, $\mathcal{P}^{-}_{r}\left(T\mathscr{T}\right)^{*}$ and $\mathcal{P}^{\ominus}_{r}\left(T\mathscr{T}\right)^{*}$ on the simplex $\mathscr{T}$. They associate the dual degrees of freedom to the subsimplices of the simplex $\mathscr{T}$. It should be emphasized that these lemmas, defined in the context of vector calculus, are based on the geometrical decomposition of the spaces dual to the spaces of polynomial differential forms as discussed in detail in \citep{ARN06}. Furthermore, similar degrees of freedom can be obtained by following the approach of \citet{RAV77}, \citet{BRE85}, and \citet{NED80,NED86}.

\begin{lemma}\label{lem:1}
If the vertex $v$ and the edge $e$ belong to the set of all simplices $\Delta\left(\mathscr{T}\right)$ of the two-simplex $\mathscr{T}$, then for the polynomial $P\in\mathcal{P}_{r}\left(\mathbb{R}^{2}\right)$ restricted to the simplex $\mathscr{T}$, the following spaces can be defined as subspaces of $\mathcal{P}_{r}\left(\mathscr{T}\right)^{*}$:
\begin{equation}\label{eqn:2.14}
\begin{aligned}
&W\left(v\right)=
\left\{\phi\in\mathcal{P}_{r}\left(\mathscr{T}\right)^{*}\Big|\,
\phi\left(P\right)=P\left(x_{v},y_{v}\right)\right\}\,,\\
&W\left(e\right)=
\left\{\phi\in\mathcal{P}_{r}\left(\mathscr{T}\right)^{*}\Big|\,
\exists Q\in\mathcal{P}_{r-2}\left(e\right)
\;\text{such that}\;
\phi\left(P\right)=\int_{e}PQ\,ds\right\}\,,\\
&W\left(\mathscr{T}\right)=
\left\{\phi\in\mathcal{P}_{r}\left(\mathscr{T}\right)^{*}\Big|\,
\exists Q\in\mathcal{P}_{r-3}\left(\mathscr{T}\right)
\;\text{such that}\;
\phi\left(P\right)=\int_{\mathscr{T}}PQ\,dA\right\}\,.
\end{aligned}
\end{equation}
Moreover, the dual space $\mathcal{P}_{r}\left(\mathscr{T}\right)^{*}$ can be represented as the following direct sum:
\begin{equation}
\label{eqn:2.15}
\mathcal{P}_{r}\left(\mathscr{T}\right)^{*}=
\underset{v\in\Delta\left(\mathscr{T}\right)}{\bigoplus}W\left(v\right)\oplus
\underset{e\in\Delta\left(\mathscr{T}\right)}{\bigoplus}W\left(e\right)\oplus
W\left(\mathscr{T}\right)\,.
\end{equation}
\begin{proof}

A natural isomorphism exists between the dual spaces of the polynomial $0$-forms and the dual spaces of polynomials. Therefore, any decomposition of the polynomial $0$-forms also induces a corresponding decomposition for polynomials. The decomposition in \eqref{eqn:2.15} follows from the geometric decomposition of the polynomial $0$-forms as shown in \cite[Section 4.2]{ARN06}.
\end{proof}
\end{lemma}

\begin{example}\label{exm:1}
This example shows the application of Lemma \ref{lem:1}. The polynomial, $P\left(x,y\right)
=a_{1}
+b_{1}x+b_{2}y
+c_{1}x^{2}+c_{2}xy+c_{3}y^{2}
+d_{1}x^{3}+d_{2}x^{2}y+d_{3}xy^{2}+d_{4}y^{3}\in\mathcal{P}_{3}\left(\mathbb{R}^{2}\right)$
is the sum of homogeneous polynomials of degrees 0 to 3. Based on \eqref{eqn:2.14}, by restricting the polynomial $P\left(x,y\right)$ to the simplex $\mathscr{T}$, the elements $\phi_{i}\in\mathcal{P}_{3}\left(\mathscr{T}\right)^{*},i=1,\ldots,10$ can be defined as follows:
\begin{equation}\label{eqn:2.16}
\begin{aligned}
&\phi_{1}=P\left(x_{1},y_{1}\right)\,,&&
 \phi_{2}=P\left(x_{2},y_{2}\right)\,,&&
 \phi_{3}=P\left(x_{3},y_{3}\right)\,,&&
 \phi_{4}=\int_{e_{1}}PQ^{e_{1}}_{1}ds\,,\\
&\phi_{5}=\int_{e_{1}}PQ^{e_{1}}_{2}ds\,,&&
 \phi_{6}=\int_{e_{2}}PQ^{e_{2}}_{1}ds\,,&&
 \phi_{7}=\int_{e_{2}}PQ^{e_{2}}_{2}ds\,,&&
 \phi_{8}=\int_{e_{3}}PQ^{e_{3}}_{1}ds\,,\\
&\phi_{9}=\int_{e_{3}}PQ^{e_{3}}_{2}ds\,,&&
 \phi_{10}=\int_{\mathscr{T}}PQ^{\mathscr{T}}dA\,,&&
\end{aligned}
\end{equation}
where $\left(x_{i},y_{i}\right)$ are the coordinates of vertices $v_{i}$, $Q^{e_{i}}_{1},Q^{e_{i}}_{2}\in\mathcal{P}_{1}\left(e_{i}\right)$ and $Q^{\mathscr{T}}\in\mathcal{P}_{0}\left(\mathscr{T}\right)$.
\end{example}

\begin{lemma}\label{lem:2}
If the edge $e$ belongs to the set of all simplices $\Delta\left(\mathscr{T}\right)$ of the two-simplex $\mathscr{T}$ and $\mathbf{t}_{e}$ is the tangent vector associated to $e$, then for the vector $\mathbf{v}\in\mathcal{P}_{r}\left(T\mathbb{R}^{2}\right)$ restricted to the simplex $\mathscr{T}$, the following spaces can be defined as subspaces of $\mathcal{P}_{r}\left(T\mathscr{T}\right)^{*}$:
\begin{equation}\label{eqn:2.17}
\begin{aligned}
&W^{c}\left(e\right)=
\left\{\phi\in\mathcal{P}_{r}\left(T\mathscr{T}\right)^{*}\Big|\,
\exists Q\in\mathcal{P}_{r}\left(e\right)
\;\text{such that}\;
\phi\left(\mathbf{v}\right)=\int_{e}\left(\mathbf{v}\cdot\mathbf{t}_{e}\right)Q\,ds\right\}\,,\\
&W^{c}\left(\mathscr{T}\right)=
\left\{\phi\in\mathcal{P}_{r}\left(T\mathscr{T}\right)^{*}\Big|\,
\exists \mathbf{w}\in\mathcal{P}^{\ominus}_{r-1}\left(T\mathscr{T}\right)
\;\text{such that}\;
\phi\left(\mathbf{v}\right)=\int_{\mathscr{T}}\left(\mathbf{v}\cdot\mathbf{w}\right)\,dA\right\}\,.
\end{aligned}
\end{equation}
Moreover, the dual space $\mathcal{P}_{r}\left(T\mathscr{T}\right)^{*}$ can be expressed as the following direct sum:
\begin{equation}\label{eqn:2.18}
\mathcal{P}_{r}\left(T\mathscr{T}\right)^{*}=
\underset{e\in\Delta\left(\mathscr{T}\right)}{\bigoplus}W^{c}\left(e\right)\oplus
W^{c}\left(\mathscr{T}\right)\,.
\end{equation}
\end{lemma}

\begin{corollary}\label{cor:2.1}
If $\mathbf{n}_{e}$ is the normal vector associated to the edge $e$ in Lemma \ref{lem:2}, then for the vector $\mathbf{v}\in\mathcal{P}_{r}\left(T\mathbb{R}^{2}\right)$ restricted to the simplex $\mathscr{T}$, the following spaces are also subspaces of $\mathcal{P}_{r}\left(T\mathscr{T}\right)^{*}$:
\begin{equation}\label{eqn:2.19}
\begin{aligned}
&W^{d}\left(e\right)=
\left\{\phi\in\mathcal{P}_{r}\left(T\mathscr{T}\right)^{*}\Big|\,
\exists Q\in\mathcal{P}_{r}\left(e\right)
\;\text{such that}\;
\phi\left(\mathbf{v}\right)=\int_{e}\left(\mathbf{v}\cdot\mathbf{n}_{e}\right)Q\,ds\right\}\,,\\
&W^{d}\left(\mathscr{T}\right)=
\left\{\phi\in\mathcal{P}_{r}\left(T\mathscr{T}\right)^{*}\Big|\,
\exists \mathbf{w}\in\mathcal{P}^{-}_{r-1}\left(T\mathscr{T}\right)
\;\text{such that}\;
\phi\left(\mathbf{v}\right)=\int_{\mathscr{T}}\left(\mathbf{v}\cdot\mathbf{w}\right)\,dA\right\}\,.
\end{aligned}
\end{equation}
The dual space $\mathcal{P}_{r}\left(T\mathscr{T}\right)^{*}$ can be written as the following direct sum:
\begin{equation}\label{eqn:2.20}
\mathcal{P}_{r}\left(T\mathscr{T}\right)^{*}=
\underset{e\in\Delta\left(\mathscr{T}\right)}{\bigoplus}W^{d}\left(e\right)\oplus
W^{d}\left(\mathscr{T}\right)\,.
\end{equation}
\end{corollary}
\begin{proof}
In $\mathbb{R}^{2}$, $1$-forms or their Hodge-duals are isomorphic to vectors. In the same way, an isomorphism exists between the dual spaces of $1$-forms or their Hodge-duals, and the dual spaces of vectors. Therefore, any decomposition of dual spaces of $1$-forms induces a corresponding decomposition of dual spaces of vectors. The decomposition in \eqref{eqn:2.18} via subspaces in \eqref{eqn:2.17} follows from the geometric decomposition of the space dual to the space of polynomial $k$-forms as shown in \cite[Section 4.5]{ARN06}. The results in Corollary \ref{cor:2.1} are obtained from the aforementioned isomorphism and the observation that in $\mathbb{R}^{2}$ the Hodge-duals of $1$-forms are $1$-forms and, hence they both belong to the same space.
\end{proof}

\begin{example}\label{exm:2}
Let us consider the following polynomial vector
\begin{equation}\label{eqn:2.21}
\mathbf{v}=
\begin{bmatrix}
  a_{1}+b_{1}x+c_{1}y+d_{1}x^{2}+e_{1}xy+f_{1}y^{2} \\
  a_{2}+b_{2}x+c_{2}y+d_{2}x^{2}+e_{2}xy+f_{2}y^{2}
\end{bmatrix}\in\mathcal{P}_{2}\left(T\mathbb{R}^{2}\right)\,.
\end{equation}
By restricting $\mathbf{v}$ to the simplex $\mathscr{T}$ in Figure \ref{fig:21}, the elements $\phi_{i}\in\mathcal{P}_{2}\left(T\mathscr{T}\right)^{*},i=1,\ldots,12$ can be defined using \eqref{eqn:2.17} as follows:
\begin{equation}\label{eqn:2.22}
\begin{aligned}
&\phi_{1}=\int_{e_{1}}\left(\mathbf{v}\cdot\mathbf{t}_{1}\right)Q^{e_{1}}_{1}ds\,,&&
 \phi_{2}=\int_{e_{1}}\left(\mathbf{v}\cdot\mathbf{t}_{1}\right)Q^{e_{1}}_{2}ds\,,&&
 \phi_{3}=\int_{e_{1}}\left(\mathbf{v}\cdot\mathbf{t}_{1}\right)Q^{e_{1}}_{3}ds\,,&&
 \phi_{4}=\int_{e_{2}}\left(\mathbf{v}\cdot\mathbf{t}_{2}\right)Q^{e_{2}}_{1}ds\,,\\
&\phi_{5}=\int_{e_{2}}\left(\mathbf{v}\cdot\mathbf{t}_{2}\right)Q^{e_{2}}_{2}ds\,,&&
 \phi_{6}=\int_{e_{2}}\left(\mathbf{v}\cdot\mathbf{t}_{2}\right)Q^{e_{2}}_{3}ds\,,&&
 \phi_{7}=\int_{e_{3}}\left(\mathbf{v}\cdot\mathbf{t}_{3}\right)Q^{e_{3}}_{1}ds\,,&&
 \phi_{8}=\int_{e_{3}}\left(\mathbf{v}\cdot\mathbf{t}_{3}\right)Q^{e_{3}}_{2}ds\,,\\
&\phi_{9}=\int_{e_{3}}\left(\mathbf{v}\cdot\mathbf{t}_{3}\right)Q^{e_{3}}_{3}ds\,,&&
 \phi_{10}=\int_{\mathscr{T}}\left(\mathbf{v}\cdot\mathbf{w}^{\mathscr{T}}_{1}\right)dA\,,&&
 \phi_{11}=\int_{\mathscr{T}}\left(\mathbf{v}\cdot\mathbf{w}^{\mathscr{T}}_{2}\right)dA\,,&&
 \phi_{12}=\int_{\mathscr{T}}\left(\mathbf{v}\cdot\mathbf{w}^{\mathscr{T}}_{3}\right)dA\,.&&
\end{aligned}
\end{equation}
where $Q^{e_{i}}_{1},Q^{e_{i}}_{2},Q^{e_{i}}_{3}\in\mathcal{P}_{2}\left(e_{i}\right)$ and $\mathbf{w}^{\mathscr{T}}_{1},\mathbf{w}^{\mathscr{T}}_{2},\mathbf{w}^{\mathscr{T}}_{3}\in\mathcal{P}^{\ominus}_{1}\left(T\mathscr{T}\right)$. By inspecting \eqref{eqn:2.19}, we can derive the alternative elements $\phi_{i}\in\mathcal{P}_{2}\left(T\mathscr{T}\right)^{*}$ as:
\begin{equation}\label{eqn:2.23}\begin{aligned}
&\phi_{1}=\int_{e_{1}}\left(\mathbf{v}\cdot\mathbf{n}_{1}\right)Q^{e_{1}}_{1}ds\,,&&
 \phi_{2}=\int_{e_{1}}\left(\mathbf{v}\cdot\mathbf{n}_{1}\right)Q^{e_{1}}_{2}ds\,,&&
 \phi_{3}=\int_{e_{1}}\left(\mathbf{v}\cdot\mathbf{n}_{1}\right)Q^{e_{1}}_{3}ds\,,&&
 \phi_{4}=\int_{e_{2}}\left(\mathbf{v}\cdot\mathbf{n}_{2}\right)Q^{e_{2}}_{1}ds\,,\\
&\phi_{5}=\int_{e_{2}}\left(\mathbf{v}\cdot\mathbf{n}_{2}\right)Q^{e_{2}}_{2}ds\,,&&
 \phi_{6}=\int_{e_{2}}\left(\mathbf{v}\cdot\mathbf{n}_{2}\right)Q^{e_{2}}_{3}ds\,,&&
 \phi_{7}=\int_{e_{3}}\left(\mathbf{v}\cdot\mathbf{n}_{3}\right)Q^{e_{3}}_{1}ds\,,&&
 \phi_{8}=\int_{e_{3}}\left(\mathbf{v}\cdot\mathbf{n}_{3}\right)Q^{e_{3}}_{2}ds\,,\\
&\phi_{9}=\int_{e_{3}}\left(\mathbf{v}\cdot\mathbf{n}_{3}\right)Q^{e_{3}}_{3}ds\,,&&
 \phi_{10}=\int_{\mathscr{T}}\left(\mathbf{v}\cdot\mathbf{w}^{\mathscr{T}}_{1}\right)dA\,,&&
 \phi_{11}=\int_{\mathscr{T}}\left(\mathbf{v}\cdot\mathbf{w}^{\mathscr{T}}_{2}\right)dA\,,&&
 \phi_{12}=\int_{\mathscr{T}}\left(\mathbf{v}\cdot\mathbf{w}^{\mathscr{T}}_{3}\right)dA\,.
\end{aligned}
\end{equation}
Regarding the degrees of freedom $\phi_{10},\phi_{11},\phi_{12}$, it is noted that $\mathbf{w}^{\mathscr{T}}_{1},\mathbf{w}^{\mathscr{T}}_{2},\mathbf{w}^{\mathscr{T}}_{3}\in\mathcal{P}^{-}_{1}\left(T\mathscr{T}\right)$.
\end{example}

\begin{lemma}\label{lem:3}
If the edge $e$ belongs to the set of all simplices $\Delta\left(\mathscr{T}\right)$ of the two-simplex $\mathscr{T}$ and $\mathbf{t}_{e}$ is the tangent vector associated to $e$, then for the vector $\mathbf{v}\in\mathcal{P}^{-}_{r}\left(T\mathbb{R}^{2}\right)$ restricted to the simplex $\mathscr{T}$, the following spaces can be defined as subspaces of $\mathcal{P}^{-}_{r}\left(T\mathscr{T}\right)^{*}$:
\begin{equation}\label{eqn:2.24}
\begin{aligned}
&W^{c}\left(e\right)=
\left\{\phi\in\mathcal{P}^{-}_{r}\left(T\mathscr{T}\right)^{*}\Big|\,
\exists Q\in\mathcal{P}_{r-1}\left(e\right)
\;\text{such that}\;
\phi\left(\mathbf{v}\right)=\int_{e}\left(\mathbf{v}\cdot\mathbf{t}_{e}\right)Q\,ds\right\}\,,\\
&W^{c}\left(\mathscr{T}\right)=
\left\{\phi\in\mathcal{P}^{-}_{r}\left(T\mathscr{T}\right)^{*}\Big|\,
\exists\mathbf{w}\in\mathcal{P}_{r-2}\left(T\mathscr{T}\right)
\;\text{such that}\;
\phi\left(\mathbf{v}\right)=\int_{\mathscr{T}}\left(\mathbf{v}\cdot\mathbf{w}\right)\,dA\right\}\,.
\end{aligned}
\end{equation}
Moreover, the dual space $\mathcal{P}^{-}_{r}\left(T\mathscr{T}\right)^{*}$ can be written as the following direct sum:
\begin{equation}\label{eqn:2.25}
\mathcal{P}^{-}_{r}\left(T\mathscr{T}\right)^{*}=
\underset{e\in\Delta\left(\mathscr{T}\right)}{\bigoplus}W^{c}\left(e\right)\oplus
W^{c}\left(\mathscr{T}\right)\,.
\end{equation}
\end{lemma}

\begin{corollary}\label{cor:3.1}
If $\mathbf{n}_{e}$ is the normal vector associated to the edge $e$ in Lemma \ref{lem:3}, then for the vector $\mathbf{v}\in\mathcal{P}^{\ominus}_{r}\left(T\mathbb{R}^{2}\right)$ restricted to the simplex $\mathscr{T}$, the following spaces are also subspaces of $\mathcal{P}^{\ominus}_{r}\left(T\mathscr{T}\right)^{*}$:
\begin{equation}\label{eqn:2.26}
\begin{aligned}
&W^{d}\left(e\right)=
\left\{\phi\in\mathcal{P}^{\ominus}_{r}\left(T\mathscr{T}\right)^{*}\Big|\,
\exists Q\in\mathcal{P}_{r-1}\left(e\right)
\;\text{such that}\;
\phi\left(\mathbf{v}\right)=\int_{e}\left(\mathbf{v}\cdot\mathbf{n}_{e}\right)Q\,ds\right\}\,,\\
&W^{d}\left(\mathscr{T}\right)=
\left\{\phi\in\mathcal{P}^{\ominus}_{r}\left(T\mathscr{T}\right)^{*}\Big|\,
\exists \mathbf{w}\in\mathcal{P}_{r-2}\left(T\mathscr{T}\right)
\;\text{such that}\;
\phi\left(\mathbf{v}\right)=\int_{\mathscr{T}}\left(\mathbf{v}\cdot\mathbf{w}\right)\,dA\right\}\,.
\end{aligned}
\end{equation}
The dual space $\mathcal{P}^{\ominus}_{r}\left(T\mathscr{T}\right)^{*}$ can be written as the following direct sum:
\begin{equation}
\label{eqn:2.27}
\mathcal{P}^{\ominus}_{r}\left(T\mathscr{T}\right)^{*}=
\underset{e\in\Delta\left(\mathscr{T}\right)}{\bigoplus}W^{d}\left(e\right)\oplus
W^{d}\left(\mathscr{T}\right)\,.
\end{equation}
\end{corollary}
\begin{proof}
Similar to Lemma \ref{lem:2}, the proof follows from the isomorphism that exists between the dual spaces of $1$-forms or their Hodge-duals in $\mathbb{R}^{2}$, and the dual spaces of vectors. The geometric decomposition of the space $\mathcal{P}^{-}_{r}\Lambda^{k}\left(\mathscr{T}\right)^{*}$ dual to the space of reduced polynomial $k$-forms on the simplex $\mathscr{T}$ is discussed in \cite[Section 4.6]{ARN06}. The decomposition in \eqref{eqn:2.25} via the subspaces in \eqref{eqn:2.24} is obtained by restricting the dual space $\mathcal{P}^{-}_{r}\Lambda^{k}\left(\mathscr{T}\right)^{*}$ to $1$-forms and applying the aforementioned isomorphism. The result in Corollary \ref{cor:3.1} follows from the notion that in $\mathbb{R}^{2}$ the Hodge-duals of $1$-forms are themselves $1$-forms. Therefore, the proper application of Hodge-duals of $1$-forms in the geometric decomposition of the space $\mathcal{P}^{-}_{r}\Lambda^{1}\left(\mathscr{T}\right)^{*}$ leads to the desired decomposition in \eqref{eqn:2.27}.
\end{proof}

\begin{example}\label{exm:3}
From \eqref{eqn:2.12} and \eqref{eqn:2.6}, a typical polynomial vector in $\mathcal{P}^{-}_{2}\left(T\mathbb{R}^{2}\right)$ can be expressed as:
\begin{equation}\label{eqn:2.28}
\mathbf{v}_{1}=
\begin{bmatrix}
a_{1}+b_{1}x+c_{1}y-\left(d_{1}x+d_{2}y\right)y\\
a_{2}+b_{2}x+c_{2}y+\left(d_{1}x+d_{2}y\right)x
\end{bmatrix}\,.
\end{equation}
Restricting this vector to the simplex $\mathscr{T}$ in Figure \ref{fig:21}, the elements $\phi_{i}\in\mathcal{P}^{-}_{2}\left(T\mathscr{T}\right)^{*},i=1,\ldots,8$ can be defined using \eqref{eqn:2.24} as follows:
\begin{equation}\label{eqn:2.29}
\begin{aligned}
&\phi_{1}=\int_{e_{1}}\left(\mathbf{v}_{1}\cdot\mathbf{t}_{1}\right)Q^{e_{1}}_{1}ds\,,&&
 \phi_{2}=\int_{e_{1}}\left(\mathbf{v}_{1}\cdot\mathbf{t}_{1}\right)Q^{e_{1}}_{2}ds\,,&&
 \phi_{3}=\int_{e_{2}}\left(\mathbf{v}_{1}\cdot\mathbf{t}_{2}\right)Q^{e_{2}}_{1}ds\,,&&
 \phi_{4}=\int_{e_{2}}\left(\mathbf{v}_{1}\cdot\mathbf{t}_{2}\right)Q^{e_{2}}_{2}ds\,,\\
&\phi_{5}=\int_{e_{3}}\left(\mathbf{v}_{1}\cdot\mathbf{t}_{3}\right)Q^{e_{3}}_{1}ds\,,&&
 \phi_{6}=\int_{e_{3}}\left(\mathbf{v}_{1}\cdot\mathbf{t}_{3}\right)Q^{e_{3}}_{2}ds\,,&&
 \phi_{7}=\int_{\mathscr{T}}\left(\mathbf{v}_{1}\cdot\mathbf{w}^{\mathscr{T}}_{1}\right)dA\,,&&
 \phi_{8}=\int_{\mathscr{T}}\left(\mathbf{v}_{1}\cdot\mathbf{w}^{\mathscr{T}}_{2}\right)dA\,,
\end{aligned}
\end{equation}
where $Q^{e_{i}}_{1},Q^{e_{i}}_{2}\in\mathcal{P}_{1}\left(e_{i}\right)$ and $\mathbf{w}^{\mathscr{T}}_{1},\mathbf{w}^{\mathscr{T}}_{2}\in
\mathcal{P}_{0}\left(T\mathscr{T}\right)$. On the other hand, a polynomial vector in $\mathcal{P}^{\ominus}_{2}\left(T\mathbb{R}^{2}\right)$ can be defined using \eqref{eqn:2.13} and \eqref{eqn:2.9} as:
\begin{equation}\label{eqn:2.30}
\mathbf{v}_{2}=
\begin{bmatrix}
a_{1}+b_{1}x+c_{1}y+\left(d_{1}x+d_{2}y\right)x\\
a_{2}+b_{2}x+c_{2}y+\left(d_{1}x+d_{2}y\right)y
\end{bmatrix}\,.
\end{equation}
Restricting this vector to the simplex $\mathscr{T}$, the elements $\phi_{i}\in\mathcal{P}^{\ominus}_{2}\left(T\mathscr{T}\right)^{*},~i=1,\ldots,8$ can be defined using \eqref{eqn:2.26} as follows:
\begin{equation}\label{eqn:2.31}
\begin{aligned}
&\phi_{1}=\int_{e_{1}}\left(\mathbf{v}_{2}\cdot\mathbf{n}_{1}\right)Q^{e_{1}}_{1}ds\,,&&
 \phi_{2}=\int_{e_{1}}\left(\mathbf{v}_{2}\cdot\mathbf{n}_{1}\right)Q^{e_{1}}_{2}ds\,,&&
 \phi_{3}=\int_{e_{2}}\left(\mathbf{v}_{2}\cdot\mathbf{n}_{2}\right)Q^{e_{2}}_{1}ds\,,&&
 \phi_{4}=\int_{e_{2}}\left(\mathbf{v}_{2}\cdot\mathbf{n}_{2}\right)Q^{e_{2}}_{2}ds\,,\\
&\phi_{5}=\int_{e_{3}}\left(\mathbf{v}_{2}\cdot\mathbf{n}_{3}\right)Q^{e_{3}}_{1}ds\,,&&
 \phi_{6}=\int_{e_{3}}\left(\mathbf{v}_{2}\cdot\mathbf{n}_{3}\right)Q^{e_{3}}_{2}ds\,,&&
 \phi_{7}=\int_{\mathscr{T}}\left(\mathbf{v}_{2}\cdot\mathbf{w}^{\mathscr{T}}_{1}\right)dA\,,&&
 \phi_{8}=\int_{\mathscr{T}}\left(\mathbf{v}_{2}\cdot\mathbf{w}^{\mathscr{T}}_{2}\right)dA\,.
\end{aligned}
\end{equation}
Similar to \eqref{eqn:2.29}, $Q^{e_{i}}_{1},Q^{e_{i}}_{2}\in\mathcal{P}_{1}\left(e_{i}\right)$ and $\mathbf{w}^{\mathscr{T}}_{1},\mathbf{w}^{\mathscr{T}}_{2}\in\mathcal{P}_{0}\left(T\mathscr{T}\right)$.
\end{example}

\section{Second-Order Shape Functions}\label{sec:3}

In this section, the second-order $\mathcal{P}_{2}\left(T\mathscr{T}\right)$ and $\mathcal{P}^{-}_{2}\left(T\mathscr{T}\right)$ shape functions are discussed in detail. It is assumed that $\mathcal{B}$ is a bounded polygonal domain, which is discretized into a finite set of two-simplices $\mathcal{T}$. The union of these simplices is the closure of $\mathcal{B}$, while the intersection of any two simplices is either empty or a common edge of them. The set of all vertices and edges of a given simplex $\mathscr{T}\in\mathcal{T}$ are denoted by $\Delta_{0}\left(\mathscr{T}\right)$ and $\Delta_{1}\left(\mathscr{T}\right)$, respectively. The set of all vertices, edges and two-simplices of $\mathcal{T}$ are denoted by $\Delta_{0}\left(\mathcal{T}\right)$, $\Delta_{1}\left(\mathcal{T}\right)$ and $\Delta_{2}\left(\mathcal{T}\right)$, respectively. Furthermore, $\Delta^{\partial}_{0}\left(\mathcal{T}\right)$ and $\Delta^{\partial}_{1}\left(\mathcal{T}\right)$ represent the set of all vertices and edges of $\mathcal{T}$ that lie on the boundary of $\mathcal{B}$ \citep{ANG17,JAH22}. As a result, the set of all interior edges of $\mathcal{T}$ can be represented by $\Delta^{i}_{1}\left(\mathcal{T}\right)=\Delta_{1}\left(\mathcal{T}\right)\setminus
\Delta^{\partial}_{1}\left(\mathcal{T}\right)$. If a unique orientation is assigned to each edge of the triangulation $\mathcal{T}$, then the jump of the tangent component of the vector $\mathbf{v}\in\mathcal{P}_{r}\left(T\mathds{R}^{2}\right)$ across an edge $e\in\Delta^{i}_{1}\left(\mathcal{T}\right)$ that is common between the simplices $\mathscr{T}_{1},\mathscr{T}_{2}\in\Delta_{2}\left(\mathcal{T}\right)$ is defined as \citep{ANG17,JAH22}:
\begin{equation}\label{eqn:3.1}
{\llbracket\mathrm{t}\mathbf{v}\rrbracket}_{e}=
{\mathbf{v}|}_{\mathscr{T}_{1}}\cdot\mathbf{t}^{e}-
{\mathbf{v}|}_{\mathscr{T}_{2}}\cdot\mathbf{t}^{e},
\end{equation}
where $\mathbf{t}^{e}$ is the global unit vector tangent to the edge $e$. Similarly, denoting $\mathbf{n}^{e}$ as the global unit vector normal to this edge, the jump of the normal component of $\mathbf{v}$ across $e$ can be defined as:
\begin{equation}\label{eqn:3.2}
{\llbracket\mathrm{n}\mathbf{v}\rrbracket}_{e}=
{\mathbf{v}|}_{\mathscr{T}_{1}}\cdot\mathbf{n}^{e}-
{\mathbf{v}|}_{\mathscr{T}_{2}}\cdot\mathbf{n}^{e}.
\end{equation}
From these definitions, the spaces of polynomial vector fields with zero jumps for tangent and normal components across an internal edge $e\in\Delta^{i}_{1}\left(\mathcal{T}\right)$ are defined, respectively, as \citep{ANG17}:
\begin{equation}\label{eqn:3.3}
\begin{aligned}
\mathcal{P}^{c}_{r}\left(T\mathcal{T}\right)
&=\left\{\mathbf{v}\in\mathcal{P}_{r}\left(T\mathcal{T}\right)\Big|
{\llbracket\mathrm{t}\mathbf{v}\rrbracket}_{e}=0\,,\,\forall e\in\Delta^{i}_{1}\left(\mathcal{T}\right)\right\}\,,\\
\mathcal{P}^{d}_{r}\left(T\mathcal{T}\right)
&=\left\{\mathbf{v}\in\mathcal{P}_{r}\left(T\mathcal{T}\right)\Big|
{\llbracket\mathrm{n}\mathbf{v}\rrbracket}_{e}=0,\,\forall e\in\Delta^{i}_{1}\left(\mathcal{T}\right)\right\}\,.
\end{aligned}
\end{equation}
It is important to mention that similar definitions are given in  \citep{NED80,NED86}. Clearly, these spaces are subspaces of the Sobolev space $H\left(T\mathcal{B}\right)$. The spaces $\mathcal{P}^{c-}_{r}\left(T\mathcal{T}\right)$ and
$\mathcal{P}^{d-}_{r}\left(T\mathcal{T}\right)$ can be defined similarly. In the following, the local shape functions belonging to $\mathcal{P}^{c}_{2}\left(T\mathscr{T}\right)$, $\mathcal{P}^{d}_{2}\left(T\mathscr{T}\right)$, $\mathcal{P}^{c-}_{2}\left(T\mathscr{T}\right)$ and $\mathcal{P}^{d-}_{2}\left(T\mathscr{T}\right)$ are defined using the lemmas presented in \S\ref{sec:2.1} and the sets defined in \eqref{eqn:3.3}. These definitions are then generalized to obtain the global shape functions that belong to $\mathcal{P}^{c}_{2}\left(T\mathcal{T}\right)$, $\mathcal{P}^{d}_{2}\left(T\mathcal{T}\right)$, $\mathcal{P}^{c-}_{2}\left(T\mathcal{T}\right)$ and $\mathcal{P}^{d-}_{2}\left(T\mathcal{T}\right)$. The results of this section are used in Appendix \ref{app:A} to determine the explicit forms of the local shape functions for the simplicial element in the natural coordinate system. The covariant and contravariant Piola transformations are used to transform the shape functions from the natural coordinate system to the element in the physical space (reference configuration).

\subsection{$\mathcal{P}^{c}_{2}(T\mathcal{T})$ shape functions}
\label{sec:3.1}

\begin{figure}
\centering
\includegraphics[width=0.85\textwidth]{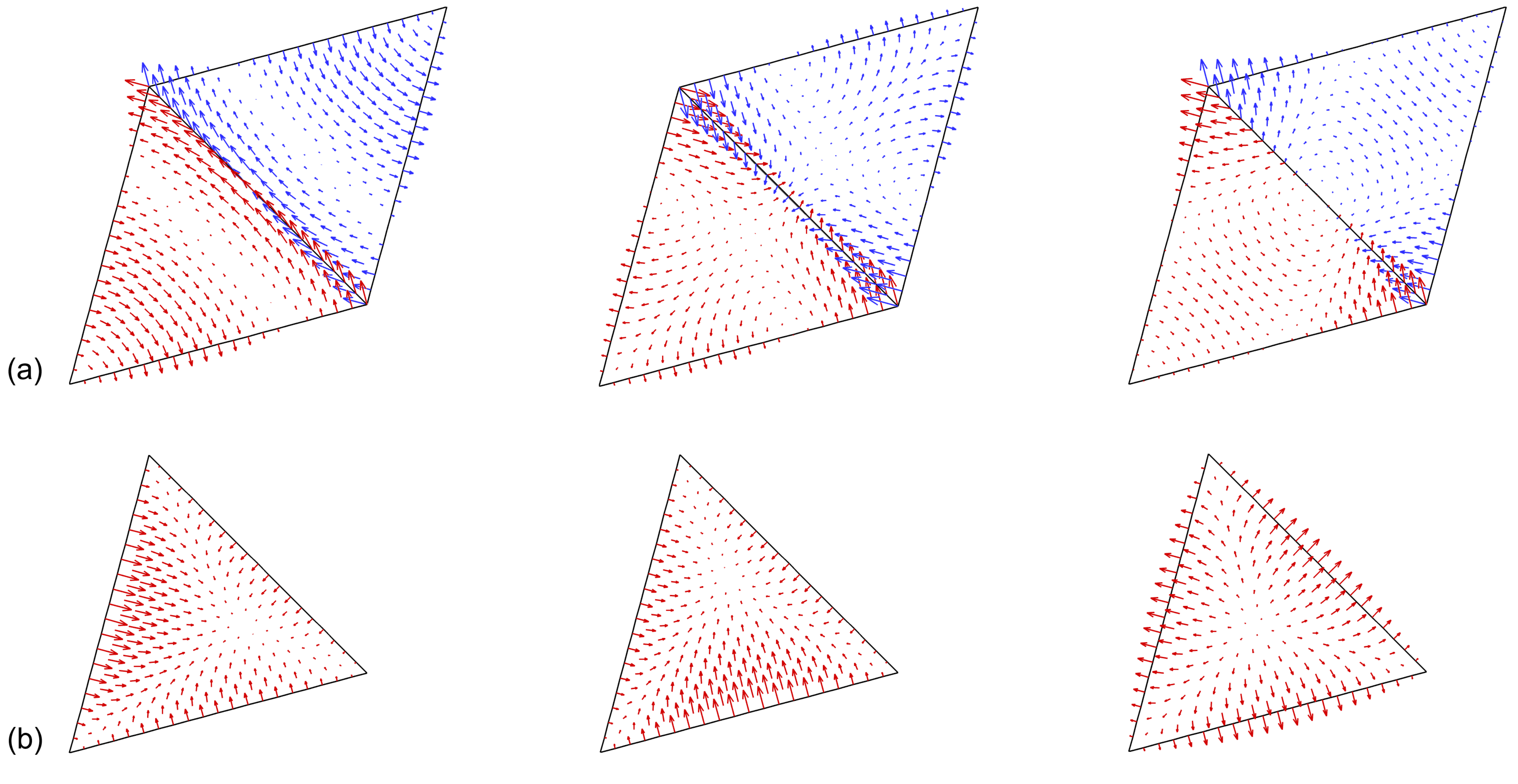}
\vskip 0.1in
\caption{(a) Global shape functions $\mathbf{V}^{e}_{1},\mathbf{V}^{e}_{2},\mathbf{V}^{e}_{3}\in\mathcal{P}^{c}_{2}\left(T\mathcal{T}\right)$ on the edge $e$ shared by two simplices $\mathscr{T}_{1}$ and $\mathscr{T}_{2}$. (b) Local shape functions $\mathbf{v}^{\mathscr{T}}_{1},\mathbf{v}^{\mathscr{T}}_{2},\mathbf{v}^{\mathscr{T}}_{3}\in\mathcal{P}^{c}_{2}\left(T\mathscr{T}\right)$ on the simplex $\mathscr{T}$. These shape functions are used for interpolating the displacement gradient.}
\label{fig:31}
\end{figure}

From Example \ref{exm:2}, the space of $\mathcal{P}^{c}_{2}(T\mathscr{T})$ shape functions can be defined based on the vector $\mathbf{v}$ in \eqref{eqn:2.21}. More specifically, the local shape functions on the edge $e_{i}$ (see Figure \ref{fig:21}) are defined as:
\begin{equation}\label{eqn:3.4}
\mathbf{v}^{\mathscr{T},e_{i}}_{j}=
\begin{bmatrix}
a^{j}_{1}+b^{j}_{1}x+c^{j}_{1}y+d^{j}_{1}x^2+e^{j}_{1}xy+f^{j}_{1}y^2 \\
a^{j}_{2}+b^{j}_{2}x+c^{j}_{2}y+d^{j}_{2}x^2+e^{j}_{2}xy+f^{j}_{2}y^2
\end{bmatrix},
\end{equation}
where $j=1,2,3$ refers to the index of shape function. The coefficients of $\mathbf{v}^{\mathscr{T},e_{i}}_{j}$ are determined by imposing certain conditions on the elements of the dual space $\mathcal{P}^{c}_{2}\left(T\mathscr{T}\right)^{*}$. Therefore, the elements of the dual space are called the degrees of freedom rather than the coefficients of $\mathbf{v}^{\mathscr{T},e_{i}}_{j}$. The degrees of freedom given in \eqref{eqn:2.22} are divided into two categories: the degrees of freedom $\phi^{\mathscr{T},e_{i}}_{j},~i,j=1,2,3$ on the edges of the simplex $\mathscr{T}$ and the degrees of freedom $\phi^{\mathscr{T}}_{k},~k=1,2,3$ on the simplex $\mathscr{T}$ itself. It is convenient to choose the base functions $Q^{e_{i}}_{j}\in\mathcal{P}_{2}\left(e_{i}\right)$ in \eqref{eqn:2.22} as $Q^{e_{i}}_{j}=s^{j-1}$. Consequently, the degrees of freedom $\phi_{m},~m=1,\ldots,12$ in that equation can be represented as:
\begin{equation}\label{eqn:3.5}
\begin{aligned}
    \phi^{\mathscr{T},e_{i}}_{j}\left(\mathbf{v}\right)&=\int_{e_{i}}\left(\mathbf{v}\cdot\mathbf{t}_{i}\right)s^{j-1}ds\,,&& i,j=1,2,3\,, \\
    \phi^{\mathscr{T}}_{k}\left(\mathbf{v}\right)&=\int_{\mathscr{T}}\left(\mathbf{v}\cdot\mathbf{w}^{\mathscr{T}}_{k}\right)dA\,,&& k=1,2,3,
\end{aligned}
\end{equation}
where the base vectors $\mathbf{w}^{\mathscr{T}}_{k}\in\mathcal{P}^{\ominus}_{1}\left(T\mathscr{T}\right)$ are defined as:
\begin{equation}\label{eqn:3.6}
\mathbf{w}^{\mathscr{T}}_{1}=
\begin{bmatrix}
1 \\
0
\end{bmatrix},\qquad
\mathbf{w}^{\mathscr{T}}_{2}=
\begin{bmatrix}
0 \\
1
\end{bmatrix},\qquad
\mathbf{w}^{\mathscr{T}}_{3}=
\begin{bmatrix}
x \\
y
\end{bmatrix}.
\end{equation}
The coefficients of the shape function $\mathbf{v}^{\mathscr{T},e_{i}}_{j}$ are obtained by requiring the following conditions to be satisfied:
\begin{equation}\label{eqn:3.7}
\phi^{\mathscr{T},e_{l}}_{m}\left(\mathbf{v}^{\mathscr{T},e_{i}}_{j}\right)=
\begin{cases}
1,&\text{if } l=i,~m=j, \\
0,&\text{otherwise},
\end{cases}
\quad\text{and}\qquad
\phi^{\mathscr{T}}_{n}\left(\mathbf{v}^{\mathscr{T},e_{i}}_{j}\right)=0.
\end{equation}
It is more efficient to consider the local shape functions $\bar{\mathbf{v}}^{\mathscr{T},e_{i}}_{j},~j=1,2,3$ instead of $\mathbf{v}^{\mathscr{T},e_{i}}_{j}$. These shape functions are defined as follows:
\begin{align}\label{eqn:3.8}
\bar{\mathbf{v}}^{\mathscr{T},e_{i}}_{1}=\mathbf{v}^{\mathscr{T},e_{i}}_{1}+\frac{l_{i}}{2}\mathbf{v}^{\mathscr{T},e_{i}}_{2}+\frac{l^{2}_{i}}{3}\mathbf{v}^{\mathscr{T},e_{i}}_{3}\,,\quad
\bar{\mathbf{v}}^{\mathscr{T},e_{i}}_{2}=-\frac{l_{i}}{6}\mathbf{v}^{\mathscr{T},e_{i}}_{2}-\frac{l^{2}_{i}}{6}\mathbf{v}^{\mathscr{T},e_{i}}_{3}\,,\quad
\bar{\mathbf{v}}^{\mathscr{T},e_{i}}_{3}=\frac{1}{3}\mathbf{v}^{\mathscr{T},e_{i}}_{1}+\frac{l_{i}}{6}\mathbf{v}^{\mathscr{T},e_{i}}_{2}+\frac{2l^{2}_{i}}{15}\mathbf{v}^{\mathscr{T},e_{i}}_{3},
\end{align}
where $l_{i}$ is the length of the edge $e_{i}$. It is straightforward to verify that the following properties hold for the shape functions $\bar{\mathbf{v}}^{\mathscr{T},e_{i}}_{j}$:
\begin{align}\label{eqn:3.9}
\left(\bar{\mathbf{v}}^{\mathscr{T},e_{i}}_{1}\cdot\mathbf{t}_{i}\right)\bigg|_{e_{i}}=\frac{1}{l_{i}}\,,\qquad
\left(\bar{\mathbf{v}}^{\mathscr{T},e_{i}}_{2}\cdot\mathbf{t}_{i}\right)\bigg|_{e_{i}}=\frac{1}{l^{2}_{i}}\left(l_{i}-2s\right)\,,\qquad
\left(\bar{\mathbf{v}}^{\mathscr{T},e_{i}}_{3}\cdot\mathbf{t}_{i}\right)\bigg|_{e_{i}}=\frac{1}{l^{3}_{i}}\left(l_{i}-2s\right)^{2}\,.
\end{align}
Eqs.~\eqref{eqn:3.5}, \eqref{eqn:3.7} and \eqref{eqn:3.8} have been used in \S\ref{sec:A.1} to derive the explicit forms of the shape functions $\bar{\mathbf{v}}^{\mathscr{T},e_{1}}_{j},~j=1,2,3$. The local shape functions $\bar{\mathbf{v}}^{\mathscr{T},e_{1}}_{j},j=1,2,3$ are shown in Figure \ref{fig:31}a using red arrows. It is important to note that the shape functions $\mathbf{v}^{\mathscr{T},e_{i}}_{j}$ and $\bar{\mathbf{v}}^{\mathscr{T},e_{i}}_{j}$ are defined on the simplex $\mathscr{T}$. It is, however, desired to have the global shape functions that satisfy the jump condition for the tangent component of $\bar{\mathbf{v}}^{\mathscr{T},e_{i}}_{j}$ over the edge $e$ that is shared by two simplices $\mathscr{T}_{1}$ and $\mathscr{T}_{2}$. In other words, it is necessary to have shape functions that satisfy \eqref{eqn:3.3}$_1$. To this end, the global shape functions $\mathbf{V}^{e}_{j},~j=1,2,3$ are defined as follows:
\begin{equation}\label{eqn:3.10}
\mathbf{V}^{e}_{j}|_{\mathscr{T}_{l}}=
\begin{cases}
c^{l,i}_{j}\bar{\mathbf{v}}^{\mathscr{T}_{l},e_{i}}_{j}\,,&
\text{if }\mathscr{T}_{l}=\mathscr{T}_{1},\mathscr{T}_{2}\,,\\
0,&\text{otherwise},
\end{cases}
\end{equation}
where $e_{i}\in\Delta_{1}\left(\mathscr{T}_{l}\right)$ is the same edge with the global identifier $e$ and,
\begin{equation}\label{eqn:3.11}
c^{l,i}_{j}=
\begin{cases}
\mathbf{t}^{e}\cdot\mathbf{t}^{l}_{i},&\text{if }j\text{ is odd},\\
1,&\text{otherwise}.
\end{cases}
\end{equation}
In \eqref{eqn:3.11}, $\mathbf{t}^{e}$ is the global unit tangent vector assigned to the edge $e$, while $\mathbf{t}^{l}_{i}$ represents the local unit vector tangent to the edge $e_{i}\in\Delta_{1}\left(\mathscr{T}_{l}\right)$ (see Figure~\ref{fig:21}). If the edge $e$ is shared by one simplex only, then $c^{l,i}_{j}$ is always equal to $1$. The global shape functions $\mathbf{V}^{e}_{j},~j=1,2,3$ on an edge common between two simplices $\mathscr{T}_{1}$ and $\mathscr{T}_{2}$ are shown in Figure~\ref{fig:31}a. We observe that the tangent component of the shape functions in two simplices is single-valued over the common edge.
On the other hand, the local shape functions on the simplex $\mathscr{T}$ are defined as:
\begin{equation}\label{eqn:3.12}
\mathbf{v}^{\mathscr{T}}_{k}=
\begin{bmatrix}
a^{k}_{1}+b^{k}_{1}x+c^{k}_{1}y+d^{k}_{1}x^2+e^{k}_{1}xy+f^{k}_{1}y^2 \\
a^{k}_{2}+b^{k}_{2}x+c^{k}_{2}y+d^{k}_{2}x^2+e^{k}_{2}xy+f^{k}_{2}y^2
\end{bmatrix}\,.
\end{equation}
The coefficients of these shape functions can be determined by imposing the following conditions:
\begin{equation}\label{eqn:3.13}
\phi^{\mathscr{T},e_{l}}_{m}\left(\mathbf{v}^{\mathscr{T}}_{k}\right)=0,\qquad
\phi^{\mathscr{T}}_{n}\left(\mathbf{v}^{\mathscr{T}}_{k}\right)=
\begin{cases}
1,&\text{if } n=k,\\
0,&\text{otherwise}.
\end{cases}
\end{equation}
The conditions in \eqref{eqn:3.13} combined with the degrees of freedom in \eqref{eqn:3.5} are used in \S\ref{sec:A.1} to derive the explicit forms of the coefficients for the local shape functions $\mathbf{v}^{\mathscr{T}}_{k},~k=1,2,3$ in $\mathcal{P}^{c}_{2}\left(T\mathscr{T}\right)$. These shape functions are shown in Figure \ref{fig:31}b. It is worthwhile to mention that the global shape functions $\mathbf{V}^{\mathscr{T}}_{k}$ are equal to the local shape functions $\mathbf{v}^{\mathscr{T}}_{k}$ for $k=1,2,3$ and, therefore, no further provisions are required to derive the global shape functions on the simplex $\mathscr{T}$. This feature follows from the notion that the shape functions $\mathbf{v}^{\mathscr{T}}_{k}$ are normal to all edges of the simplex.

\subsection{$\mathcal{P}^{d}_{2}(T\mathcal{T})$ shape functions}
\label{sec:3.2}

\begin{figure}
\centering
\includegraphics[width=0.85\textwidth]{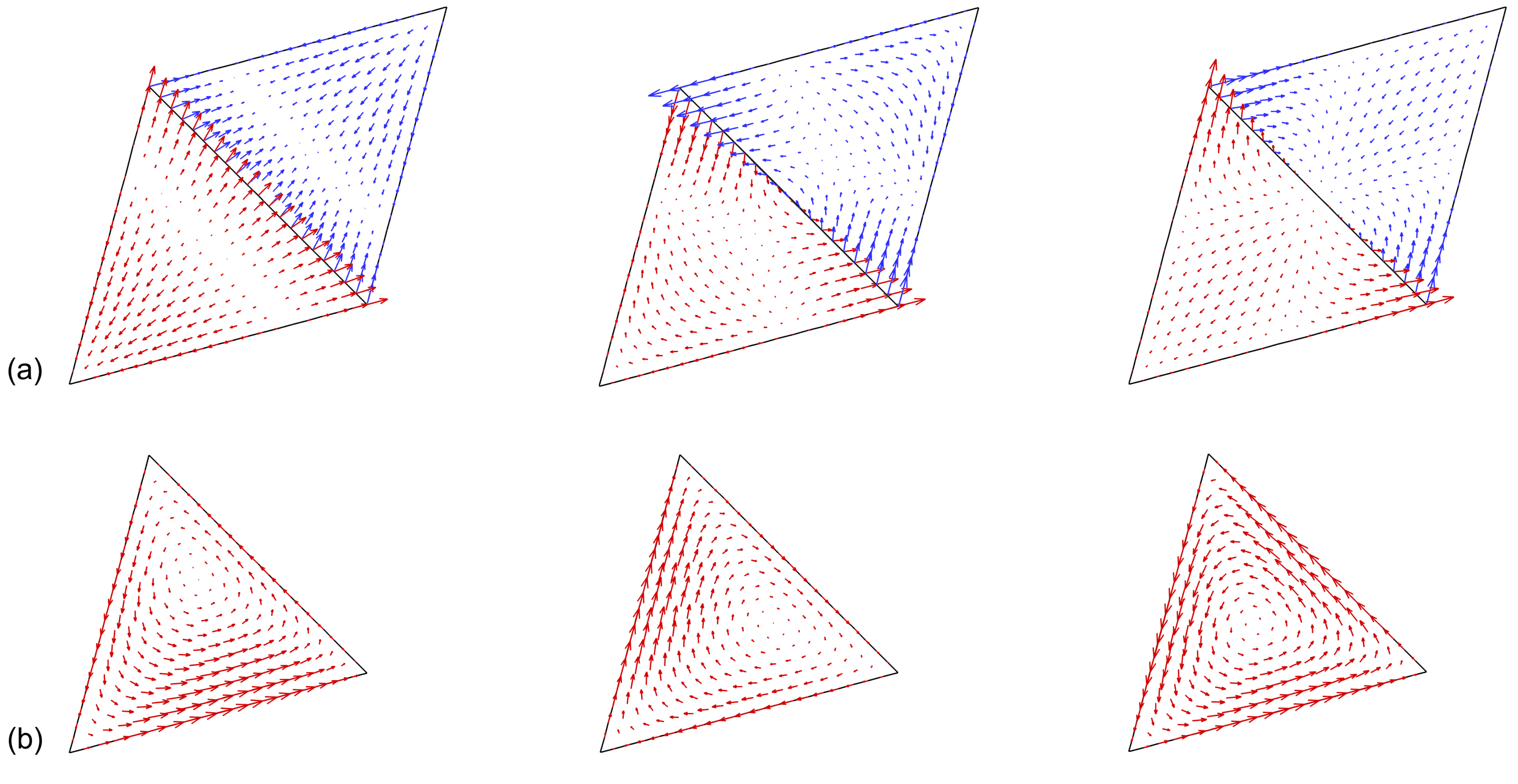}
\vskip 0.1in
\caption{(a) Global shape functions $\mathbf{V}^{e}_{1},\mathbf{V}^{e}_{2},\mathbf{V}^{e}_{3}\in\mathcal{P}^{d}_{2}\left(T\mathcal{T}\right)$ on the edge $e$ shared by two simplices $\mathscr{T}_{1}$ and $\mathscr{T}_{2}$. (b) Local shape functions $\mathbf{v}^{\mathscr{T}}_{1},\mathbf{v}^{\mathscr{T}}_{2},\mathbf{v}^{\mathscr{T}}_{3}\in\mathcal{P}^{d}_{2}\left(T\mathscr{T}\right)$ on the simplex $\mathscr{T}$. These shape functions are used for interpolating the stress tensor.}
\label{fig:32}
\end{figure}

The space of $\mathcal{P}^{d}_{2}\left(T\mathscr{T}\right)$ shape functions can be represented by the vectors used for $\mathcal{P}^{c}_{2}\left(T\mathscr{T}\right)$. The local shape functions $\mathbf{v}^{\mathscr{T},e_{i}}_{j},~j=1,2,3$ on the edges $e_{i},~i=1,2,3$ of the simplex $\mathscr{T}$ are defined through \eqref{eqn:3.4}. Similarly, the degrees of freedom in \eqref{eqn:2.23} are divided into the degrees of freedom on the edges of the simplex $\mathscr{T}$ and the degrees of freedom on the simplex itself. These degrees of freedom can be conveniently expressed as:
\begin{equation} \label{eqn:3.14}
\begin{aligned}
    \phi^{\mathscr{T},e_{i}}_{j}\left(\mathbf{v}\right)&=\int_{e_{i}}\left(\mathbf{v}\cdot\mathbf{n}_{i}\right)s^{j-1}ds\,,&& i,j=1,2,3\,,\\
    \phi^{\mathscr{T}}_{k}\left(\mathbf{v}\right)&=\int_{\mathscr{T}}\left(\mathbf{v}\cdot\mathbf{w}^{\mathscr{T}}_{k}\right)dA\,,&& k=1,2,3\,,
\end{aligned}
\end{equation}
where the base vectors $\mathbf{w}^{\mathscr{T}}_{k}\in\mathcal{P}^{-}_{1}\left(T\mathscr{T}\right)$ are here  defined as:
\begin{equation}\label{eqn:3.15}
\mathbf{w}^{\mathscr{T}}_{1}=
\begin{bmatrix}
1 \\
0
\end{bmatrix},\qquad
\mathbf{w}^{\mathscr{T}}_{2}=
\begin{bmatrix}
0 \\
1
\end{bmatrix},\qquad
\mathbf{w}^{\mathscr{T}}_{3}=
\begin{bmatrix}
-y \\
x
\end{bmatrix}.
\end{equation}
The coefficients of the shape functions in \eqref{eqn:3.4} are determined by enforcing the conditions in \eqref{eqn:3.7}. However, \eqref{eqn:3.14} should be applied in those conditions for evaluating the degrees of freedom. It is more efficient to consider the local shape functions $\bar{\mathbf{v}}^{\mathscr{T},e_{i}}_{j},~j=1,2,3$ as defined in \eqref{eqn:3.8}. It is observed that these shape functions satisfy the relations in \eqref{eqn:3.9} when the tangent vectors $\mathbf{t}_{i},~i=1,2,3$ in those relations are replaced with the normal vectors $\mathbf{n}_{i}$. The explicit forms of the local shape functions $\bar{\mathbf{v}}^{\mathscr{T},e_{i}}_{j},~j=1,2,3$ are derived in \S\ref{sec:A.2}. The variations of these shape functions are shown for the edge $e_{1}$ in Figure \ref{fig:32}a using red arrows. We observe that the shape functions are tangent to the edges $e_{2}$ and $e_{3}$. In order to have shape functions that satisfy the jump condition for the normal component of the shape function over the edge $e$ shared by two simplices $\mathscr{T}_{1}$ and $\mathscr{T}_{2}$, i.e., \eqref{eqn:3.3}$_2$, the global shape functions are defined as per \eqref{eqn:3.10} and \eqref{eqn:3.11}. Figure \ref{fig:32}a shows these shape functions whose normal component  is single-valued over the common edge of two simplices.

The local shape functions on the simplex $\mathscr{T}$ are defined using the vectors $\mathbf{v}^{\mathscr{T}}_{k},~k=1,2,3$ in \eqref{eqn:3.12}. The coefficients of these shape functions can be obtained by satisfying \eqref{eqn:3.13}, where again the degrees of freedom are computed from \eqref{eqn:3.14}. The explicit forms of the coefficients for the shape functions $\mathbf{v}^{\mathscr{T}}_{k},~k=1,2,3$ are given in \S\ref{sec:A.2}. These shape functions are shown in Figure \ref{fig:32}b. It is worth mentioning that the global shape functions on the simplex $\mathscr{T}$ are identical to the local shape functions and, therefore, it is not necessary to define them separately. This feature follows from the consideration that the shape functions on the simplex $\mathscr{T}$ are tangent to all edges.

\subsection{$\mathcal{P}^{c-}_{2}(T\mathcal{T})$ shape functions}
\label{sec:3.3}

\begin{figure}
\centering
\vskip 0.1in
\includegraphics[width=0.6\textwidth]{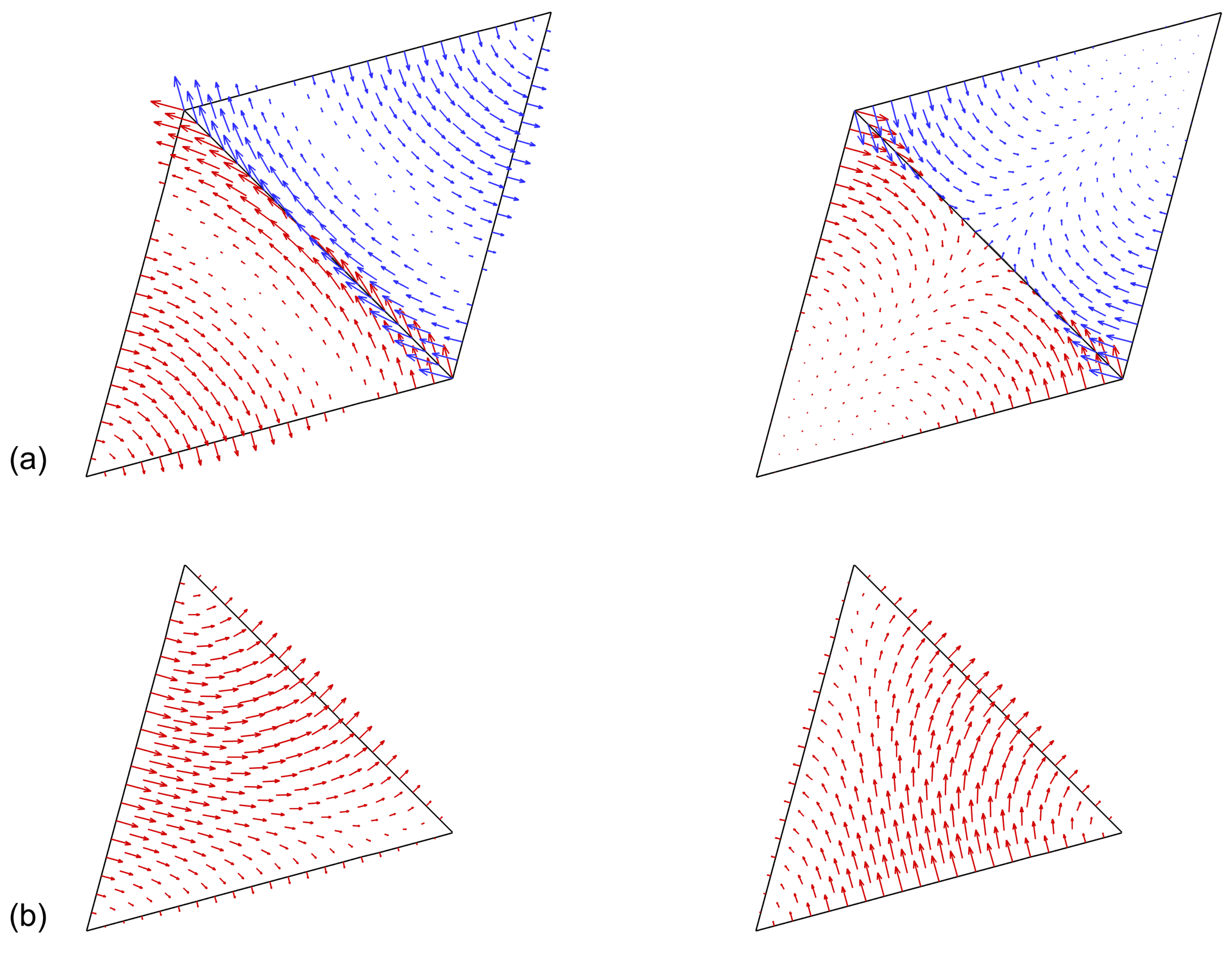}
\caption{(a) Global shape functions $\mathbf{V}^{e}_{1},\mathbf{V}^{e}_{2}\in\mathcal{P}^{c-}_{2}\left(T\mathcal{T}\right)$ on the edge $e$ shared by two simplices $\mathscr{T}_{1}$ and $\mathscr{T}_{2}$. (b) Local shape functions $\mathbf{v}^{\mathscr{T}}_{1},\mathbf{v}^{\mathscr{T}}_{2}\in\mathcal{P}^{c-}_{2}\left(T\mathscr{T}\right)$ on the simplex $\mathscr{T}$. These shape functions are used for interpolating the displacement gradient.}
\label{fig:33}
\end{figure}

The space of $\mathcal{P}^{c-}_{2}\left(T\mathscr{T}\right)$ shape functions can be defined from Example \ref{exm:3}. Eq.~\eqref{eqn:2.28} is used to define the local shape functions on the edge $e_{i}$ of the simplex $\mathscr{T}$ (see Figure~\ref{fig:21}) as:
\begin{equation}\label{eqn:3.16}
\mathbf{v}^{\mathscr{T},e_{i}}_{j}=
\begin{bmatrix}
a^{j}_{1}+b^{j}_{1}x+c^{j}_{1}y-(d^{j}_{1}x+d^{j}_{2}y)y\\
a^{j}_{2}+b^{j}_{2}x+c^{j}_{2}y+(d^{j}_{1}x+d^{j}_{2}y)x
\end{bmatrix},
\end{equation}
where $j=1,2$ refers to the index of the shape function. The degrees of freedom in \eqref{eqn:2.29} are divided into the degrees of freedom $\phi^{\mathscr{T},e_{i}}_{j},~j=1,2$ on the edges $e_{i},~i=1,2,3$ and the degrees of freedom $\phi^{\mathscr{T}}_{k},~k=1,2$ on the simplex $\mathscr{T}$. Choosing $s^{j-1}$ for the base functions $Q^{e_{i}}_{j}\in\mathcal{P}_{1}\left(e_{i}\right)$, the degrees of freedom $\phi_{m},~m=1,\ldots,8$ in \eqref{eqn:2.29} can be represented by the following equations:
\begin{equation} \label{eqn:3.17}
\begin{aligned}
    \phi^{\mathscr{T},e_{i}}_{j}\left(\mathbf{v}\right)&=\int_{e_{i}}\left(\mathbf{v}\cdot\mathbf{t}_{i}\right)s^{j-1}ds\,,&& i=1,2,3\,,~j=1,2,\\
    \phi^{\mathscr{T}}_{k}\left(\mathbf{v}\right)&=\int_{\mathscr{T}}\left(\mathbf{v}\cdot\mathbf{w}^{\mathscr{T}}_{k}\right)dA\,,&& k=1,2,
\end{aligned}
\end{equation}
where the base vectors $\mathbf{w}^{\mathscr{T}}_{k}$ are given by:
\begin{equation}\label{eqn:3.18}
\mathbf{w}^{\mathscr{T}}_{1}=
\begin{bmatrix}
1 \\
0
\end{bmatrix},\qquad
\mathbf{w}^{\mathscr{T}}_{2}=
\begin{bmatrix}
0 \\
1
\end{bmatrix}\,.
\end{equation}
Instead of the shape functions $\mathbf{v}^{\mathscr{T},e_{i}}_{j},~j=1,2$, it is possible to consider the alternative local shape functions $\bar{\mathbf{v}}^{\mathscr{T},e_{i}}_{j}$ defined as follows:
\begin{equation}\label{eqn:3.19}
\bar{\mathbf{v}}^{\mathscr{T},e_{i}}_{1}=\mathbf{v}^{\mathscr{T},e_{i}}_{1}+\frac{l_{i}}{2}\mathbf{v}^{\mathscr{T},e_{i}}_{2},\qquad
\bar{\mathbf{v}}^{\mathscr{T},e_{i}}_{2}=-\frac{l_{i}}{6}\mathbf{v}^{\mathscr{T},e_{i}}_{2}.
\end{equation}
It is easily verified that these shape functions satisfy the following relations:
\begin{equation}\label{eqn:3.20}
\left(\bar{\mathbf{v}}^{\mathscr{T},e_{i}}_{1}\cdot\mathbf{t}_{i}\right)\bigg|_{e_{i}}=\frac{1}{l_{i}},\qquad
\left(\bar{\mathbf{v}}^{\mathscr{T},e_{i}}_{2}\cdot\mathbf{t}_{i}\right)\bigg|_{e_{i}}=\frac{1}{l^{2}_{i}}\left(l_{i}-2s\right).
\end{equation}
Eqs. \eqref{eqn:3.7}, \eqref{eqn:3.17} and \eqref{eqn:3.19} are used in \S\ref{sec:A.3} to derive the explicit forms of the local shape functions $\bar{\mathbf{v}}^{\mathscr{T},e_{i}}_{j},~j=1,2$. The variations of these shape functions for the edge $e_{1}$ are shown in Figure \ref{fig:33}a using red arrows. In order to have shape functions that satisfy \eqref{eqn:3.3}$_1$ for the edges shared by two simplices, the global shape functions on the edge with global identifier $e$ are defined using \eqref{eqn:3.10}. The global shape functions $\mathbf{V}^{e}_{1}$ and $\mathbf{V}^{e}_{2}$ over the edge $e$ common between two simplices are shown in Figure~\ref{fig:33}a. It is clear that the tangent component of shape functions is single-valued over the edge.
On the other hand, the local shape functions on the simplex $\mathscr{T}$ are defined using the vector $\mathbf{v}^{\mathscr{T}}_{k},~k=1,2$ given bellow:
\begin{equation}\label{eqn:3.21}
\mathbf{v}^{\mathscr{T}}_{k}=
\begin{bmatrix}
a^{k}_{1}+b^{k}_{1}x+c^{k}_{1}y-\left(d^{k}_{1}x+d^{k}_{2}y\right)y\\
a^{k}_{2}+b^{k}_{2}x+c^{k}_{2}y+\left(d^{k}_{1}x+d^{k}_{2}y\right)x
\end{bmatrix}.
\end{equation}
The coefficients of these shape functions can be obtained by satisfying the conditions in \eqref{eqn:3.13}. The conditions in \eqref{eqn:3.13} combined with the degrees of freedom in \eqref{eqn:3.17} are used in \S\ref{sec:A.3} to derive the explicit forms of the coefficients for the local shape functions $\mathbf{v}^{\mathscr{T}}_{k},~k=1,2$, which are shown in Figure \ref{fig:33}b. It should be mentioned that the global shape functions on the simplex $\mathscr{T}$ are identical to the local shape functions $\mathbf{v}^{\mathscr{T}}_{k}$.

\subsection{$\mathcal{P}^{d-}_{2}(T\mathcal{T})$ shape functions}
\label{sec:3.4}

\begin{figure}
\centering
\includegraphics[width=0.6\textwidth]{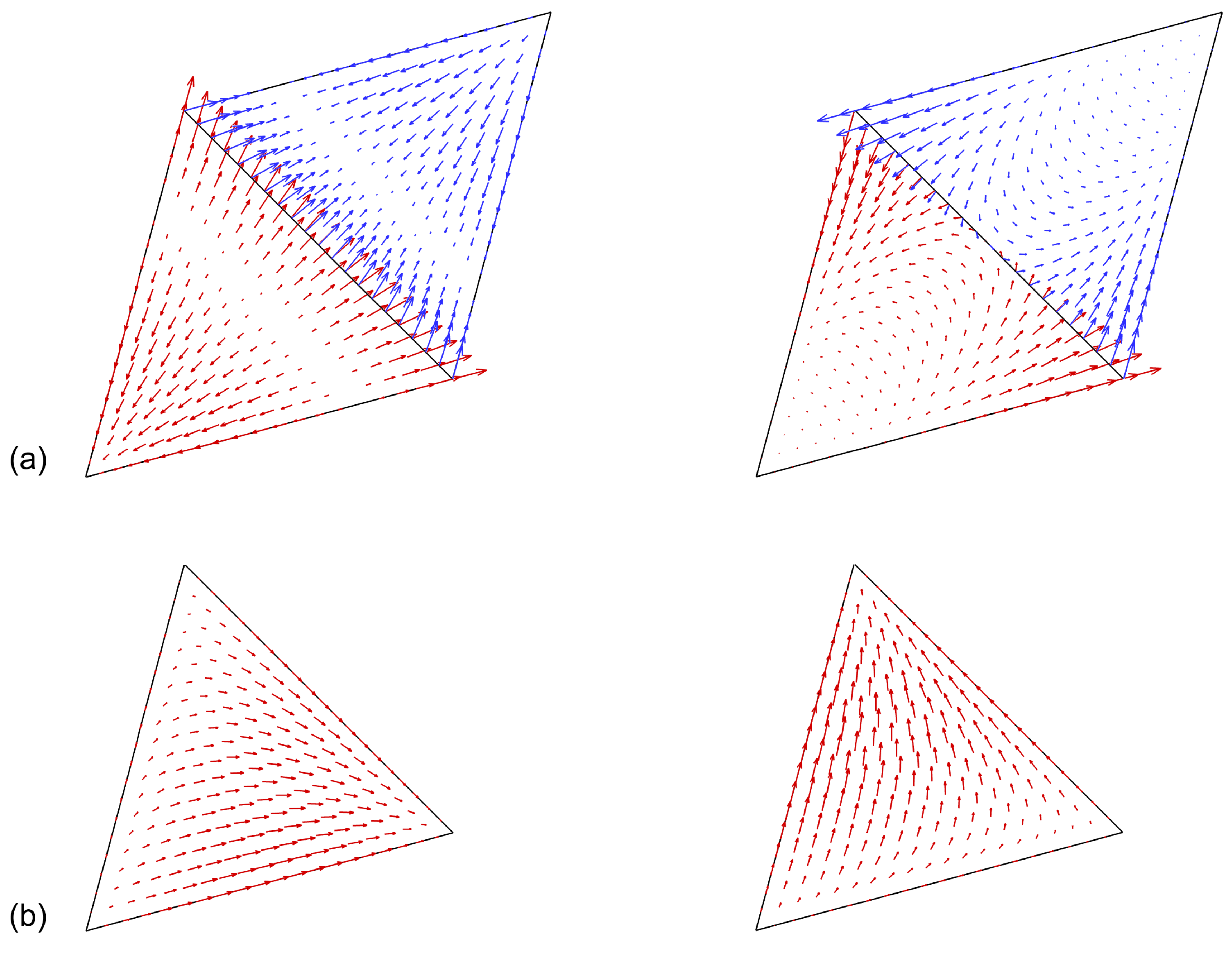}
\vskip 0.1in
\caption{(a) Global shape functions $\mathbf{V}^{e}_{1},\mathbf{V}^{e}_{2}\in\mathcal{P}^{d-}_{2}\left(T\mathcal{T}\right)$ on the edge $e$ shared by two simplices $\mathscr{T}_{1}$ and $\mathscr{T}_{2}$. (b) Local shape functions $\mathbf{v}^{\mathscr{T}}_{1},\mathbf{v}^{\mathscr{T}}_{2}\in\mathcal{P}^{d-}_{2}\left(T\mathscr{T}\right)$ on the simplex $\mathscr{T}$. These shape functions are used for interpolating the stress tensor.}
\label{fig:34}
\end{figure}

Example \ref{exm:3} can be used to define the space of $\mathcal{P}^{d-}_{2}\left(T\mathscr{T}\right)$ shape functions. From \eqref{eqn:2.30}, the local shape functions on the edge $e_{i}$ are expressed in the following form:
\begin{equation}\label{eqn:3.22}
\mathbf{v}^{\mathscr{T},e_{i}}_{j}=
\begin{bmatrix}
a^{j}_{1}+b^{j}_{1}x+c^{j}_{1}y+(d^{j}_{1}x+d^{j}_{2}y)x\\
a^{j}_{2}+b^{j}_{2}x+c^{j}_{2}y+(d^{j}_{1}x+d^{j}_{2}y)y
\end{bmatrix}\,.
\end{equation}
The degrees of freedom in \eqref{eqn:2.31}, which can be divided into the degrees of freedom on the edges $e_{i},~i=1,2,3$ of the simplex $\mathscr{T}$ and the degrees of freedom on the simplex itself, are expressed in the following form:
\begin{equation} \label{eqn:3.23}
\begin{aligned}
    \phi^{\mathscr{T},e_{i}}_{j}\left(\mathbf{v}\right)&=\int_{e_{i}}\left(\mathbf{v}\cdot\mathbf{n}_{i}\right)s^{j-1}ds,\;i=1,2,3\,,&& j=1,2\,, \\
    \phi^{\mathscr{T}}_{k}\left(\mathbf{v}\right)&=\int_{\mathscr{T}}\left(\mathbf{v}\cdot\mathbf{w}^{\mathscr{T}}_{k}\right)dA\,,&& k=1,2\,,
\end{aligned}
\end{equation}
where the vectors $\mathbf{w}^{\mathscr{T}}_{k},~k=1,2$ are given in \eqref{eqn:3.18}. It is, however, more efficient to consider the local shape functions $\bar{\mathbf{v}}^{\mathscr{T},e_{i}}_{j}$ as defined in \eqref{eqn:3.19}. It is straightforward to verify that these shape functions satisfy the relations in \eqref{eqn:3.20} when the tangent vectors $\mathbf{t}_{i}$ are replaced with the normal vectors $\mathbf{n}_{i}$. Eqs.~\eqref{eqn:3.7}, \eqref{eqn:3.19} and \eqref{eqn:3.23} have been utilized in \S\ref{sec:A.4} to derive the explicit forms of the shape functions $\bar{\mathbf{v}}^{\mathscr{T},e_{i}}_{j}$. The variations of these shape functions are shown for the edge $e_{1}$ in Figure \ref{fig:34}a using red arrows. It is desired to have global shape functions that satisfy \eqref{eqn:3.3}$_2$ concerning the jump for the normal component of the shape function over an edge common between two simplices. To achieve this goal, the global shape functions are defined from \eqref{eqn:3.10} and \eqref{eqn:3.11}. The shape functions so defined are shown for the edge $e$ shared by two simplices $\mathscr{T}_{1}$ and $\mathscr{T}_{2}$ in Figure \ref{fig:34}a. It is clear that the normal component of the shape functions is single-valued over the edge $e$. The local shape functions on the simplex $\mathscr{T}$ are represented by the vector $\mathbf{v}^{\mathscr{T}}_{k},k=1,2$ defined as:
\begin{equation}\label{eqn:3.24}
\mathbf{v}^{\mathscr{T}}_{k}=
\begin{bmatrix}
a^{k}_{1}+b^{k}_{1}x+c^{k}_{1}y+\left(d^{k}_{1}x+d^{k}_{2}y\right)x\\
a^{k}_{2}+b^{k}_{2}x+c^{k}_{2}y+\left(d^{k}_{1}x+d^{k}_{2}y\right)y
\end{bmatrix}\,.
\end{equation}
The coefficients of these shape functions can be determined by enforcing \eqref{eqn:3.13} and evaluating the degrees of freedom based on \eqref{eqn:3.23}. The explicit forms of the coefficients are given in \S\ref{sec:A.4}. Figure \ref{fig:34}b shows the shape functions $\mathbf{v}^{\mathscr{T}}_{k},~k=1,2$. It is to be noted that the global shape functions on the simplex $\mathscr{T}$ are identical to the local shape functions. Therefore, it is not necessary to compute them separately.

\section{Second-Order Compatible-Strain Mixed Finite Elements}
\label{sec:4}

This section discusses the finite element formulation in detail. Section~\ref{sec:4.1} outlines the variational formulation of the finite element method. The linearization and discretization of the governing equations are discussed, including the triangulation of the domain and the selection of suitable spaces for various shape functions. Section~\ref{sec:4.2} discusses the implementation of the mixed finite element method. Key topics, such as the minimum rank of submatrices required to obtain an invertible stiffness matrix, are also discussed. Furthermore, it is explained how the adoption of pseudo-nodes in the middle of edges facilitates enforcing the continuity constraints across element boundaries. These pseudo-nodes correspond to the displacement gradient and stress degrees of freedom associated with the edges to which they are connected.


\subsection{Variational formulation}
\label{sec:4.1}

An elastic body occupies a subset $\mathcal{B}\subset\mathbb{R}^2$ in the reference configuration and its surface where the traction forces are applied is denoted by $\partial_{\tau}\mathcal{B}$. The material coordinates $X^{i},~i=1,2$ refer to the position of material points in the reference configuration and the spatial coordinates $x^{i},~i=1,2$ indicate the position of the same material points in the spatial configuration. The displacement vector $\bm{\varphi}=\mathbf{x}-\mathbf{X}$, the displacement gradient $\mathbf{H}$, and the first Piola-Kirchhoff stress tensor $\mathbf{P}$ are considered as the independent fields in the Hu-Washizu functional defined as:
\begin{align}\label{eqn:4.1}
\Pi\left(\bm{\varphi},\mathbf{H},\mathbf{P}\right)=
\int_{\mathcal{B}}\widehat{W}\left(\mathbf{C}\right)dV+
\int_{\mathcal{B}}\bm{\tau}:\left(\nabla\bm{\varphi}-\mathbf{h}\right)dV
-\int_{\mathcal{B}}\rho_{0}\left(\mathbf{b}\cdot\bm{\varphi}\right)dV-
\int_{\partial_{\tau}\mathcal{B}}\mathbf{t}\cdot\bm{\varphi}\,dA,
\end{align}
where $\widehat{W}$ is the stored energy function, $\mathbf{C}=\mathbf{F}^{\mathsf{T}}\mathbf{F}$ is the right Cauchy-Green tensor, $\mathbf{F}=\mathbf{I}+\mathbf{H}$ is the deformation gradient, $\mathbf{h}=\mathbf{H}\mathbf{F}^{-1}$ is the push forward of the displacement gradient to the current configuration, $\bm{\tau}=\mathbf{P}\mathbf{F}^{\mathsf{T}}$ is the Kirchhoff stress tensor, $\rho_{0}$ is the material mass density, $\mathbf{b}$ is the body force and $\mathbf{t}$ is the traction force on a portion of the boundary where the traction is specified. The material and spatial derivative operators are, respectively, defined as $\nabla_{0}=\frac{\partial}{\partial\mathbf{X}}$ and $\nabla=\left(\frac{\partial}{\partial\mathbf{X}}\right)\mathbf{F}^{-1}$. It is noted that the Kirchhoff stress $\bm{\tau}$ plays the role of a Lagrange multiplier that enforces the constraint $\nabla\bm{\varphi}=\mathbf{h}$, which is equivalent to $\nabla_{0}\bm{\varphi}=\mathbf{H}$.

To derive the governing equations, the first variation of the functional \eqref{eqn:4.1} is set to zero. Let us denote the independent variations of $\bm{\varphi}$, $\mathbf{H}$ and $\mathbf{P}$ by $\delta\bm{\varphi}$, $\delta\mathbf{H}$ and $\delta\mathbf{P}$, respectively. The corresponding variations of the functional are denoted by $(\delta_{\bm{\varphi}}\Pi,\delta_{\mathbf{H}}\Pi,\delta_{\mathbf{P}}\Pi)=(D\Pi\cdot\delta\bm{\varphi},D\Pi\cdot\delta\mathbf{H},D\Pi\cdot\delta\mathbf{P})$, where $D\Pi$ is the derivative of the functional. These three variations are calculated as:

\begin{equation}\label{eqn:4.2}
D\Pi\cdot\delta\bm{\varphi}=
\int_{\mathcal{B}}\bm{\tau}:\nabla\delta\bm{\varphi}\,dV-
\int_{\mathcal{B}}\rho_{0}\left(\mathbf{b}\cdot\delta\bm{\varphi}\right)dV-
\int_{\partial_{\tau}\mathcal{B}}\mathbf{t}\cdot\delta\bm{\varphi}\,dA=0\,,
\end{equation}
\begin{equation}
\label{eqn:4.3}
D\Pi\cdot\delta\mathbf{H}=
\int_{\mathcal{B}}\bm{\widehat{\tau}}:\delta\mathbf{h}^{s}\,dV-
\int_{\mathcal{B}}\bm{\tau}:\delta\mathbf{h}\,dV=0,
\end{equation}
and
\begin{equation}\label{eqn:4.4}
D\Pi\cdot\delta\mathbf{P}=
\int_{\mathcal{B}}\delta\bm{\tau}:\nabla\bm{\varphi}\,dV-
\int_{\mathcal{B}}\delta\bm{\tau}:\mathbf{h}\,dV=0,
\end{equation}
where $\bm{\widehat{\tau}}=2\mathbf{F}\frac{\partial\widehat{W}}{\partial\mathbf{C}}\mathbf{F}^{\mathsf{T}}$
is the Kirchhoff stress tensor obtained from the stored energy function $\widehat{W}$ (compared with $\bm{\tau}$ that is obtained from $\mathbf{P}$) and $\delta\mathbf{h}^s=\frac{1}{2}\left(\delta\mathbf{h}+\delta\mathbf{h}^{\mathsf{T}}\right)$. Eqs.~\eqref{eqn:4.2}-\eqref{eqn:4.4} are nonlinear, and hence an iterative procedure must be employed to solve them. It is convenient to assign the residuals $R_{1}\left(\delta\bm{\varphi}\right)=\delta_{\delta\bm{\varphi}}\,\Pi$, $R_{2}\left(\delta\mathbf{H}\right)=\delta_{\delta\mathbf{H}}\,\Pi$ and $R_{3}\left(\delta\mathbf{P}\right)=\delta_{\delta\mathbf{P}}\,\Pi$. Then, considering $\Delta\bm{\varphi}$, $\Delta\mathbf{H}$ and $\Delta\mathbf{P}$ as another set of independent variations in $\bm{\varphi}$, $\mathbf{H}$ and $\mathbf{P}$, Eqs.~\eqref{eqn:4.2}-\eqref{eqn:4.4} can be linearized as follows:
\begin{equation}\label{eqn:4.5}
R_{1}\left(\delta\bm{\varphi}\right)+\delta R_{1}\left(\delta\bm{\varphi}\right)=0,\,\qquad
R_{2}\left(\delta\mathbf{H}\right)+\delta R_{2}\left(\delta\mathbf{H}\right)=0,\,\qquad
R_{3}\left(\delta\mathbf{P}\right)+\delta R_{3}\left(\delta\mathbf{P}\right)=0.
\end{equation}
It is easily verified that,
\begin{equation} \label{eqn:4.6}
\begin{aligned}
\delta R_{1}\left(\delta\bm{\varphi}\right)&=
DR_{1}\left(\delta\bm{\varphi}\right)\cdot\Delta\mathbf{P}
=\int_{\mathcal{B}}\Delta\bm{\tau}:\nabla\delta\bm{\varphi}\,dV \\
\delta R_{2}\left(\delta\mathbf{H}\right)&=
DR_{2}\left(\delta\mathbf{H}\right)\cdot\Delta\mathbf{H}+DR_{2}\left(\delta\mathbf{H}\right)\cdot\Delta\mathbf{P} \\
&=\int_{\mathcal{B}}\delta\mathbf{h}^{s}:\left(\boldsymbol{\mathbb{C}}:\Delta\mathbf{h}^{s}\right)dV+
\int_{\mathcal{B}}\left(\Delta\mathbf{h}\,\bm{\widehat{\tau}}\right):\delta\mathbf{h}\,dV
-\int_{\mathcal{B}}\Delta\bm{\tau}:\delta\mathbf{h}\,dV \\
\delta R_{3}\left(\delta\mathbf{P}\right)&=
DR_{3}\left(\delta\mathbf{P}\right)\cdot\Delta\bm{\varphi}+DR_{3}\left(\delta\mathbf{P}\right)\cdot\Delta\mathbf{H}
=\int_{\mathcal{B}}\delta\bm{\tau}:\nabla\Delta\bm{\varphi}\,dV-
\int_{\mathcal{B}}\delta\bm{\tau}:\Delta\mathbf{h}\,dV.
\end{aligned}
\end{equation}
In the above equations, $\Delta\bm{\tau}=\Delta\mathbf{P}\,\mathbf{F}^{\mathsf{T}}$, $\boldsymbol{\mathbb{C}}$ is the spatial fourth-order elasticity tensor, $\Delta\mathbf{h}=\Delta\mathbf{H}\,\mathbf{F}^{-1}$ and $\Delta\mathbf{h}^{s}$ is the symmetric part of $\Delta\mathbf{h}$. In light of \eqref{eqn:4.2}-\eqref{eqn:4.4} and \eqref{eqn:4.6}, the linearized equations in \eqref{eqn:4.5} can be cast into the following set of equations that should be solved at each iteration:
\begin{align}\label{eqn:4.7}
\int_{\mathcal{B}}\Delta\bm{\tau}:\nabla\delta\bm{\varphi}\,dV&=
\int_{\mathcal{B}}\rho_{0}\left(\mathbf{b}\cdot\delta\bm{\varphi}\right)dV+
\int_{\partial_{\tau}\mathcal{B}}\mathbf{t}\cdot\delta\bm{\varphi}\,dA-
\int_{\mathcal{B}}\bm{\tau}:\nabla\delta\bm{\varphi}\,dV,\nonumber\\
\int_{\mathcal{B}}\delta\mathbf{h}^{s}:\left(\boldsymbol{\mathbb{C}}:\Delta\mathbf{h}^{s}\right)dV&+
\int_{\mathcal{B}}\left(\Delta\mathbf{h}\,\bm{\widehat{\tau}}\right):\delta\mathbf{h}\,dV-
\int_{\mathcal{B}}\Delta\bm{\tau}:\delta\mathbf{h}\,dV=\nonumber\\
&-\left(\int_{\mathcal{B}}\bm{\widehat{\tau}}:\delta\mathbf{h}^{s}\,dV
-\int_{\mathcal{B}}\bm{\tau}:\delta\mathbf{h}\,dV\right),\\
\int_{\mathcal{B}}\delta\bm{\tau}:\nabla\Delta\bm{\varphi}\,dV&-
\int_{\mathcal{B}}\delta\bm{\tau}:\Delta\mathbf{h}\,dV=
-\left(\int_{\mathcal{B}}\delta\bm{\tau}:\nabla\bm{\varphi}\,dV-
\int_{\mathcal{B}}\delta\bm{\tau}:\mathbf{h}\,dV\right).\nonumber
\end{align}

In order for $\mathbf{H}$ to represent a displacement gradient, the shape functions that interpolate it should belong to the space $\mathcal{P}^{c}_{2}\left(\otimes^{2}\mathcal{T}\right)$ of
$\bigl(
\begin{smallmatrix}
2\\
0
\end{smallmatrix}
\bigr)$-tensors on the triangulation $\mathcal{T}$ \citep{ANG17}. The following form is assumed for the independent field $\mathbf{H}$:
\begin{equation}\label{eqn:4.8}
\mathbf{H}=\alpha_{a}\mathbf{G}^{c}_{a}\,,\qquad a=1,\ldots,24\,,
\end{equation}
where the summation is implied on $a$, and $\alpha_{a}$ denotes the degrees of freedom that correspond to the displacement gradient. Since three degrees of freedom are defined both for each edge $e_{i},~i=1,2,3$ and for the simplex $\mathscr{T}$, and, furthermore, independent degrees of freedom are required for each row of the interpolating tensors, the index $a$ can be computed as:
\begin{equation}\label{eqn:4.9}
a=
\begin{cases}
6(i-1)+2(j-1)+l, & j\text{th shape function on edge } e_{i},\\
2(k-1)+l+18,     & k\text{th shape function on simplex }\mathscr{T}.
\end{cases}
\end{equation}
In this relation, $l=1,2$ refers to the interpolating tensor $\mathbf{G}^{c}_{a_{l}}$, which contains the vector $\mathbf{V}^{e_{i}}_{j}$ on its $l$th row and the null vector on the other row. More specifically, for given $i$ and $j$ for the edge $e_{i}$ (or given $k$ for the simplex $\mathscr{T}$), one can compute $a_{1}$ and $a_{2}$ by substituting $l=1,2$ in the above relation. The submatrices $\mathbf{G}^{c}_{a_{1}}$ and $\mathbf{G}^{c}_{a_{2}}$ are then written as:
\begin{equation}\label{eqn:4.10}
\mathbf{G}^{c}_{a_{1}}=
\begin{bmatrix}
\left(\mathbf{V}^{e_{i}}_{j}\right)^{\mathsf{T}} \\
\mathbf{0}_{1\times2}
\end{bmatrix}\,,\qquad
\mathbf{G}^{c}_{a_{2}}=
\begin{bmatrix}
\mathbf{0}_{1\times2}\\
\left(\mathbf{V}^{e_{i}}_{j}\right)^{\mathsf{T}}
\end{bmatrix}\,.
\end{equation}
The vector $\mathbf{V}^{e_{i}}_j$ in these matrices is the $j$th global shape function on the edge $e_{i}$ that belongs to the space $\mathcal{P}^{c}_{2}\left(T\mathcal{T}\right)$. It is noted that for the $k$th shape function on the simplex $\mathscr{T}$, the matrices $\mathbf{G}^{c}_{a_{1}}$ and $\mathbf{G}^{c}_{a_{2}}$ are obtained by substituting the global shape function $\mathbf{V}^{\mathscr{T}}_{k}$ for $\mathbf{V}^{e_{i}}_j$ in \eqref{eqn:4.10}. The push-forward of $\mathbf{H}$ to the current configuration is given as:
\begin{equation}\label{eqn:4.11}
\mathbf{h}=\alpha_{a}\mathbf{g}^{c}_{a}\,,\qquad a=1,\ldots,24\,,
\end{equation}
where $\mathbf{g}^{c}_{a}=\mathbf{G}^{c}_{a}\mathbf{F}^{-1}$. On the other hand, the shape functions that interpolate $\mathbf{P}$ belong to the space $\mathcal{P}^{d-}_{2}\left(\otimes^{2}\mathcal{T}\right)$ of
$\bigl(
\begin{smallmatrix}
2\\
0
\end{smallmatrix}
\bigr)$-tensors. Hence, the following form is assumed for the independent field $\mathbf{P}$:
\begin{equation}\label{eqn:4.12}
\mathbf{P}=\gamma_{b}\mathbf{G}^{d}_{b}\,,\qquad b=1,\ldots,16\,.
\end{equation}
The parameters $\gamma_{b}$ in this equation are the stress degrees of freedom on the edges $e_{i},~i=1,2,3$ and the simplex $\mathscr{T}$. In this case, two degrees of freedom are required for each edge and another two degrees of freedom for the simplex. Moreover, since $\mathbf{P}$ is a tensor field, separate degrees of freedom are required for each row. Therefore, the index $b$ can be calculated as:
\begin{equation}\label{eqn:4.13}
b=
\begin{cases}
4(i-1)+2(j-1)+l, & j\text{th shape function on edge } e_{i},\\
2(k-1)+l+12,     & k\text{th shape function on simplex }\mathscr{T}.
\end{cases}
\end{equation}
Similarly, the parameter $l=1,2$ in \eqref{eqn:4.13} refers to the interpolating tensor $\mathbf{G}^{d}_{b_{l}}$, which contains the vector $\mathbf{W}^{e_{i}}_{j}$ on its $l$th row and the null vector on the other row. For given $i$ and $j$ for the edge $e_{i}$, it is convenient to compute $b_{1}$ and $b_{2}$ via setting $l=1,2$ in the above relation. Then, the submatrices $\mathbf{G}^{d}_{b_{1}}$ and $\mathbf{G}^{d}_{b_{2}}$ are computed as:
\begin{equation}\label{eqn:4.14}
\mathbf{G}^{d}_{b_{1}}=
\begin{bmatrix}
\left(\mathbf{W}^{e_{i}}_{j}\right)^{\mathsf{T}} \\
\mathbf{0}_{1\times2}
\end{bmatrix},\qquad
\mathbf{G}^{d}_{b_{2}}=
\begin{bmatrix}
\mathbf{0}_{1\times2}\\
\left(\mathbf{W}^{e_{i}}_{j}\right)^{\mathsf{T}}
\end{bmatrix}\,.
\end{equation}
The vector $\mathbf{W}^{e_{i}}_{j}$ in these matrices is the $j$th global shape function on the edge $e_{i}$ that belongs to the space $\mathcal{P}^{d-}_{2}\left(T\mathcal{T}\right)$. Similarly, for given $k$ for the simplex $\mathscr{T}$, the matrices $\mathbf{G}^{d}_{b_{1}}$ and $\mathbf{G}^{d}_{b_{2}}$ are obtained by substituting the global shape function $\mathbf{W}^{\mathscr{T}}_{k}$ for $\mathbf{W}^{e_{i}}_{j}$ in \eqref{eqn:4.14}. From \eqref{eqn:4.12} and the definition of the Kirchhoff stress tensor, it is clear that $\bm{\tau}$ can be interpolated as:
\begin{equation} \label{eqn:4.15}
    \bm{\tau}=\gamma_{b}\mathbf{g}^{d}_{b}\,,\qquad b=1,\ldots,16\,,
\end{equation}
where $\mathbf{g}^{d}_{b}=\mathbf{G}^{d}_{b}\mathbf{F}^{\mathsf{T}}$.

To write the governing equations in matrix form, the components of the tensors $\mathbf{h}$ and $\bm{\tau}$ in \eqref{eqn:4.11} and \eqref{eqn:4.15} are written in vector form as:
\begin{equation}\label{eqn:4.16}
\bm{\varepsilon}=
\begin{bmatrix}
h_{11}\\
h_{12}\\
h_{21}\\
h_{22}
\end{bmatrix},\qquad
\bm{\Gamma}=
\begin{bmatrix}
\tau_{11}\\
\tau_{12}\\
\tau_{21}\\
\tau_{22}
\end{bmatrix}\,.
\end{equation}
It is important to note that the interpolated Kirchhoff stress $\bm{\tau}$ in \eqref{eqn:4.15}, which is obtained via interpolation using shape functions (see Eq. \eqref{eqn:4.12}), is not necessarily symmetric.\footnote{This, however, does not imply that the balance of angular momentum is violated; the balance of angular momentum is enforced weakly through \eqref{eqn:4.28}$_2$.} Thus, it is required to keep both components $\tau_{12}$ and $\tau_{21}$ in the vector $\Gamma$. In accordance with the definition of the vectors $\bm{\varepsilon}$ and $\bm{\Gamma}$, the matrices $\mathbf{\bar{g}}^{c}$ and $\mathbf{\bar{g}}^{d}$ are defined as:
\begin{equation}\label{eqn:4.17}
\mathbf{\bar{g}}^{c}=
\begin{bmatrix}
\left(\mathbf{g}^{c}_{1}\right)_{11}&\left(\mathbf{g}^{c}_{2}\right)_{11}&\cdots&\left(\mathbf{g}^{c}_{24}\right)_{11}\\
\left(\mathbf{g}^{c}_{1}\right)_{12}&\left(\mathbf{g}^{c}_{2}\right)_{12}&\cdots&\left(\mathbf{g}^{c}_{24}\right)_{12}\\
\left(\mathbf{g}^{c}_{1}\right)_{21}&\left(\mathbf{g}^{c}_{2}\right)_{21}&\cdots&\left(\mathbf{g}^{c}_{24}\right)_{21}\\
\left(\mathbf{g}^{c}_{1}\right)_{22}&\left(\mathbf{g}^{c}_{2}\right)_{22}&\cdots&\left(\mathbf{g}^{c}_{24}\right)_{22}
\end{bmatrix}_{4\times24}\,,
\end{equation}
and
\begin{equation}\label{eqn:4.18}
\mathbf{\bar{g}}^{d}=
\begin{bmatrix}
\left(\mathbf{g}^{d}_{1}\right)_{11}&\left(\mathbf{g}^{d}_{2}\right)_{11}&\cdots&\left(\mathbf{g}^{d}_{16}\right)_{11}\\
\left(\mathbf{g}^{d}_{1}\right)_{12}&\left(\mathbf{g}^{d}_{2}\right)_{12}&\cdots&\left(\mathbf{g}^{d}_{16}\right)_{12}\\
\left(\mathbf{g}^{d}_{1}\right)_{21}&\left(\mathbf{g}^{d}_{2}\right)_{21}&\cdots&\left(\mathbf{g}^{d}_{16}\right)_{21}\\
\left(\mathbf{g}^{d}_{1}\right)_{22}&\left(\mathbf{g}^{d}_{2}\right)_{22}&\cdots&\left(\mathbf{g}^{d}_{16}\right)_{22}
\end{bmatrix}_{4\times16}\,.
\end{equation}
The $a$th column in the matrix $\mathbf{\bar{g}}^{c}$ lists the components of the matrix $\mathbf{g}^{c}_{a}$ in \eqref{eqn:4.11}, while the $b$th column in the matrix $\mathbf{\bar{g}}^{d}$ lists the components of $\mathbf{g}^{d}_{b}$ in \eqref{eqn:4.15}. In view of Eqs.~\eqref{eqn:4.16}-\eqref{eqn:4.18}, it is possible to express \eqref{eqn:4.11} and \eqref{eqn:4.15} in the following matrix form as $\bm{\varepsilon}=\mathbf{\bar{g}}^{c}\bm{\alpha}$ and $\bm{\Gamma}=\mathbf{\bar{g}}^{d}\bm{\gamma}$, where the vectors $\bm{\alpha}$ and $\bm{\gamma}$ list the degrees of freedom $\alpha_{a},~a=1,\ldots,24$ and $\gamma_{b},~b=1,\ldots,16$, respectively. To be consistent with the definition of vector $\bm{\varepsilon}$ in \eqref{eqn:4.16}, the gradient of the displacement field $\bm{\varphi}$ is expressed in the following form:
\begin{equation}\label{eqn:4.19}
\bm{\epsilon}^{\mathsf{T}}=
\begin{bmatrix}
\left(\nabla\bm{\varphi}\right)_{11} & \left(\nabla\bm{\varphi}\right)_{12} &
\left(\nabla\bm{\varphi}\right)_{21} & \left(\nabla\bm{\varphi}\right)_{22}
\end{bmatrix}\,.
\end{equation}
Considering a six-node triangular element, the displacement field $\bm{\varphi}$ is interpolated using the Lagrange shape functions $N_{i}\in\mathcal{P}_{2}\left(\mathds{R}^{2}\right),~i=1,\ldots,6$. Hence, the matrix $\mathbf{B}$ takes the form:
\begin{equation}\label{eqn:4.20}
\mathbf{B}=
\begin{bmatrix}
\nabla N_{1}          & \mathbf{0}_{2\times1} &
\ldots                &
\nabla N_{6}          & \mathbf{0}_{2\times1} \\
\mathbf{0}_{2\times1} & \nabla N_{1}          &
\ldots                &
\mathbf{0}_{2\times1} & \nabla N_{6}
\end{bmatrix}_{4\times12}\,,
\end{equation}
where the gradient of shape functions in the current configuration are obtained from $\nabla N_{i}=\left(\nabla_{0}N_{i}\right)\mathbf{F}^{-1}$. In view of \eqref{eqn:4.19} and \eqref{eqn:4.20}, the vector $\bm{\epsilon}$ can be computed using the relation $\bm{\epsilon}=\mathbf{B}\mathbf{u}$. The vector $\mathbf{u}$ in this equation lists the displacement degrees of freedom at the nodes of a given triangular element. Clearly, defining the matrix $\mathbf{N}$ as:
\begin{equation}\label{eqn:4.21}
\mathbf{N}=
\begin{bmatrix}
N_{1} & 0     & \ldots & N_{6} & 0     \\
0     & N_{1} & \ldots & 0     & N_{6}
\end{bmatrix}_{2\times12}\,,
\end{equation}
enables one to interpolate the displacement vector $\bm{\varphi}$ via the equation, $\bm{\varphi}=\mathbf{N}\mathbf{u}$.

To discretize \eqref{eqn:4.7}, the following vector forms are considered for $\delta\bm{\varphi}$, $\nabla\delta\bm{\varphi}$, $\delta\mathbf{h}$ and $\delta\bm{\tau}$:
\begin{equation}\label{eqn:4.22}
    \delta\bm{\varphi}=\mathbf{N}\,\delta\mathbf{u}\,,\qquad
    \delta\bm{\epsilon}=\mathbf{B}\,\delta\mathbf{u}\,,\qquad
    \delta\bm{\varepsilon}=\mathbf{\bar{g}}^{c}\delta\bm{\alpha},\qquad
    \delta\bm{\Gamma}=\mathbf{\bar{g}}^{d}\delta\bm{\gamma}\,.
\end{equation}
On the other hand, the vector forms for $\Delta\bm{\varphi}$, $\nabla\Delta\bm{\varphi}$, $\Delta\mathbf{h}$ and $\Delta\bm{\tau}$ are expressed as:
\begin{equation}\label{eqn:4.23}
    \Delta\bm{\varphi}=\mathbf{N}\,\Delta\mathbf{u}\,,\qquad
    \Delta\bm{\epsilon}=\mathbf{B}\,\Delta\mathbf{u},\qquad
    \Delta\bm{\varepsilon}=\mathbf{\bar{g}}^{c}\Delta\bm{\alpha}\,,\qquad
    \Delta\bm{\Gamma}=\mathbf{\bar{g}}^{d}\Delta\bm{\gamma}\,.
\end{equation}
Given that the vectors $\delta\mathbf{u}$, $\delta\bm{\alpha}$ and $\delta\bm{\gamma}$ are arbitrary, it is straightforward to show that \eqref{eqn:4.7} reduces to the following matrix equations:
\begin{align}\label{eqn:4.24}
 \mathbf{K}^{u\gamma}_{\iota}\Delta\bm{\gamma}=
\mathbf{F}_{b}+\mathbf{F}_{t}-\mathbf{F}^{u}_{i}\,,\qquad
 \mathbf{K}^{\alpha\alpha}_{m}\Delta\bm{\alpha}-
\mathbf{K}^{\alpha\gamma}_{\iota}\Delta\bm{\gamma}=-\mathbf{F}^{\alpha}_{i}\,,\qquad
\mathbf{K}^{\gamma u}_{\iota}\Delta\mathbf{u}-
\mathbf{K}^{\gamma\alpha}_{\iota}\Delta\bm{\alpha}=-\mathbf{F}^{\gamma}_{i}\,,
\end{align}
where the various stiffness matrices are defined as:
\begin{equation} \label{eqn:4.25}
\begin{aligned}
    \mathbf{K}^{u\gamma}_{\iota}&=\mathbf{K}^{\gamma u\,\mathsf{T}}_{\iota}=
    \int_{\mathcal{B}}\mathbf{B}^{\mathsf{T}}\mathbf{\bar{g}}^{d}dV\,,\qquad
    \mathbf{K}^{\alpha\alpha}_{m}=
    \int_{\mathcal{B}}\mathbf{\bar{g}}^{c\,\mathsf{T}}\mathbf{D}\,\mathbf{\bar{g}}^{c}dV+
    \int_{\mathcal{B}}\mathbf{\bar{g}}^{c\,\mathsf{T}}\mathbf{\widehat{T}}\,\mathbf{\bar{g}}^{c}dV\,,\\
    \mathbf{K}^{\alpha\gamma}_{\iota} &=\mathbf{K}^{\gamma\alpha\,\mathsf{T}}_{\iota}=
    \int_{\mathcal{B}}\mathbf{\bar{g}}^{c\,\mathsf{T}}\mathbf{\bar{g}}^{d}dV.
\end{aligned}
\end{equation}
In these relations, the matrices identified with the subindex $\iota$ are obtained from the interpolation of the independent fields. On the other hand, the first integral in the definition of $\mathbf{K}^{\alpha\alpha}_{m}$ characterizes the material behavior, while the second integral emerges as a result of geometric nonlinearity. The matrix $\mathbf{D}$ represents the elasticity tensor $\boldsymbol{\mathbb{C}}$ in matrix form and the matrix $\mathbf{\widehat{T}}$ is defined as follows:
\begin{equation}\label{eqn:4.26}
\mathbf{\widehat{T}}=
\begin{bmatrix}
\bm{\widehat{\tau}}&\mathbf{0}_{4\times4}\\
\mathbf{0}_{4\times4}&\bm{\widehat{\tau}}
\end{bmatrix}\,.
\end{equation}
The vectors $\mathbf{F}_{b}$ and $\mathbf{F}_{t}$ in \eqref{eqn:4.24} are the body and traction force vectors, which are defined, respectively, as:
\begin{equation}\label{eqn:4.27}
    \mathbf{F}_{b}=\int_{\mathcal{B}}\rho_{0}\left(\mathbf{N}^{\mathsf{T}}\mathbf{b}\right)dV\,,\qquad
    \mathbf{F}_{t}=\sum_{b=1}^{n_{b}}\int_{\partial_{\tau}M_{b}}\mathbf{N}^{\mathsf{T}}_{b}\mathbf{t}\,dA\,.
\end{equation}
The matrix $\mathbf{N}_{b}$ is defined similarly to the matrix $\mathbf{N}$ in \eqref{eqn:4.21} and is used to interpolate the displacement on the portion of the boundary where traction is specified. From \eqref{eqn:4.7}, it is evident that the internal load vectors, $\mathbf{F}^{u}_{i}$, $\mathbf{F}^{\alpha}_{i}$ and $\mathbf{F}^{\gamma}_{i}$ can be defined as:
\begin{equation}\label{eqn:4.28}
    \mathbf{F}^{u}_{i}=
    \int_{\mathcal{B}}\mathbf{B}^{\mathsf{T}}\bm{\Gamma}\,dV\,,\qquad
    \mathbf{F}^{\alpha}_{i}=
    \int_{\mathcal{B}}\mathbf{\bar{g}}^{c\,\mathsf{T}}\left(\bm{\widehat{\Gamma}}-\bm{\Gamma}\right)dV\,,\qquad
    \mathbf{F}^{\gamma}_{i}=
    \int_{\mathcal{B}}\mathbf{\bar{g}}^{d\,\mathsf{T}}\left(\bm{\epsilon}-\bm{\varepsilon}\right)dV\,.
\end{equation}
It is important to note that while the role of the first internal load vector is to balance the internal forces with the externally applied loads, the role of the second and third internal load vectors are, respectively, to enforce the constraints $\bm{\widehat{\tau}}=\bm{\tau}$ and $\nabla\bm{\varphi}=\mathbf{h}$. It is convenient to cast \eqref{eqn:4.24} into the form $\mathbf{K}^{\mathscr{T}}_{t}\Delta\mathbf{U}^{\mathscr{T}}=\mathbf{F}^{\mathscr{T}}_{e}-\mathbf{F}^{\mathscr{T}}_{i}$, where the superindex $\mathscr{T}$ implies that the equation corresponds to a given simplex $\mathscr{T}\in\mathcal{T}$. In this equation, the tangent stiffness matrix is expressed as:
\begin{equation}\label{eqn:4.29}
\mathbf{K}^{\mathscr{T}}_{t}=
\begin{bmatrix}
\mathbf{0}_{12\times12}       &  \mathbf{0}_{12\times24}          &  \mathbf{K}^{u\gamma}_{\iota}      \\
\mathbf{0}_{24\times12}       &  \mathbf{K}^{\alpha\alpha}_{m}    & -\mathbf{K}^{\alpha\gamma}_{\iota} \\
\mathbf{K}^{\gamma u}_{\iota} & -\mathbf{K}^{\gamma\alpha}_{\iota}&  \mathbf{0}_{16\times16}
\end{bmatrix},
\end{equation}
while the vectors $\Delta\mathbf{U}^\mathscr{T}$, $\mathbf{F}^\mathscr{T}_{e}$ and $\mathbf{F}^\mathscr{T}_{i}$ have the following representations:
\begin{equation}\label{eqn:4.30}
    \Delta\mathbf{U}^{\mathscr{T}}=
    \begin{bmatrix}
    \Delta\mathbf{u}\\
    \Delta\bm{\alpha}\\
    \Delta\bm{\gamma}
    \end{bmatrix}\,,\qquad
    \mathbf{F}^{\mathscr{T}}_{e}=
    \begin{bmatrix}
    \mathbf{F}_{b}+\mathbf{F}_{t}\\
    \mathbf{0}_{24\times1}\\
    \mathbf{0}_{16\times1}
    \end{bmatrix}\,,\qquad
    \mathbf{F}^{\mathscr{T}}_{i}=
    \begin{bmatrix}
    \mathbf{F}^{u}_{i}\\
    \mathbf{F}^{\alpha}_{i}\\
    \mathbf{F}^{\gamma}_{i}
    \end{bmatrix}\,.
\end{equation}
The global stiffness matrix, and the global external and internal load vectors can be assembled using the stiffness matrix, external and internal load vectors for each simplex $\mathscr{T}$ as:
\begin{equation}\label{eqn:4.31}
    \mathbb{K}_{t}=\underset{\mathscr{T}\in\mathcal{T}}{\mathbf{\scalebox{1.5}{A}}}\mathbf{K}^{\mathscr{T}}_{t}\,,\qquad
    \mathbb{F}_{e}=\underset{\mathscr{T}\in\mathcal{T}}{\mathbf{\scalebox{1.5}{A}}}\mathbf{F}^{\mathscr{T}}_{e}\,,\qquad
    \mathbb{F}_{i}=\underset{\mathscr{T}\in\mathcal{T}}{\mathbf{\scalebox{1.5}{A}}}\mathbf{F}^{\mathscr{T}}_{i}\,,
\end{equation}
where $\mathbf{\scalebox{1.5}{A}}$ is the assembly operator. The vector $\Delta\mathbb{U}$, obtained by solving the assembled system of equations, is used to update the global degrees of freedom using the relation $\mathbb{U}_{i+1} = \mathbb{U}_{i} + \Delta\mathbb{U}_{i}$. The subindex $i$ in this equation refers to the iteration number. The iterative procedure continues until the norm of $\Delta\mathbb{U}_{i}$ is negligible within a given tolerance value.

\subsection{Implementation of the finite element method}
\label{sec:4.2}

\begin{figure}
\centering
\includegraphics[width=0.5\textwidth]{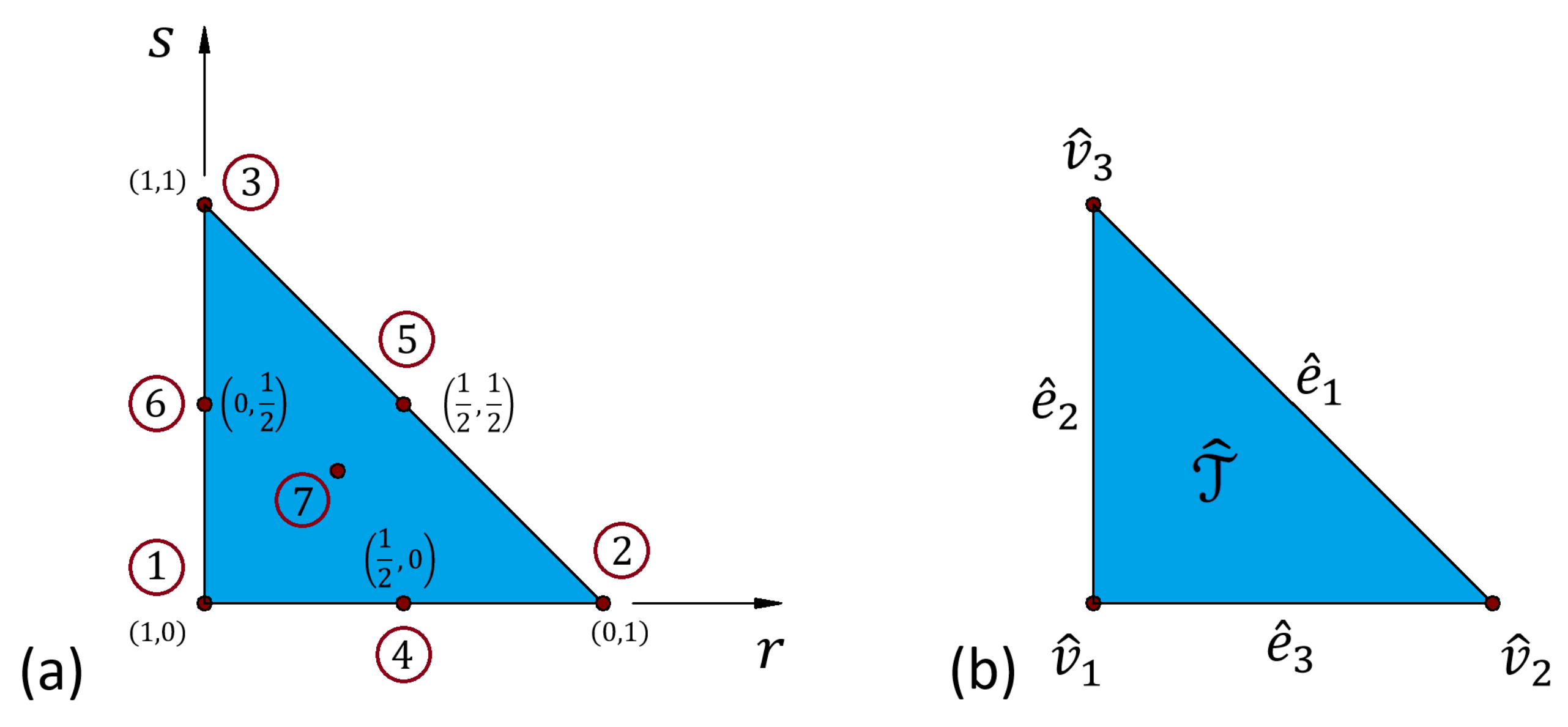}
\vskip 0.1in
\caption{Correspondence between (a) the six-node triangular element, and (b) the simplex $\widehat{\mathscr{T}}$ in the natural coordinate system.}
\label{fig:421}
\end{figure}

In \S\ref{sec:3}, it was shown that the degrees of freedom for the displacement gradient and the first Piola-Kirchhoff stress tensor are specified on the edges of the simplex $\mathscr{T}$ and the simplex itself. Since the displacement field is of second order, a six-node triangular element, as shown in Figure \ref{fig:421}a, is considered in its natural coordinates $r$ and $s$. The edges $142$, $253$ and $361$ of the triangular element correspond, respectively, to the edges $\widehat{e}_{3}$, $\widehat{e}_{1}$ and $\widehat{e}_{2}$ of the simplex $\widehat{\mathscr{T}}$ depicted in Figure \ref{fig:421}b, also in the natural coordinate system. The element in Figure \ref{fig:421}a contains a pseudo-node that is labeled as node $7$. This node represents the displacement gradient and stress degrees of freedom associated to the simplex itself. The application of pseudo-nodes in the formulation of compatible elements is discussed in detail in \citep{JAH22}. We recall also that the nodes $1$, $2$ and $3$ have displacement degrees of freedom only, while the nodes $4$, $5$ and $6$ have displacement, displacement gradient, and stress degrees of freedom. As discussed above, the pseudo-node $7$ has displacement gradient and stress degrees of freedom. Therefore, no displacement degrees of freedom are defined for node $7$. The degrees of freedom for nodes $1$ to $7$ are summarized in Table~\ref{tab:4.1}. It is important to note that including the displacement gradient and stress degrees of freedom for nodes $4$, $5$ and $6$ helps to satisfy the Hadamard jump condition for displacement gradient and the continuity of traction for the stress tensor (see Eq. \eqref{eqn:3.3}) over the edges $\widehat{e}_{3}$, $\widehat{e}_{1}$ and $\widehat{e}_{2}$ of the simplex $\widehat{\mathscr{T}}$ in Figure~\ref{fig:421}b. In other words, nodes $4$, $5$ and $6$ represent the displacement gradient and stress degrees of freedom on the edges $\widehat{e}_{3}$, $\widehat{e}_{1}$ and $\widehat{e}_{2}$. This point is further discussed in \citep{JAH22}.

\begin{table}
\centering
\caption{Degrees of freedom for nodes $1$ to $7$ of the six-node triangular element.}
\label{tab:4.1}
\renewcommand{\arraystretch}{1.5}
\renewcommand{\tabcolsep}{0.2cm}
\begin{tabular}{c c c c}
\hline
Node            & \multicolumn{3}{c}{degrees of freedom}                                                      \\
\cline{2-4}
                & Displacement    & Displacement Gradient            & Stress                                 \\
\hline
1               & $u_{1},u_{2}$   &                                  &                                        \\
2               & $u_{3},u_{4}$   &                                  &                                        \\
3               & $u_{5},u_{6}$   &                                  &                                        \\
4               & $u_{7},u_{8}$   & $\alpha_{1},\alpha_{2},\alpha_{3},\alpha_{4},\alpha_{5},\alpha_{6}$       & $\gamma_{1},\gamma_{2},\gamma_{3},\gamma_{4}$     \\
5               & $u_{9},u_{10}$  & $\alpha_{7},\alpha_{8},\alpha_{9},\alpha_{10},\alpha_{11},\alpha_{12}$    & $\gamma_{5},\gamma_{6},\gamma_{7},\gamma_{8}$     \\
6               & $u_{11},u_{12}$ & $\alpha_{13},\alpha_{14},\alpha_{15},\alpha_{16},\alpha_{17},\alpha_{18}$ & $\gamma_{9},\gamma_{10},\gamma_{11},\gamma_{12}$  \\
7               &                 & $\alpha_{19},\alpha_{20},\alpha_{21},\alpha_{22},\alpha_{23},\alpha_{24}$ & $\gamma_{13},\gamma_{14},\gamma_{15},\gamma_{16}$ \\
\hline
\end{tabular}
\end{table}

In \S\ref{sec:3}, it was explained that the spaces $\mathcal{P}^{c}_{r}\left(T\mathcal{T}\right)$, $\mathcal{P}^{d}_{r}\left(T\mathcal{T}\right)$ and their variants $\mathcal{P}^{c-}_{r}\left(T\mathcal{T}\right)$, $\mathcal{P}^{d-}_{r}\left(T\mathcal{T}\right)$ are subspaces of the Sobolev space $H\left(T\mathcal{B}\right)$. For the finite element solution to converge, it is necessary that the approximate solutions for the displacement gradient and stress tensor lie within these subspaces. Therefore, the spaces of the global shape functions that interpolate the displacement gradient and stress tensor must be sufficiently large. As the degrees of freedom for the displacement gradient and stress tensor are defined on the edges of the simplicial elements, this can be achieved by using meshes with randomly distributed simplex elements. The deviation of internal angles from an equilateral triangle must be kept within acceptable tolerances. Numerical experiments with meshes that possess this property have shown to accelerate convergence and enhance numerical stability. It was also discussed in \S\ref{sec:4.1} that the interpolated Kirchhoff stress $\bm{\tau}$ is not symmetric. Therefore, both components $\tau_{12}$ and $\tau_{21}$ are retained in the vector $\bm{\Gamma}$. In contrast, the Kirchhoff stress tensor $\bm{\widehat{\tau}}$ derived from the stored energy function is symmetric. The internal load vector $\mathbf{F}^{\alpha}_{i}$, computed using \eqref{eqn:4.28}$_2$ is utilized to implement the constraint $\bm{\widehat{\tau}}=\bm{\tau}$ in weak form. In other words, the constraints $\tau_{12}=\widehat{\tau}_{12}$ and $\tau_{21}=\widehat{\tau}_{21}$ are imposed to weakly symmetrize the tensor $\bm{\tau}$. According to the definitions of vectors $\bm{\Gamma}$ and $\bm{\widehat{\Gamma}}$, and the matrix $\mathbf{\bar{g}}^{c}$ in \eqref{eqn:4.17}, the material matrix $\mathbf{D}$ must be defined as follows:
\begin{equation}\label{eqn:4.32}
\mathbf{D}=
\begin{bmatrix}
\mathbb{C}_{1111}&\mathbb{C}_{1112}&\mathbb{C}_{1121}&\mathbb{C}_{1122}\\
\mathbb{C}_{1211}&\mathbb{C}_{1212}&\mathbb{C}_{1221}&\mathbb{C}_{1222}\\
\mathbb{C}_{2111}&\mathbb{C}_{2112}&\mathbb{C}_{2121}&\mathbb{C}_{2122}\\
\mathbb{C}_{2211}&\mathbb{C}_{2212}&\mathbb{C}_{2221}&\mathbb{C}_{2222}
\end{bmatrix}\,.
\end{equation}
Similarly, \eqref{eqn:4.28}$_3$ is used to enforce the constraint $\nabla\bm{\varphi}=\mathbf{h}$.
Similar to $\bm{\tau}$, the matrix $\mathbf{h}$ interpolated using \eqref{eqn:4.11} is not symmetric, and hence, all its components are required to satisfy this constraint. The above points elucidate the forms used for the vectors $\bm{\varepsilon}$ and $\bm{\Gamma}$ in \eqref{eqn:4.16}, as well as the vector $\bm{\epsilon}$ in \eqref{eqn:4.19}. All integrals in the computation of the internal forces and stiffness matrices can be computed using numerical integration. According to \citep{ZIE05}, if the number of the degrees of freedom is larger than the number of independent equations provided at Gauss points, the stiffness matrix must be singular. Each triangular element has $52$ degrees of freedom. On the other hand, for $2$D elasticity problems three independent equations are available at each Gauss point. Furthermore, three equations are obtained from enforcing the constraint $\nabla\bm{\varphi}=\mathbf{h}$ (assuming the symmetry of $\nabla\bm{\varphi}$) and another three equations are derived from the constraint $\bm{\widehat{\tau}}=\bm{\tau}$. Therefore, a minimum of $6$ Gauss points are required to evaluate the integrals. However, to integrate the stiffness matrix $\mathbf{K}^{\alpha\alpha}_{m}$ with higher accuracy, $13$ integration points have been used in the solution of the numerical examples that follow.

The existence and uniqueness of solutions for saddle point problems are thoroughly examined in \citep{BOF13}. Adopting their approach, one can alternatively express the tangent stiffness matrix in \eqref{eqn:4.29} as:
\begin{equation}\label{eqn:4.33}
\mathbf{K}^{\mathscr{T}}_{t}=
\begin{bmatrix}
  \mathbf{A} & \mathbf{B}^{\mathsf{T}} \\
  \mathbf{B} & \mathbf{0}_{16\times16}
\end{bmatrix}\,,
\end{equation}
where the matrices $\mathbf{A}$, $\mathbf{B}$ and $\mathbf{B}^{\mathsf{T}}$ have the following representations:
\begin{equation}\label{eqn:4.34}
\mathbf{A}=
\begin{bmatrix}
  \mathbf{K}^{\alpha\alpha}_{m} & \mathbf{0}_{24\times12} \\
  \mathbf{0}_{12\times24}       & \mathbf{0}_{12\times12}
\end{bmatrix}\,,\qquad
\mathbf{B}=
\begin{bmatrix}
  -\mathbf{K}^{\gamma\alpha}_{i} & \mathbf{K}^{\gamma u}_{i}
\end{bmatrix}\,,\qquad
\mathbf{B}^{\mathsf{T}}=
\begin{bmatrix}
  -\mathbf{K}^{\alpha\gamma}_{i} \\
  \mathbf{K}^{u\gamma}_{i}
\end{bmatrix}\,.
\end{equation}
If for a given element the degrees of freedom for the displacement gradient, the stress tensor, and the displacement field are denoted, respectively, by $r=24$, $s=16$ and $k=12$ (see Table \ref{tab:4.1}), then the dimensions of the matrices $\mathbf{K}^{\alpha\alpha}_{m}$, $-\mathbf{K}^{\gamma\alpha}_{i}$ and $\mathbf{K}^{\gamma u}_{i}$ are $r\times r$, $s\times r$ and $s\times k$, respectively. Obviously, the condition, $r+k>s>k$, which is necessary for the invertibility of the tangent stiffness matrix $\mathbf{K}^{\mathscr{T}}_{t}$ is satisfied \cite[Remark 3.2.1]{BOF13}. However, other conditions should be imposed on the matrices $\mathbf{K}^{\alpha\alpha}_{m}$, $-\mathbf{K}^{\gamma\alpha}_{i}$ and $\mathbf{K}^{\gamma u}_{i}$ so that the tangent stiffness matrix is non-singular. These conditions can be summarized as follows \cite[Section 3.2.5]{BOF13}:
\begin{align}\label{eqn:4.35}
\ker\left(\mathbf{K}^{\gamma u}_{i}\right)=\{\mathbf{0}_{k}\}\,,\qquad
\ker\left(-\mathbf{K}^{\alpha\gamma}_{i}\right)\cap\ker\left(\mathbf{K}^{u\gamma}_{i}\right)=\{\mathbf{0}_{s}\}\,,\qquad
\ker\left(\mathbf{K}^{\alpha\alpha}_{m}\right)\cap K=\{\mathbf{0}_{r}\}\,,
\end{align}
where the set $K$ is defined as:
\begin{equation}\label{eqn:4.36}
    K=\left\{\mathbf{x}\in\mathbb{R}^{r}\;\text{such that}\;\mathbf{y}^{\mathsf{T}}\mathbf{K}^{\gamma\alpha}_{i}\mathbf{x}=0\,,\quad\forall~\mathbf{y}\in\ker\left(\mathbf{K}^{u\gamma}_{i}\right)\right\}\,.
\end{equation}
The conditions given in \eqref{eqn:4.35} are both necessary and sufficient for ensuring that the tangent stiffness matrix $\mathbf{K}^{\mathscr{T}}_{t}$ is non-singular \citep{BOF13}.
It is shown in \citep{SHO18} that these conditions can be expressed in the form of inf-sup or the Ladyzhenskaya-Babu\v{s}ka-Brezzi conditions. Numerical experiments using the singular value decomposition of the matrices $\mathbf{K}^{\alpha\alpha}_{m}$, $-\mathbf{K}^{\gamma\alpha}_{i}$ and $\mathbf{K}^{\gamma u}_{i}$ applied for simple meshes with few elements show that the conditions given in \eqref{eqn:4.35} are also satisfied for the global stiffness matrix $\mathbb{K}_{t}$.

\section{Numerical Examples}
\label{sec:5}

This section presents several numerical examples to evaluate the efficiency and stability of the second-order CSMFE introduced in the previous sections. The performance of our second-order element, mainly the deformation and load-displacement curves it yields, is compared with that of existing compatible finite elements. Specifically, we compare it with the first-order compatible elements developed in \citep{DHA22a,JAH22}, as well as with both the first and second-order CSMFE elements of types $\mathrm{H}1\mathrm{c}1\mathrm{d}\bar{1}$ and $\mathrm{H}2\mathrm{c}2\mathrm{d}\bar{2}$ of \citet{ANG17}. It is worth mentioning that the shape functions for $\mathrm{H}2\mathrm{c}2\mathrm{d}\bar{2}$ elements in \citep{ANG17} are computed using the numerical integration, while in this work they are evaluated via the explicit forms given in Appendix \ref{app:A}. The problems discussed in this section have previously been used to examine different types of numerical instabilities for various finite elements. For example, some of these problems can reveal that the first-order elements may introduce spurious energy modes. Therefore, it is crucial to develop methods that eliminate these unrealistic modes \citep{REE00}. In other problems, shear locking and numerical instabilities may occur in the near-incompressible regime, requiring the use of stabilization techniques to enhance the element performance. The numerical results obtained with the second-order compatible-strain element in the subsequent examples demonstrate its excellent performance and the absence of numerical instabilities. All numerical examples have been solved using the full Newton-Raphson method employing the tangent stiffness matrix derived in \eqref{eqn:4.29}. A convergence norm of $10^{-9}$ has been used in all the numerical examples. Different load steps used to carry out the simulations are discussed for each specific example. Numerical experiments show that on average $3$ equilibrium iterations are required for each load step.

Three material types are used in the numerical examples. Two of these materials are of the compressible neo-Hookean type and the third is the Ogden material. These material types are characterized by the following strain energy functions:
\begin{equation}\label{eqn:5.1}
\begin{aligned}
\widehat{W}_{1}&=\frac{\mu}{2}\Bigl[\left(I_{1}-3\right)-\ln J^{2}\Bigr]+
\frac{\kappa}{2}\left(J-1\right)^{2}\,,\qquad
\widehat{W}_{2}=\frac{\mu}{2}\left(I_{1}-3\right)-\mu\ln J+
\frac{\kappa}{2}\left(\ln J\right)^{2}\,,\\
\widehat{W}_{3}&=\sum_{i=1}^{m}\frac{\mu_{i}}{\alpha_{i}}\left(\bar{\lambda}^{\alpha_{i}}_{1}+\bar{\lambda}^{\alpha_{i}}_{2}+\bar{\lambda}^{\alpha_{i}}_{3}-3\right)+
\frac{\kappa}{2}\left(J-1\right)^{2}\,,
\end{aligned}
\end{equation}
where $\mu$ is the shear modulus and $\kappa$ is the bulk modulus; $I_{1}=\mathrm{tr}\,\mathbf{C}$ is the first invariant of the right Cauchy-Green deformation tensor $\mathbf{C}=\mathbf{F}^{\mathsf{T}}\mathbf{F}$ and $J$ is the Jacobian of deformation.\footnote{Suppose the local coordinates $\{X^A\}$ and $\{x^a\}$ are used for the reference configuration $\mathcal{B}$ and the Euclidean ambient space $\mathcal{S}$. The metric of the ambient space is denoted as $\mathbf{g}$, which induces the flat metric $\mathbf{G}=\mathbf{g}\big|_{\mathcal{B}}$ on the reference configuration. The Jacobian of deformation has the following expression: $J=\sqrt{\frac{\det\mathbf{g}}{\det\mathbf{G}}}\,\det\mathbf{F}$. When Cartesian coordinates are used for both the reference and current configurations, this expression is simplified to $J=\det\mathbf{F}$.} The deviatoric stretches $\bar{\lambda}_{1}$, $\bar{\lambda}_{2}$ and $\bar{\lambda}_{3}$ are the eigenvalues of $\bar{\mathbf{F}}=J^{-\frac{1}{3}}\mathbf{F}$, and $\alpha_{i},\mu_{i},i=1,\ldots,m$ are the material parameters. We note that for plane strain problems $\bar{\lambda}_{3}=J^{-\frac{1}{3}}$. The Kirchhoff stress tensor for the three material models has the following representations:
\begin{equation}\label{eqn:5.2}
\begin{aligned}
\bm{\widehat{\tau}}_{1}&=\mu\mathbf{B}+
\Bigl[\kappa J\left(J-1\right)-\mu\Bigr]\mathbf{I}\,,\qquad
\bm{\widehat{\tau}}_{2}=\mu\mathbf{B}+\Bigl[\kappa\left(\ln J\right)-\mu\Bigr]\mathbf{I}\,,\\
\bm{\widehat{\tau}}_{3}&=\sum_{i=1}^{m}\mu_{i}J^{-\frac{\alpha_{i}}{3}}\left(\text{dev }\mathbf{B}^{\frac{\alpha_{i}}{2}}\right)+\kappa J\left(J-1\right)\,,
\end{aligned}
\end{equation}
where $\mathbf{B}=\mathbf{F}\mathbf{F}^{\mathsf{T}}$ is the left Cauchy-Green deformation tensor and $\mathbf{I}$ is the identity tensor of second-order. It is important to note that $\bm{\widehat{\tau}}_{1}$, $\bm{\widehat{\tau}}_{2}$ and $\bm{\widehat{\tau}}_{3}$ are the Kirchhoff stress tensors that are derived constitutively. This is in contrast to the interpolated Kirchhoff stress $\bm{\tau}$ that is obtained using \eqref{eqn:4.15}. The spatial elasticity tensors that correspond to the Kirchhoff stress tensors in \eqref{eqn:5.2} are
\begin{equation}\label{eqn:5.3}
\begin{aligned}
\boldsymbol{\mathbb{C}}_{1}&=\kappa J\left(2J-1\right)\mathbf{I}\otimes\mathbf{I}+
2\Bigl[\mu-\kappa J\left(J-1\right)\Bigr]\boldsymbol{\mathcal{I}}\,,\qquad
\boldsymbol{\mathbb{C}}_{2}=\kappa\mathbf{I}\otimes\mathbf{I}-2\Bigl[\kappa\left(\ln J\right)-\mu\Bigr]\boldsymbol{\mathcal{I}}\,,\\
\boldsymbol{\mathbb{C}}_{3}&=\sum_{i=1}^{m}\mu_{i}J^{-\frac{\alpha_{i}}{3}}
\Bigg\{-\frac{\alpha_{i}}{3}\left[\left(\text{dev }\mathbf{B}^{\frac{\alpha_{i}}{2}}\right)\otimes\mathbf{I}+
\mathbf{I}\otimes\left(\text{dev }\mathbf{B}^{\frac{\alpha_{i}}{2}}\right)\right]\\
&+\left(\alpha_{i}-2\right)\Bigl[\lambda^{\alpha_{i}}_{1}\left(\mathbf{E}_{1}\otimes\mathbf{E}_{1}\right)
+\lambda^{\alpha_{i}}_{2}\left(\mathbf{E}_{2}\otimes\mathbf{E}_{2}\right)\Bigr]\\
&+2\left(\frac{\lambda^{\alpha_{i}-2}_{1}-\lambda^{\alpha_{i}-2}_{2}}{\lambda^{2}_{1}-\lambda^{2}_{2}}\right)
\Bigl[\boldsymbol{\mathcal{I}}_{B}-\lambda^{4}_{1}\left(\mathbf{E}_{1}\otimes\mathbf{E}_{1}\right)-\lambda^{4}_{2}\left(\mathbf{E}_{2}\otimes\mathbf{E}_{2}\right)\Bigr]\Bigg\}\\
&+\left[\kappa J\left(2J-1\right)-\frac{1}{9}\sum_{i=1}^{m}\alpha_{i}\mu_{i}J^{-\frac{\alpha_{i}}{3}}\left(\text{tr }\mathbf{B}^{\frac{\alpha_{i}}{2}}\right)\right]\left(\mathbf{I}\otimes\mathbf{I}\right)\\
&-2\Bigl[\kappa J\left(J-1\right)-\frac{1}{3}\sum_{i=1}^{m}\mu_{i}J^{-\frac{\alpha_{i}}{3}}\left(\text{tr }\mathbf{B}^{\frac{\alpha_{i}}{2}}\right)\Bigr]\boldsymbol{\mathcal{I}}\,,
\end{aligned}
\end{equation}
where the symbol $\boldsymbol{\mathcal{I}}$ refers to the identity tensor of fourth-order\footnote{If the metric of the ambient space $\mathcal{S}$ is denoted by $\mathbf{g}$, then the components of the second and fourth-order identity tensors with respect to a local chart $\{x^a\}$ are $I^{ab}=g^{ab}$ and $\mathcal{I}^{abcd}=\frac{1}{2}\left(g^{ac}g^{bd}+g^{ad}g^{bc}\right)$, respectively.} and the tensor $\boldsymbol{\mathcal{I}}_B$ has the components $\left(\boldsymbol{\mathcal{I}}_B\right)^{abcd}=\frac{1}{2}\left(B^{ac}B^{bd}+B^{ad}B^{bc}\right)$ in the ambient space $\mathcal{S}$. The second-order tensors $\mathbf{E}_{1}$ and $\mathbf{E}_{2}$ are the in-plane eigenprojections of $\mathbf{B}$.
\begin{remark}\label{rem:5.1}
In the following numerical examples, the stability of the element is examined for the near-incompressible limit. This limit is attained by choosing large values for the ratio of the bulk modulus $\kappa$ to the rigidity modulus $\mu$.
\end{remark}
\begin{remark}\label{rem:5.2}
In each numerical example, several meshes have been considered for solving the problem. The number of elements and the number of degrees of freedom for each mesh are provided in tables to give an idea regarding the computational effort used to solve the examples. It is, however, important to note that the number of the degrees of freedom was not available for other references. Therefore, only the number of elements has been used to compare the performance of our second-order element with other references.
\end{remark}

\subsection{Homogeneous compression test}
\label{sec:5.1}

In this example, a near-incompressible block of $1 \times 1$~mm, homogeneously compressed to $80\%$ of its height, is considered. It is known that several eigenvalues of the stiffness matrix become negative at certain values of the applied deformation \citep{GLA97}. Using symmetry, only the right half of the system is analyzed. The bottom edge of the block is constrained in the vertical direction, while the left edge, which lies along the axis of symmetry, is constrained in the horizontal direction. The specified vertical displacement is applied to the top edge. The first and second material models in \eqref{eqn:5.1} are assigned parameters $\mu = 80.2$~MPa and $\kappa = 40,000$~MPa. It is worthwhile to mention that the second material model coincides with the one considered by \citet{GLA97}. These authors used $200$ Q1/ET4 elements to discretize the domain. \citet{JAH22} solved the problem using irregular meshes composed of $22$, $44$ and $72$ first-order compatible strain elements and showed that even coarse meshes are capable of accurately modeling the deformed configuration.

\begin{figure}
\centering
\includegraphics[width=0.6\textwidth]{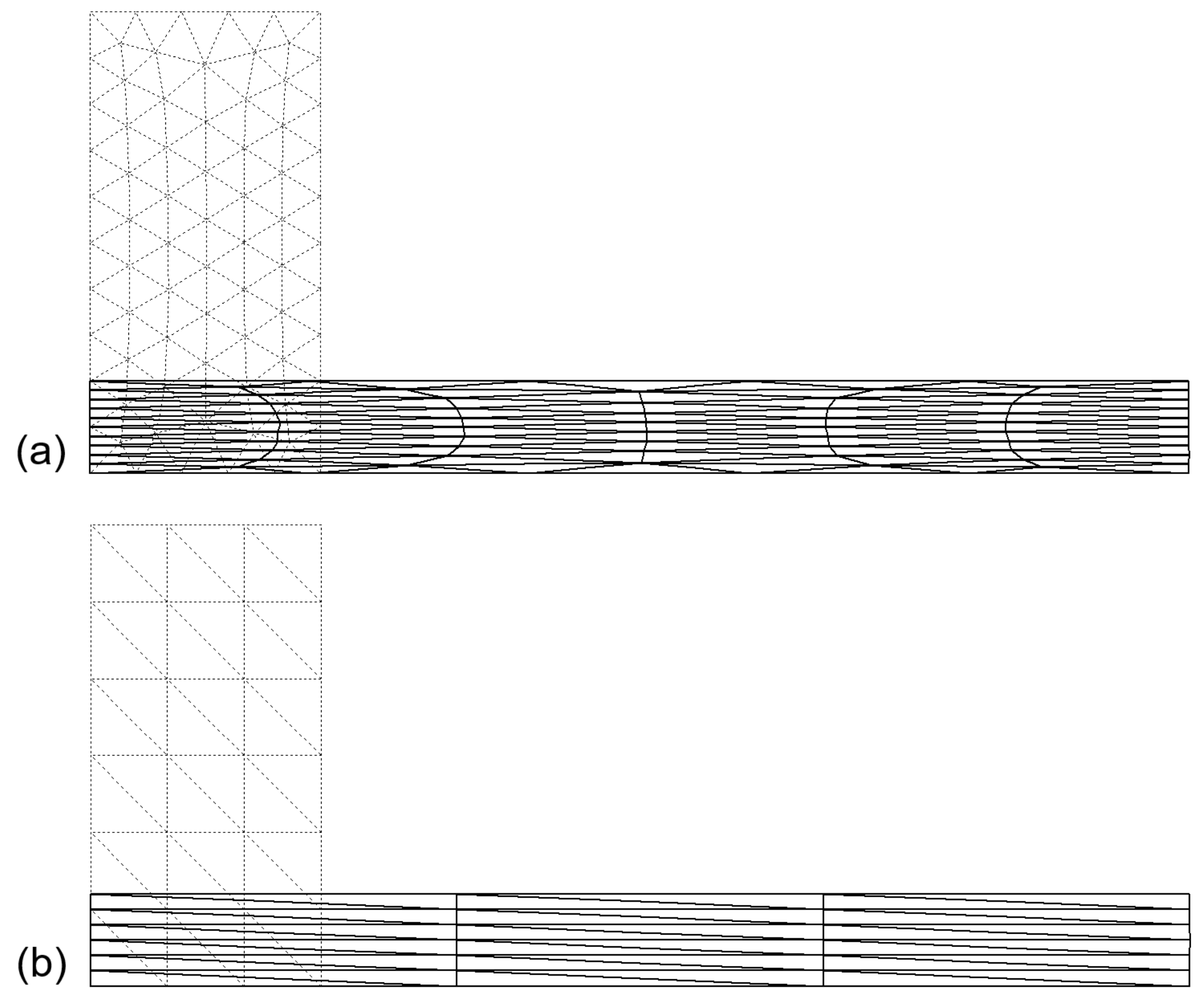}
\vskip 0.1in
\caption{Homogeneous compression test using: (a) An irregular mesh of $120$ second-order compatible strain elements and (b) a regular mesh of $36$ second-order compatible strain elements.}
\label{fig:511}
\end{figure}

\begin{figure}
\centering
\includegraphics[width=0.50\textwidth]{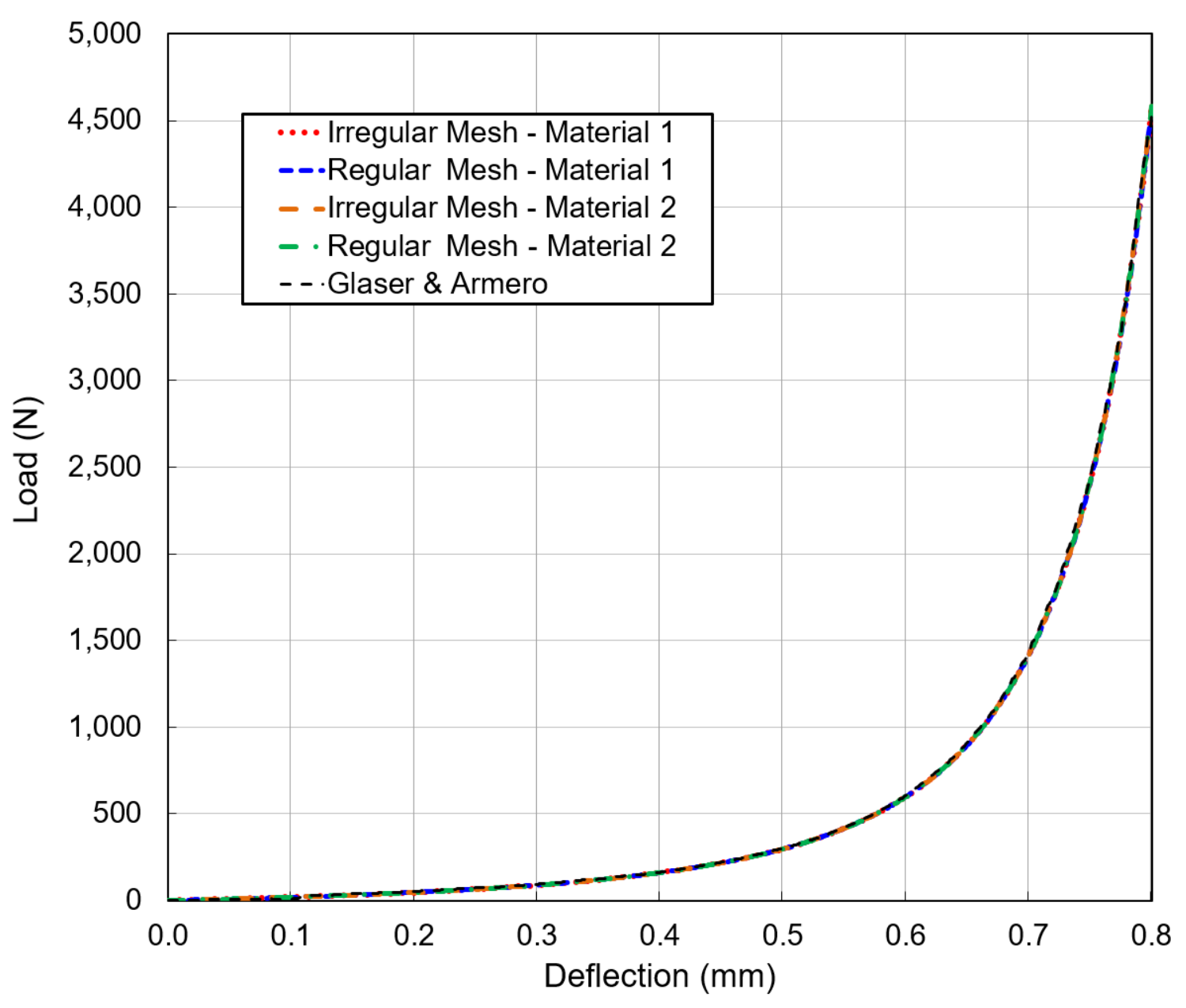}
\vskip 0.1in
\caption{Load-deflection curves for the homogeneous compression test using irregular and regular meshes of, respectively, $120$ and $36$ second-order compatible strain elements. The first and second material models are used in generating the curves and they are in good agreement with \citep{GLA97}.}
\label{fig:512}
\end{figure}

We use irregular meshes using $22$, $44$, $72$ and $120$ second-order compatible strain elements to discretize the domain. Furthermore, we also consider a regular mesh consisting of $36$ elements of the same type. The specified displacement is applied in $800$ load steps. The deformed configurations for irregular and regular meshes consisting, respectively, of $120$ and $36$ second-order compatible strain elements are shown in Figure~\ref{fig:511}. Figure~\ref{fig:512} illustrates that the load-deflection curves generated for the irregular and regular meshes using the first and second material models coincide with \citep{GLA97} confirming the stability of the element. As discussed in \citep{GLA97}, several first-order elements exhibit hourglass instability at certain deformation values. However, no such instability was observed when using the first and second-order compatible strain elements.

\subsection{Shearing of a block}
\label{sec:5.2}

\begin{figure}
\centering
\includegraphics[width=0.3\textwidth]{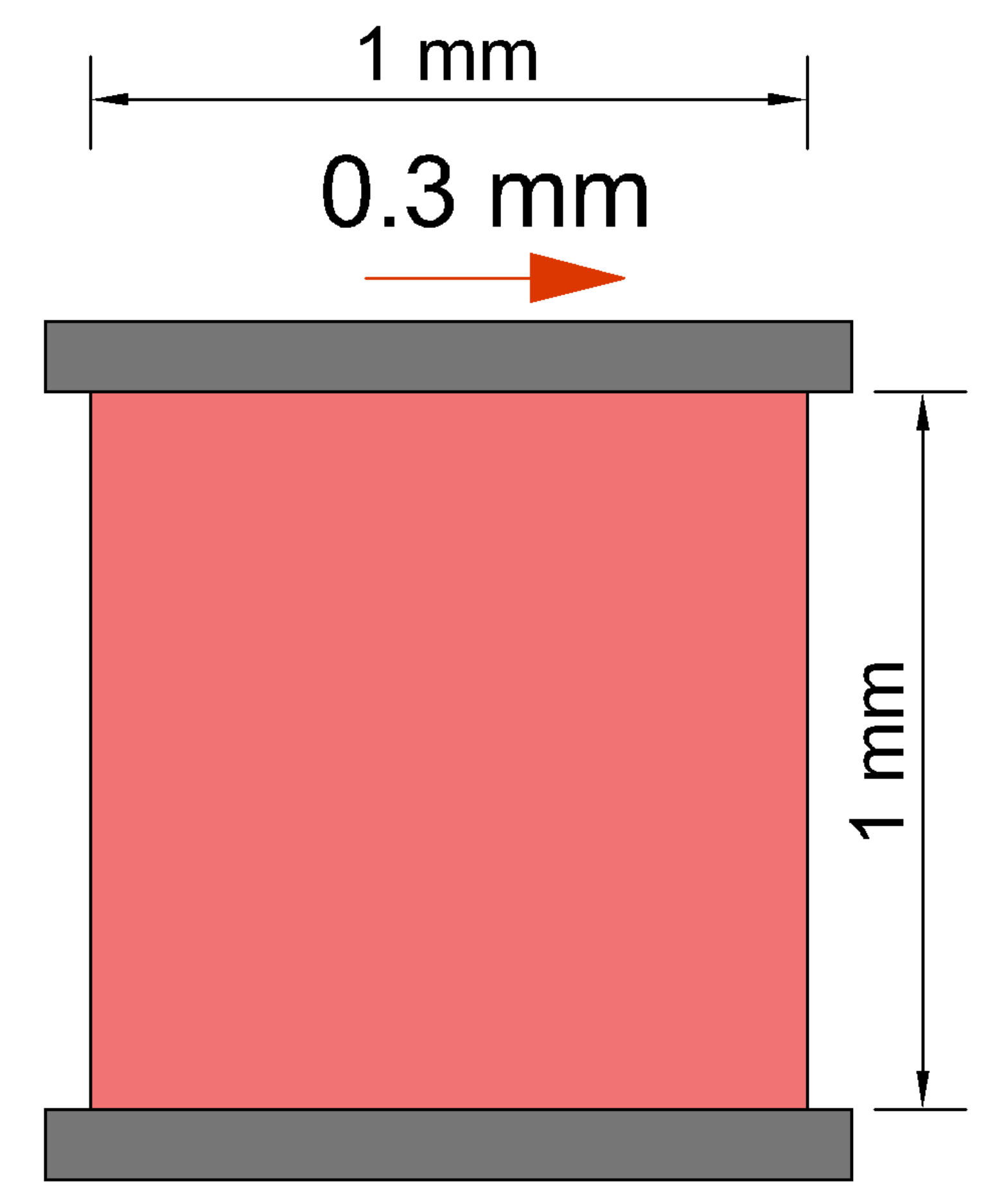}
\vskip -0.1in
\caption{Geometry and boundary conditions for the shearing of a block: The bottom edge is restrained both in the horizontal and vertical directions. The top edge, which is restrained in the vertical direction, undergoes a horizontal deformation of $0.3$~mm.}
\label{fig:521}
\end{figure}

\begin{figure}
\centering
\includegraphics[width=0.6\textwidth]{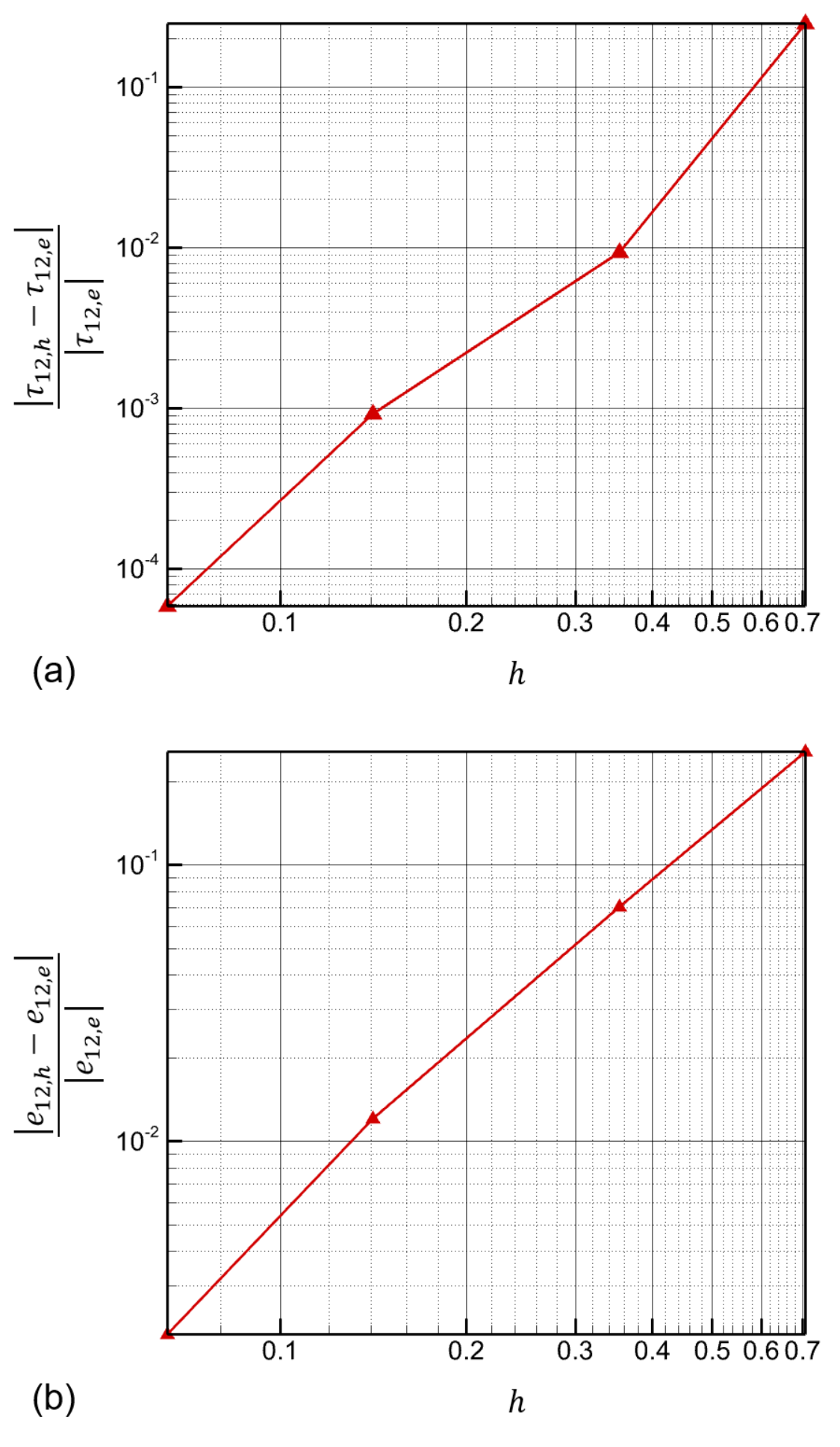}
\vskip -0.1in
\caption{Relative errors in (a) the Kirchhoff stress component $\tau_{12}$ and (b) the Eulerian strain component $e_{12}$ for the shearing of a block versus the average mesh size $h$.}
\label{fig:522}
\end{figure}

This example is used to study the convergence properties of the second-order CSMFEs developed in this work. A unit block whose geometry is shown in Figure~\ref{fig:521} is subject to a shearing deformation at its top edge. The bottom edge is restrained in both horizontal and vertical directions. The top edge, which is restrained in the vertical direction, undergoes a horizontal deformation of $0.3$~mm in $30$ load steps. The material parameters of $\mu = 80.2$~MPa and $\kappa = 40,000$~MPa are assigned to the first material model in \eqref{eqn:5.1}.
In order to study the convergence of the second-order element, the values of the Kirchhoff stress component $\tau_{12}$ and the Eulerian strain component $e_{12}$ at the center of the block are compared with a reference solution at that point. Since a closed-form solution for the problem is not available, a highly refined mesh consisting of $10,000$ Q2/P1 elements is used to obtain a numerically converged solution, which serves as the reference. Several meshes of the second-order CSMFEs have been considered.
The number of elements, the number of degrees of freedom and the average mesh size for these meshes are given in Table \ref{tab:5.1}. The relative error for the stress and strain components is plotted versus the average mesh size in Figures~\ref{fig:522}a and \ref{fig:522}b.
In these figures, $\tau_{12,h}$ and $e_{12,h}$ denote the values obtained from the meshes of second-order CSMFEs, while $\tau_{12,e}$ and $e_{12,e}$ represent the values computed using the highly refined Q2/P1 mesh. As shown in the error plots in Figure~\ref{fig:522}, the second-order CSMFEs exhibit good convergence behavior, even for a problem characterized by a high ratio of bulk modulus to shear modulus.

\begin{table}
\centering
\caption{The number of elements, the number of degrees of freedom, and the average mesh size for the meshes of second-order compatible-strain elements used to model the shearing of a block.}
\label{tab:5.1}
\renewcommand{\arraystretch}{1.5}
\renewcommand{\tabcolsep}{0.2cm}
\begin{tabular}{c c c c}
\hline
Mesh & No. Elements & DOFs  & Mesh Size (mm) \\
\hline
1    &   8          &   270 & 0.7071         \\
2    &  32          &  1006 & 0.3536         \\
3    & 200          &  5998 & 0.1414         \\
4    & 720          & 21238 & 0.0657         \\
\hline
\end{tabular}
\end{table}

\subsection{Inhomogeneous compression test}
\label{sec:5.3}

This example has been used as a benchmark test to study the behavior of the first-order elements, such as Q1, Q1/E4, Q1/ET4 and Q1SP, in relation to hourglass instability \citep{REE00,JAH22}. Many first-order elements without a stabilization technique for hourglassing would fail in the early stages of loading or show severe locking. This test problem has also been used to verify the stability of various first and second-order compatible strain elements \citep{ANG17,DHA22a,JAH22}.

\begin{figure}
\centering
\includegraphics[width=0.4\textwidth]{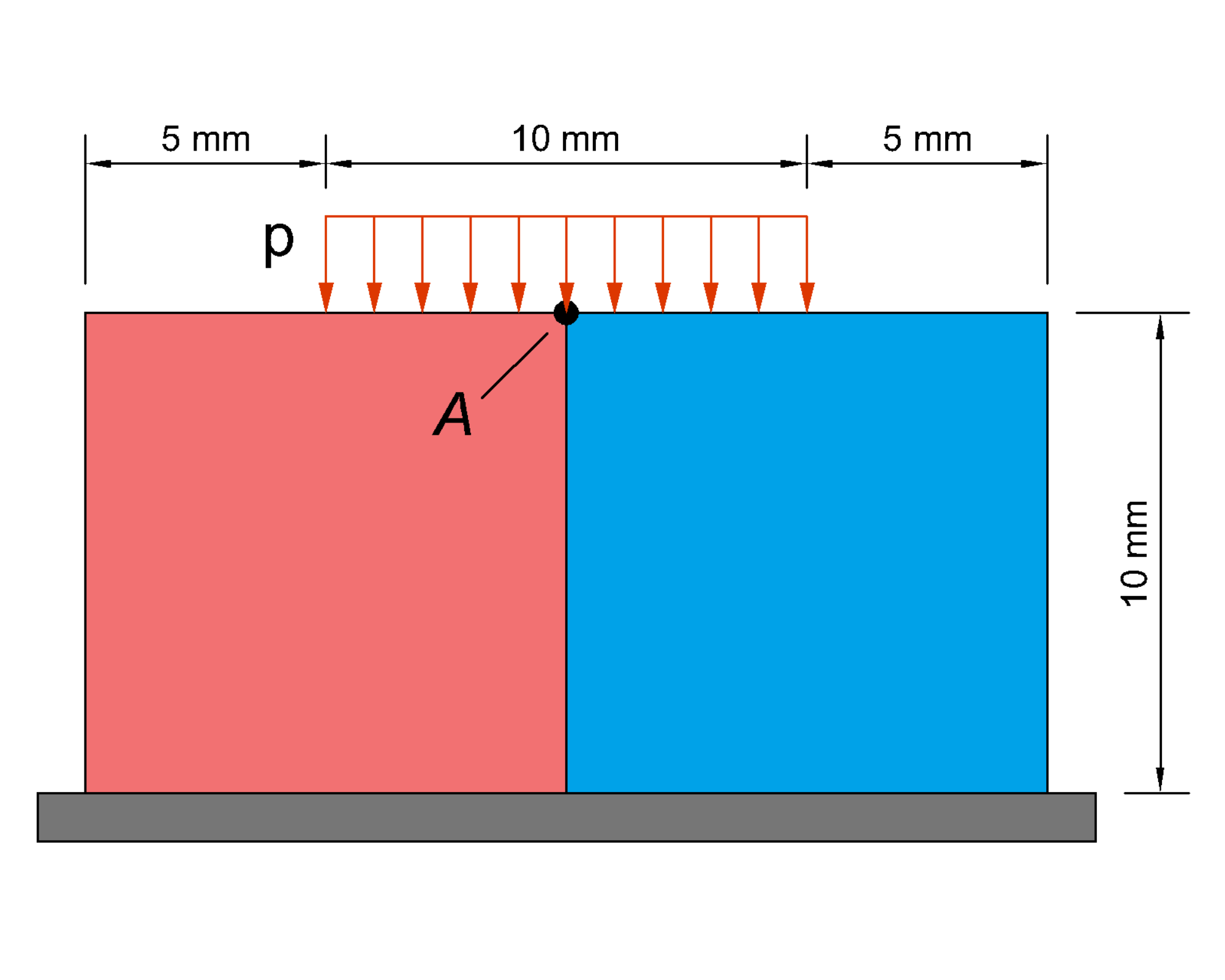}
\vskip -0.1in
\caption{Geometry and boundary conditions for the inhomogeneous compression test: The bottom edge is restrained in the vertical direction and the top edge is restrained in the horizontal direction. Due to symmetry, only the right half of the model is considered.}
\label{fig:531}
\end{figure}

The geometry of the problem is depicted in Figure~\ref{fig:531}. Due to symmetry, only the right half of the model is considered. The bottom edge is constrained to prevent vertical displacement. The nodes along the left edge, which lie on the axis of symmetry, and those on the top edge are restrained in the horizontal direction. The pressure $p$ is increased up to $600$~MPa, and the displacement of point $A$ is monitored during the application of $p$. The second material model in \eqref{eqn:5.1}, with parameters $\mu = 80.194$ MPa and $\kappa = 400,889.8$ MPa, is used for the simulations. In \citep{JAH22}, three different meshes consisting of $86$, $144$ and $324$ first-order compatible strain elements were used to solve the problem. It was demonstrated that a final displacement of $60\%$ could be achieved with all three meshes without encountering any numerical instabilities. This contrasts with the first-order $\mathrm{H}1\mathrm{c}1\mathrm{d}\bar{1}$ element \citep{ANG17}, which could only model a displacement of up to $30\%$.

\begin{figure}
\centering
\includegraphics[width=0.60\textwidth]{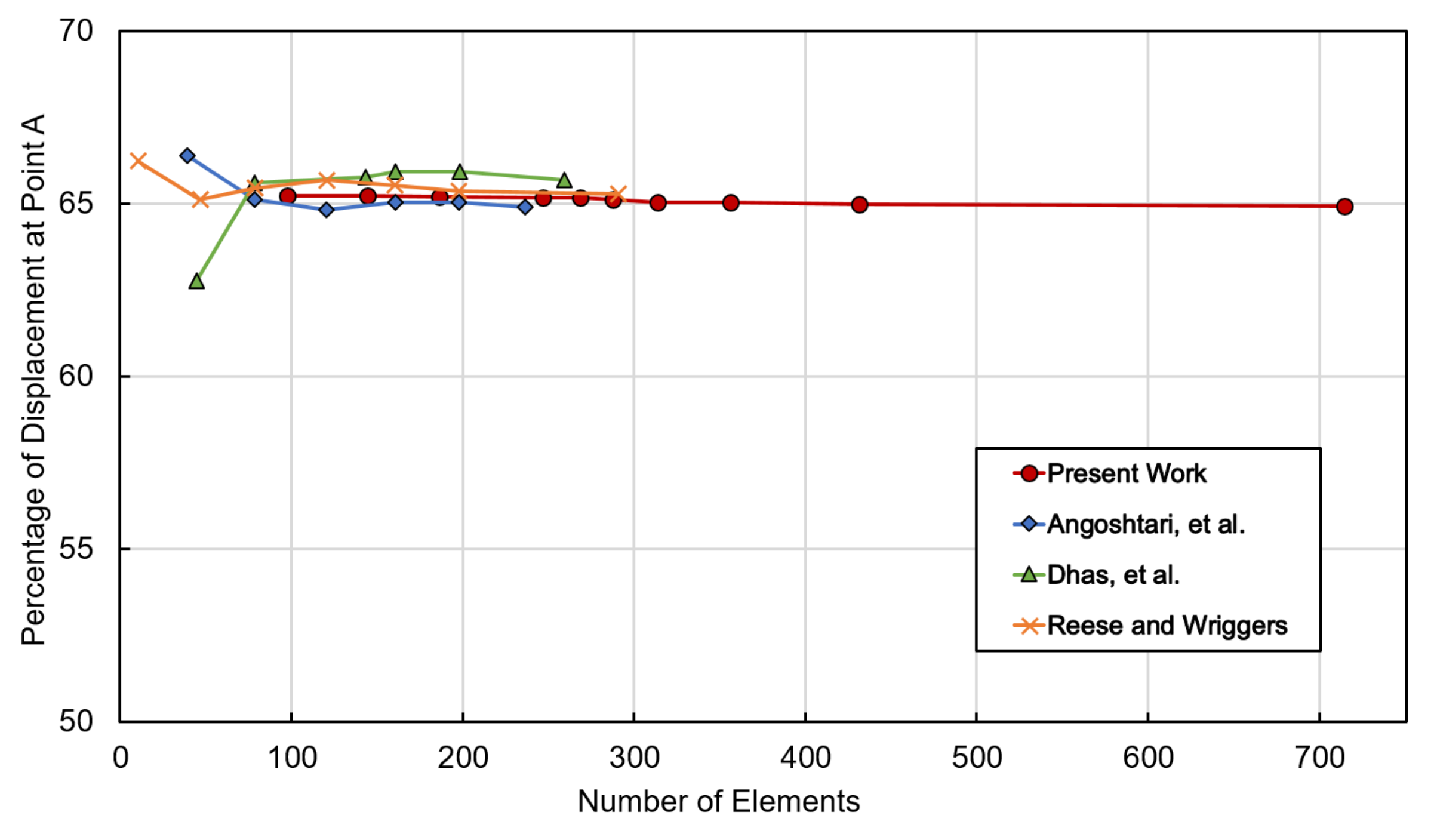}
\vskip 0.1in
\caption{Comparison of the vertical displacement of point $A$ for the inhomogeneous compression test with \cite{REE00,ANG17} using different number of elements and $p=600$ MPa.}
\label{fig:532}
\end{figure}

\begin{figure}
\centering
\includegraphics[width=0.75\textwidth]{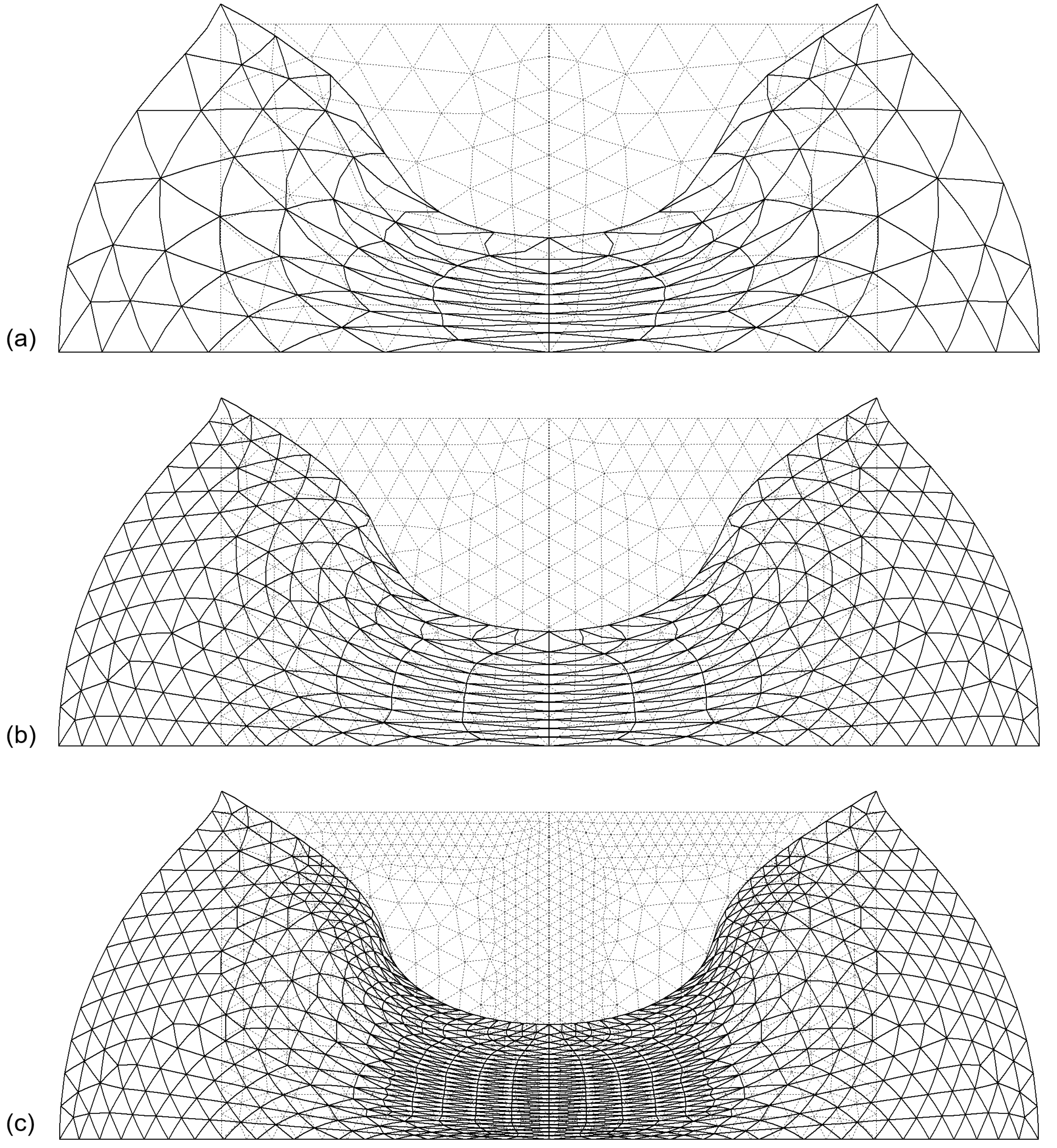}
\vskip 0.1in
\caption{Deformed configuration for the inhomogeneous compression test using meshes consisting of: (a) $97$, (b) $268$ and (c) $714$ second-order compatible strain elements.}
\label{fig:533}
\end{figure}

\begin{figure}
\centering
\includegraphics[width=0.75\textwidth]{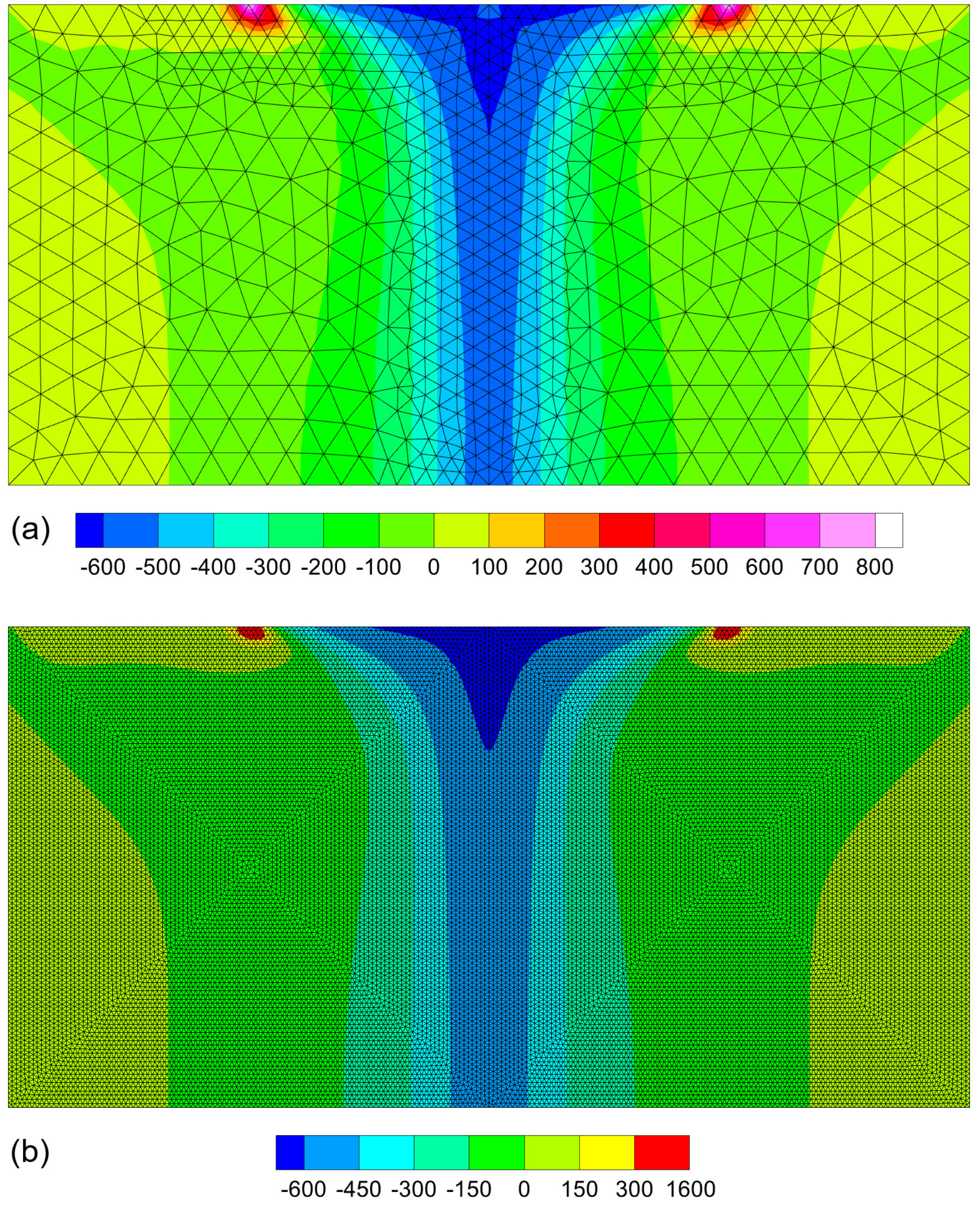}
\vskip 0.1in
\caption{Contours of the Kirchhoff stress component $\tau_{22}$ (MPa) for the inhomogeneous compression test using (a) a mesh consisting of $714$ second-order compatible strain elements developed in this work and (b) a mesh consisting of $22910$ second order triangular elements with U/P mixed formulation.}
\label{fig:534}
\end{figure}

Here, several meshes of second-order compatible strain elements are considered for solving the problem. The number of elements, the number of degrees of freedom and the displacement of point $A$ are given in Table~\ref{tab:5.2}. The pressure $p$ is applied in $1000$ load steps. The displacement of this point for various meshes is compared in Figure~\ref{fig:532} with the results from \citep{REE00,ANG17,DHA22a}. The value of $64.93\%$ for the mesh containing $714$ second-order compatible strain elements is comparable with the displacements of $64.9\%$ and $65.68\%$ in \citep{ANG17} and \citep{DHA22a}, respectively. To verify the results of the second-order CSMFEs developed in this work, a refined mesh consisting of $22910$ second-order triangular elements with U/P mixed formulation is also used to solve the problem. The displacement of point $A$ using this last mesh is obtained as $64.91\%$, which is very close to the value of $64.93\%$ computed using our element. The deformed configuration for the meshes consisting of $87$, $268$ and $714$ second-order CSMFEs is shown in Figure \ref{fig:533}. The contour plots for the Kirchhoff stress component $\tau_{22}$ are shown in Figures \ref{fig:534}a and \ref{fig:534}b, respectively, for the second-order CSMFEs developed in this work and the second-order triangular element with U/P mixed formulation. Regarding the contour plots, it is observed that the values of the compressive stress for the two meshes are $625$~MPa and $620$~MPa. Clearly, the two values nearly coincide.

\begin{table}
\centering
\caption{The number of elements, the number of degrees of freedom, and the displacement of point  $A$  for the meshes of second-order compatible-strain elements used to model the inhomogeneous compression test. The displacement value at point  $A$  corresponds to  $p = 600$ MPa.}
\label{tab:5.2}
\renewcommand{\arraystretch}{1.5}
\renewcommand{\tabcolsep}{0.2cm}
\begin{tabular}{c c c c}
\hline
Mesh & No. Elements & DOFs  & Displ. A (mm) \\
\hline
1    &  97          &  2958 & 6.523         \\
2    & 144          &  4348 & 6.524         \\
3    & 186          &  5590 & 6.519         \\
4    & 246          &  7362 & 6.518         \\
5    & 268          &  8014 & 6.518         \\
6    & 287          &  8570 & 6.513         \\
7    & 313          &  9334 & 6.503         \\
8    & 356          & 10598 & 6.504         \\
9    & 431          & 12800 & 6.498         \\
10   & 714          & 21122 & 6.493         \\
\hline
\end{tabular}
\end{table}

\subsection{Cook's membrane problem}
\label{sec:5.4}

\begin{figure}
\centering
\includegraphics[width=0.30\textwidth]{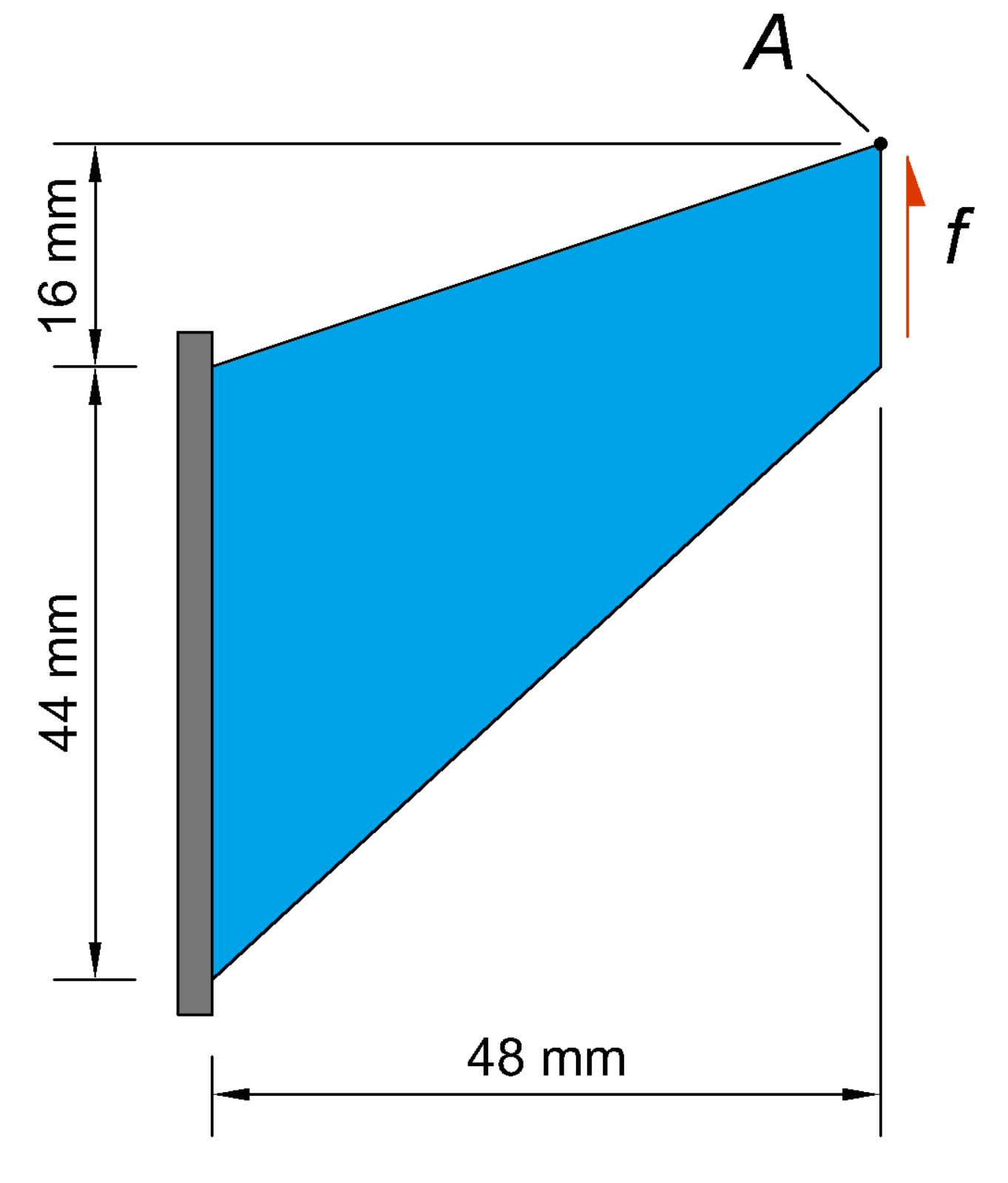}
\vskip 0.1in
\caption{Geometry and boundary conditions for Cook's membrane problem: The left edge is fixed, while the right edge is subject to a uniform shearing load $f$.}
\label{fig:541}
\end{figure}

The classical Cook's membrabe problem has been studied in many references \citep{BRI96,GLA97,REE00,ANG17,DHA22a,JAH22}. It serves as a benchmark test for bending-dominated problems in plane strain. A tapered panel clamped at one end is subjected to a uniform shearing load at the opposite end. The geometry of the problem is depicted in Figure \ref{fig:541}. The first material model in \eqref{eqn:5.1}, with parameters $\mu = 80.194$ MPa and $\kappa = 400,889.8$ MPa, is used for the simulations. It is known that certain first-order elements, such as Q1 and Q1/P0, as well as some enhanced strain-based elements, such as Q1/ES4, exhibit stiff behavior when solving Cook’s membrane problem \citep{GLA97,REE00}.

\begin{figure}
\centering
\includegraphics[width=0.60\textwidth]{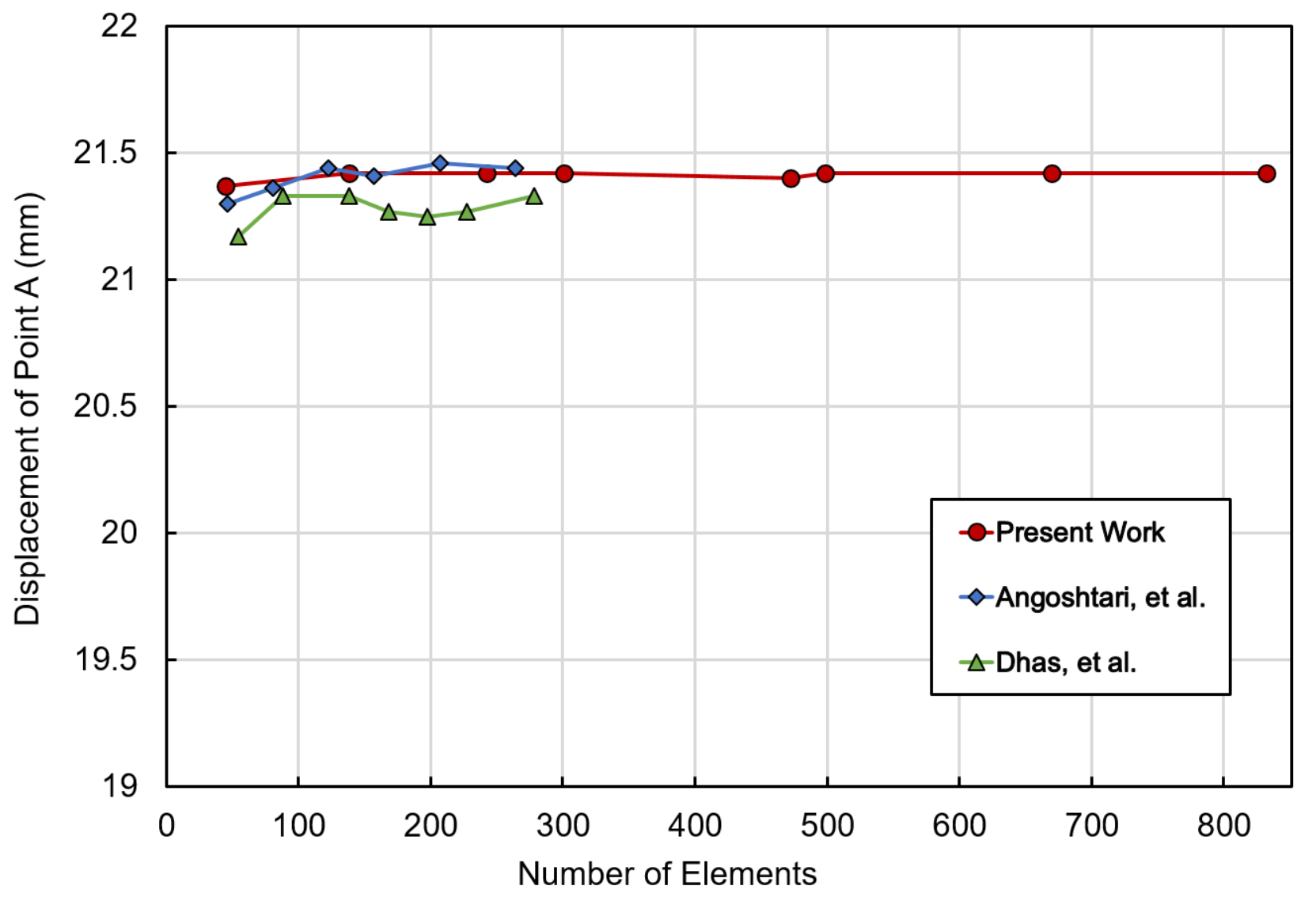}
\vskip 0.1in
\caption{Comparison of the vertical displacement of point $A$ for Cook's membrane problem with \cite{ANG17,DHA22a} using different number of elements and $f=32$ MPa.}
\label{fig:542}
\end{figure}

\begin{figure}
\centering
\includegraphics[width=0.70\textwidth]{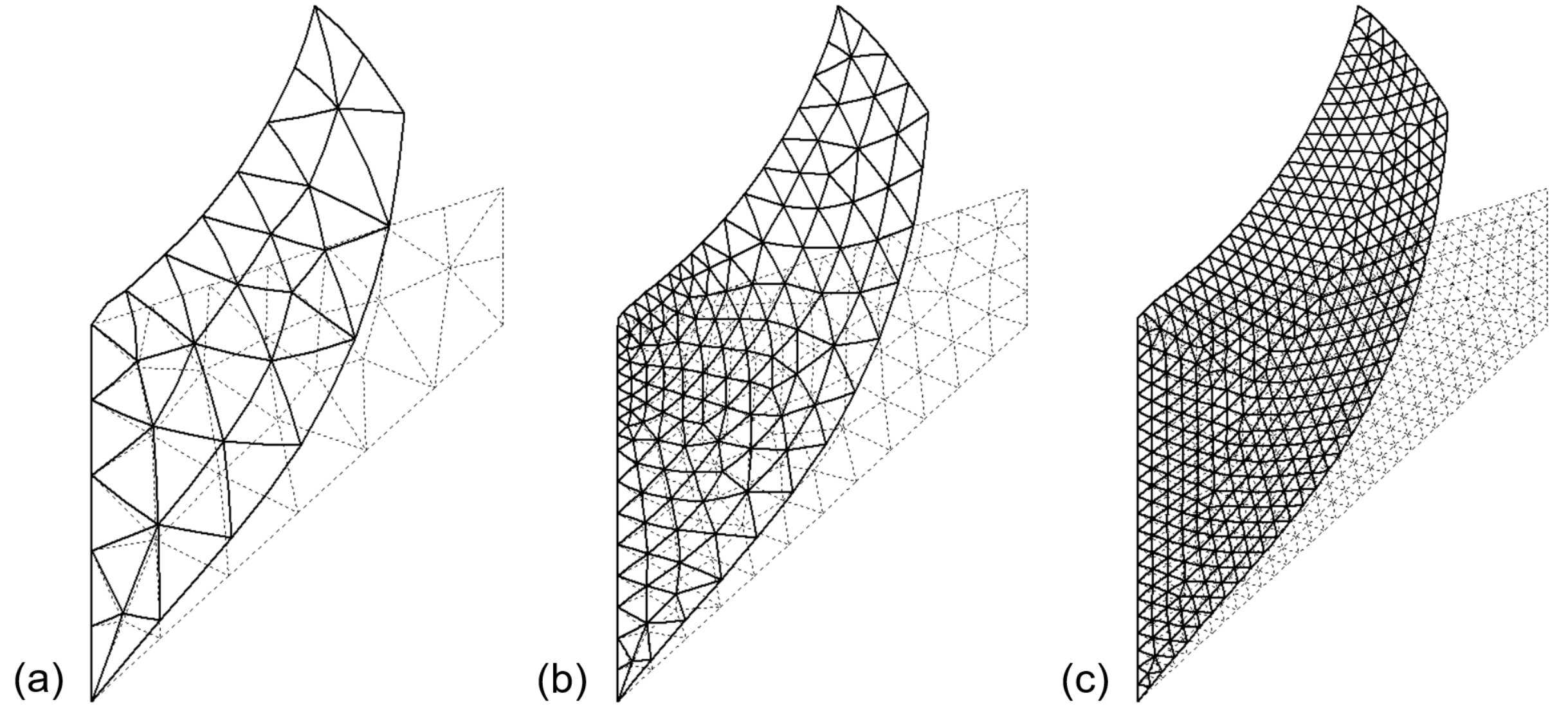}
\vskip 0.1in
\caption{Deformed configuration of Cook's membrane problem using: (a) $44$, (b) $242$ and (c) $832$ second-order compatible strain elements. The shearing load is $f=32$ MPa for all meshes.}
\label{fig:543}
\end{figure}

\begin{figure}
\centering
\includegraphics[width=0.45\textwidth]{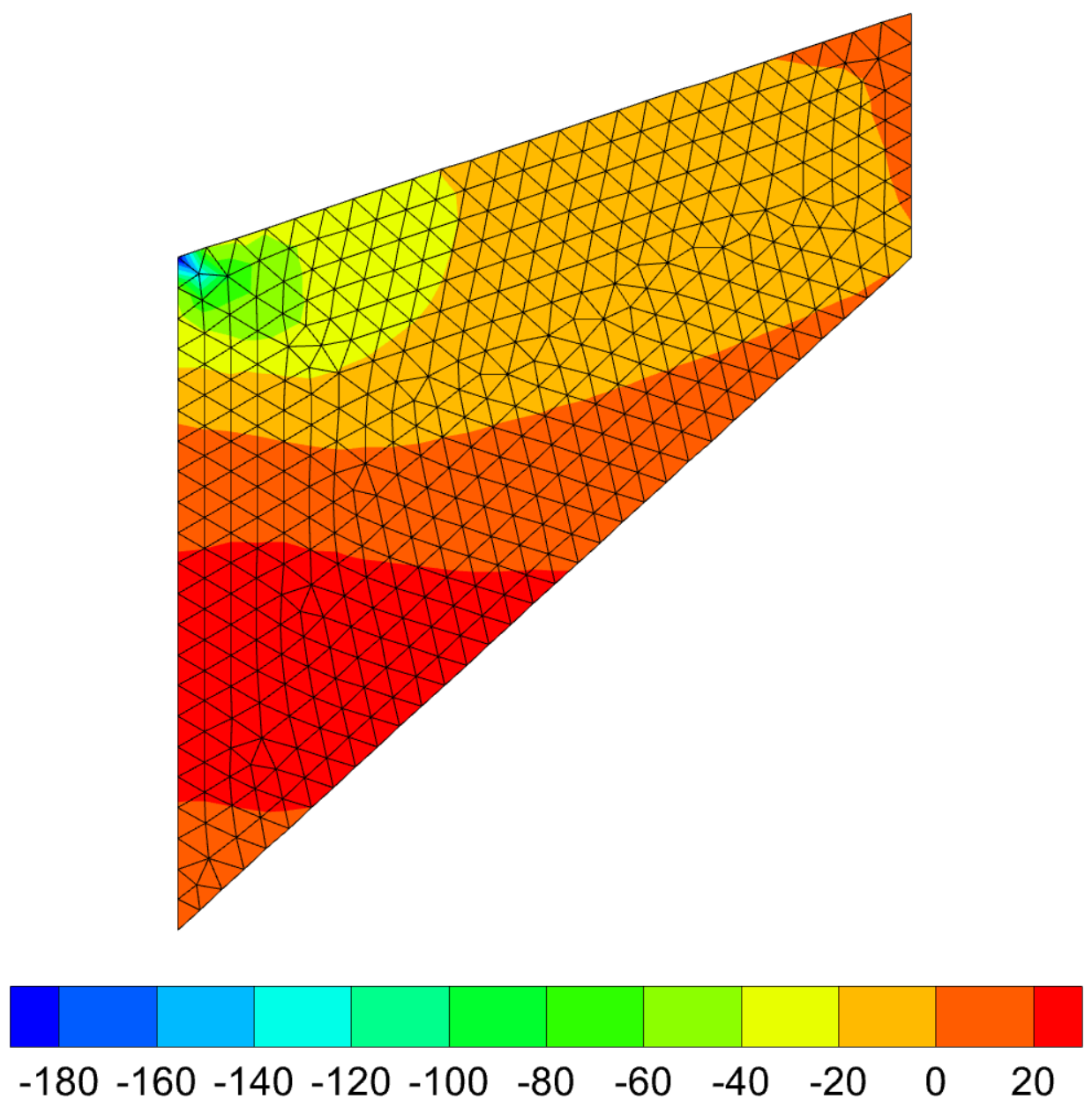}
\vskip 0.1in
\caption{Contour of the Kirchhoff stress component $\tau_{11}$ (MPa) for the Cook's membrane problem using a mesh consisting of 832 second-order compatible strain elements developed in this work.}\label{fig:544}
\end{figure}

In \citep{JAH22}, the domain was discretized using three meshes that consist of $44$, $138$ and $242$ first-order compatible strain elements. Here, we consider several meshes. The number of elements and the number of degrees of freedom for these meshes are given in Table \ref{tab:5.3}. Assuming the free end is subjected to a uniform shearing load of $f=32$ MPa, which is applied in $1000$ load steps, the vertical displacement of point $A$ computed using the meshes in this table is compared in Figure \ref{fig:542} with the results from the first-order elements of \citet{DHA22a} and the second-order elements of \citep{ANG17}. The displacement of point $A$ is also provided in Table \ref{tab:5.3} for the various meshes used in our simulations. The deformed and undeformed configurations for selected meshes are shown in Figure~\ref{fig:543}. From Figure~\ref{fig:542}, it can be observed that the displacement approximated by $832$ second-order elements developed in this work is $21.42$ mm. This value is comparable with the displacement of $21.44$ mm obtained by \citet{ANG17} using an element of type $\mathrm{H}2\mathrm{c}2\mathrm{d}\bar{2}$. The contour for the stress component $\tau_{11}$ is shown in Figure \ref{fig:544}.

\begin{table}
\centering
\caption{The number of elements, the number of degrees of freedom, and the displacement of point $A$ for the meshes of second-order compatible elements used to model Cook's membrane problem. The value of displacement at point $A$ corresponds to a uniform shearing load of $f=32$ MPa.}
\label{tab:5.3}
\renewcommand{\arraystretch}{1.5}
\renewcommand{\tabcolsep}{0.2cm}
\begin{tabular}{c c c c}
\hline
Mesh & No. Elements & DOFs  & Displ. A (mm) \\
\hline
1    & 44           &  1396 & 21.37         \\
2    & 138          &  4206 & 21.42         \\
3    & 242          &  7290 & 21.42         \\
4    & 300          &  9000 & 21.42         \\
5    & 472          & 14096 & 21.40         \\
6    & 498          & 14850 & 21.42         \\
7    & 669          & 19874 & 21.42         \\
8    & 832          & 24656 & 21.42         \\
\hline
\end{tabular}
\end{table}

\subsection{Rubber sealing problem}
\label{sec:5.5}

This problem involves the compression of a rubber sealing and was studied in \citep{BRI98} to provide a posteriori error estimation for finite element computations in finite elasticity. The Q2/P1 element and the Ogden material model were used, with the applied deformation limited to $2.0$ mm. Solving this problem becomes challenging when using the first-order elements, such as Q1 and Q1/P0. Using meshes of Q1/ET4 elements and the compressible neo-Hookean material model, it was demonstrated that displacements of up to $1.8$~mm can be modeled \citep{JAH22}. The same problem was solved in \citep{ANG17,JAH22} to study the behavior of the first and second-order CSMFEs for a displacement of $2.2$ mm.

\begin{figure}
\centering
\includegraphics[width=0.50\textwidth]{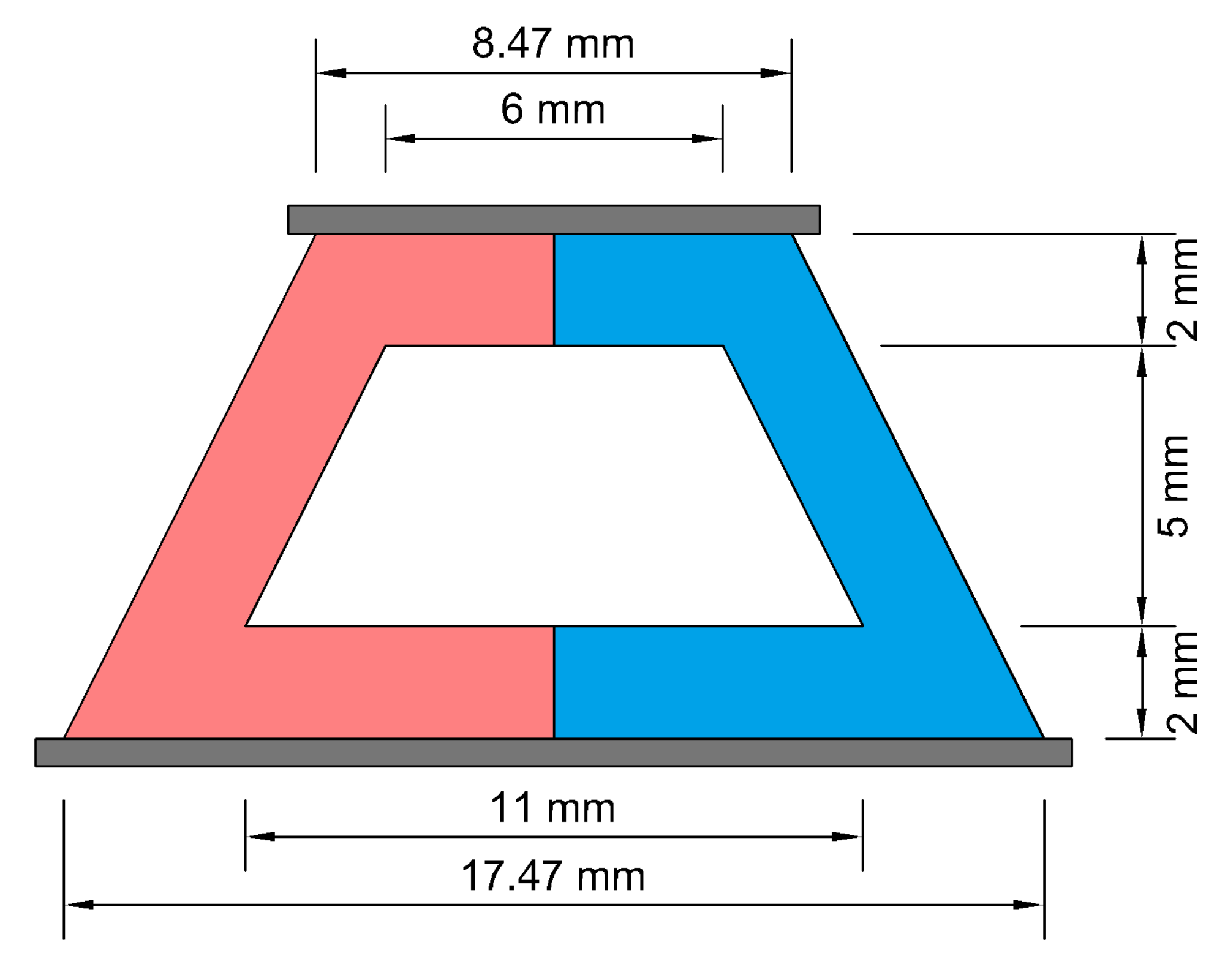}
\vskip 0.1in
\caption{Geometry and boundary conditions for the rubber sealing problem: The bottom edge is fixed, while a specified displacement of $2.2$ mm is imposed on the top edge. Due to symmetry, only the right half of the model is considered.}
\label{fig:551}
\end{figure}

\begin{figure}
\centering
\includegraphics[width=0.6\textwidth]{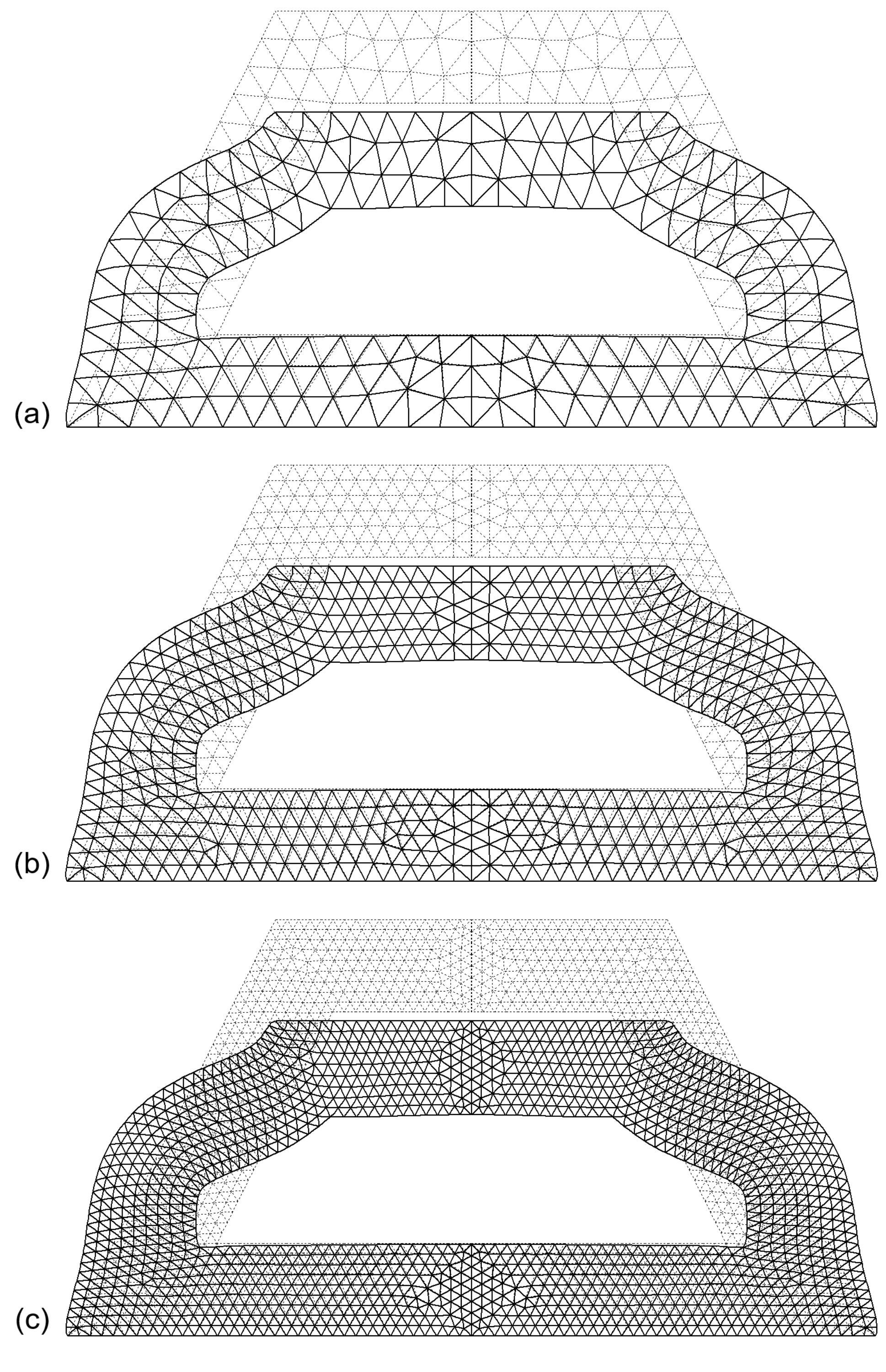}
\vskip 0.1in
\caption{Deformed configuration for the meshes used to discretize the rubber sealing problem: (a) $166$, (b) $546$ and (c) $1326$ second-order compatible strain elements.}
\label{fig:552}
\end{figure}

\begin{figure}
\centering
\includegraphics[width=0.6\textwidth]{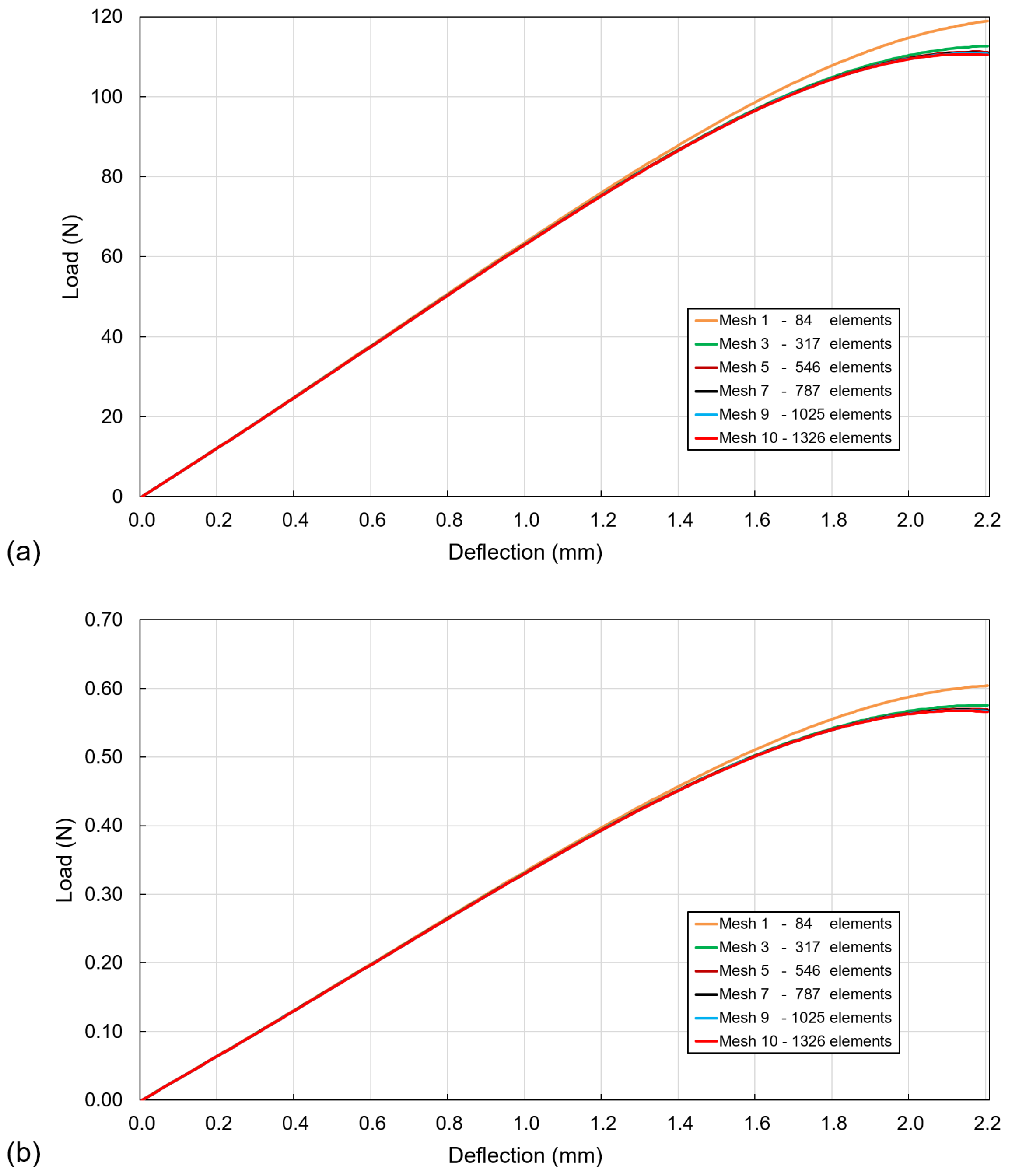}
\vskip 0.1in
\caption{Load-deflection curves for the rubber sealing problem using (a) the first and (b) the third material models, and different meshes of second-order compatible strain elements.}
\label{fig:553}
\end{figure}

\begin{figure}
\centering
\includegraphics[width=0.60\textwidth]{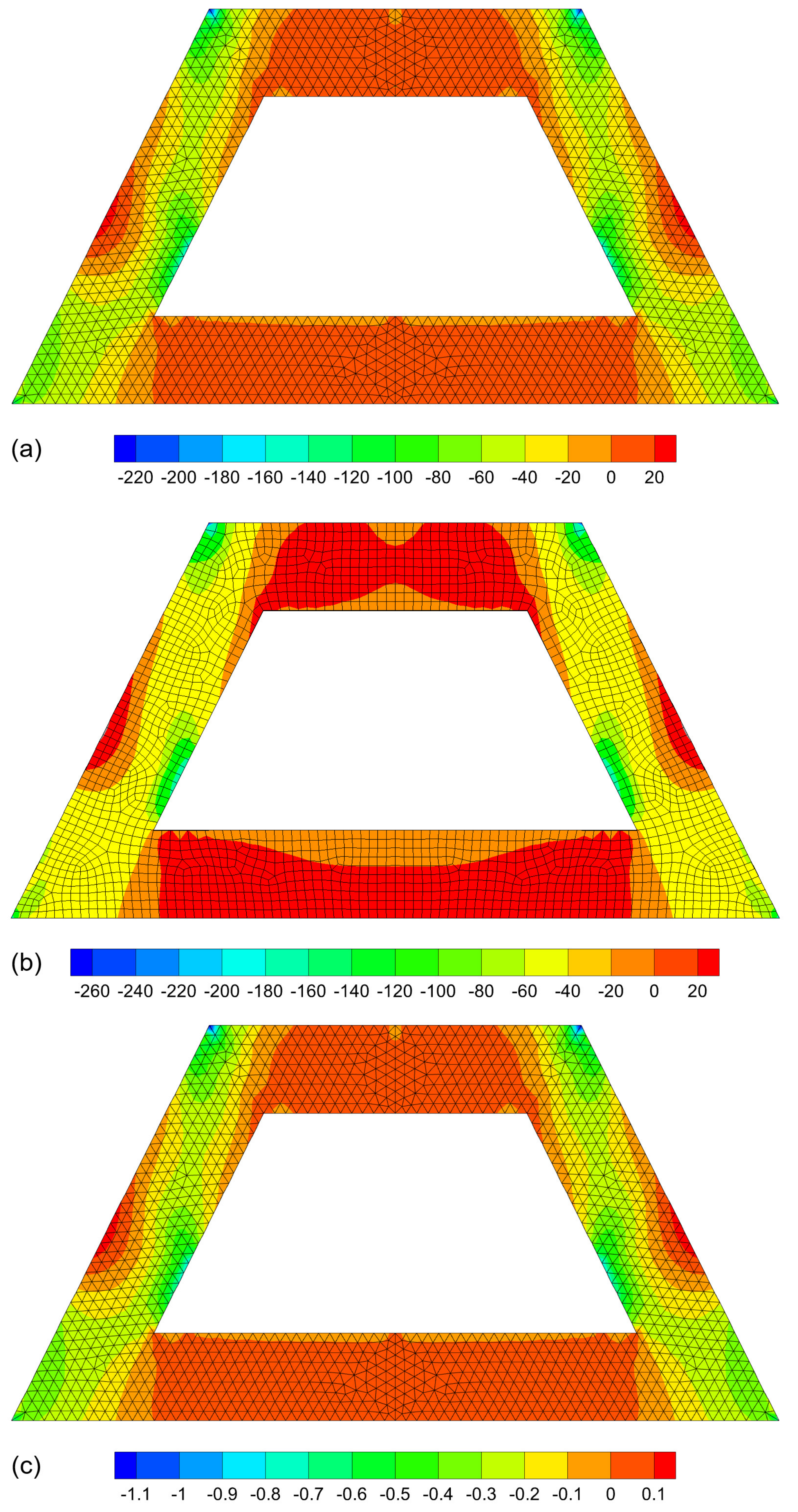}
\caption{Contours of the Kirchhoff stress component $\tau_{22}$ (MPa) for the rubber sealing problem; (a) $1326$ second-order CSMFEs developed in this work and the first material type are used to model the problem for an applied displacement of $2.2$~mm; (b) $981$ Q2/P1 elements and the first material type are applied to model the problem for a specified displacements of $2.07$~mm; (c) $1326$ CSMFEs are used to model the problem for an applied displacement of $2.2$~mm using the third material type.}
\label{fig:554}
\end{figure}

The geometry and boundary conditions of the problem are shown in Figure~\ref{fig:551}. The first and third material models in \eqref{eqn:5.1} are used for simulations. The parameters $\mu=80.194$ and $\kappa=400,889.8$~MPa are applied for the first material, while for the third material model the bulk modulus is set to $\kappa=1000$~MPa and the following parameters are used with $m=3$:
\begin{alignat*}{3}\label{eqn:5.4}
\alpha_{1}&=1.3\,,&\quad
\alpha_{2}&=5\,,&\quad
\alpha_{3}&=-2\,,\\
\mu_{1}&=0.63\text{ MPa}\,,&\quad
\mu_{2}&=0.0012\text{ MPa}\,,&\quad
\mu_{3}&=-0.01\text{ MPa}\,.
\end{alignat*}
These material parameters are the same as those used in the original work of \citep{BRI98}. A downward displacement of $2.2$~mm is imposed on the top edge, while the bottom edge remains fixed. The specified displacement is applied in $220$ load steps. Due to symmetry, only the right half of the model is considered. Meshes consisting of $84$, $166$ and $317$ first-order CSMFEs were used in \citep{JAH22} to solve the same problem. It was shown that the displacements up to of $2.2$~mm can be modeled using all three meshes. Here, we consider several meshes of second-order CSMFEs to study the behavior of the new element for this problem. The number of elements, the number of degrees of freedom and the final value of the load applied to the top edge for these meshes are given in Table \ref{tab:5.4}. The deformed and undeformed configurations for the selected meshes are shown in Figure~\ref{fig:552}. The load applied to the top edge is plotted against the displacement of the same edge for various meshes and the two material types in Figure~\ref{fig:553}. It is clear that the results converge toward the final solution as the mesh is refined. It should be mentioned that using meshes of Q2/P1 element it is only possible to apply a maximum displacement of $2.07$~mm to the top edge of the specimen, while a maximum displacement of $2.2$~mm has been applied to all meshes of the CSMFEs developed in this work. This point becomes more important when it is noted that the Q2/P1 element is developed specifically for modeling incompressible media. The pressure is considered as an additional degree of freedom and it is interpolated linearly over the domain of the element. However, no degrees of freedom are considered for the pressure in our formulation. Figures \ref{fig:554}a and \ref{fig:554}b present the contour plots of the Kirchhoff stress component $\tau_{22}$ for the meshes of $1326$ CSMFEs and $981$ Q2/P1 elements using the first material model. The applied displacements for these meshes were $2.2$~mm and $2.07$~mm, respectively. Regarding the contours, it is observed that the maximum values of the compressive stress for the two meshes are $235$~MPa and $261$~MPa. While the applied displacement for the Q2/P1 mesh was smaller than the mesh of CSMFEs, the maximum compressive stress was higher for that mesh. The contour plots of the Kirchhoff stress component $\tau_{22}$ for a mesh of $1326$ CSMFEs and an applied displacement of $2.2$~mm using the third material model is shown in Figure \ref{fig:554}c.

\begin{table}
\centering
\caption{The number of elements, the number of degrees of freedom, and the final value of the load applied to the top edge at the target displacement of $2.2$~mm for the meshes of second-order compatible strain elements, and the first and third material types used to model the rubber sealing problem.}
\label{tab:5.4}
\renewcommand{\arraystretch}{1.5}
\renewcommand{\tabcolsep}{0.2cm}
\begin{tabular}{c c c c c}
\hline
\multirow{2}{*}{Mesh} & \multirow{2}{*}{No. Elements} & \multirow{2}{*}{DOFs} & \multicolumn{2}{c}{Load (N)} \\
\cline{4-5}
                      &                               &                       & Material 1 &    Material 3   \\
\hline
          1           &            84                 &         2678          &  118.89    &      0.6042     \\
          2           &            166                &         5154          &  115.00    &      0.5862     \\
          3           &            317                &         9662          &  112.64    &      0.5755     \\
          4           &            489                &         14762         &  111.61    &      0.5711     \\
          5           &            546                &         16442         &  111.15    &      0.5692     \\
          6           &            626                &         18806         &  111.11    &      0.5690     \\
          7           &            787                &         23536         &  110.84    &      0.5678     \\
          8           &            851                &         25436         &  110.82    &      0.5677     \\
          9           &            1025               &         30526         &  110.65    &      0.5669     \\
          10          &            1326               &         39360         &  110.43    &      0.5658     \\
\hline
\end{tabular}
\end{table}

\subsection{Tension of a perforated block}
\label{sec:5.6}

Multiple versions of this problem have been solved in the literature to investigate the behavior of the first and second-order compatible strain elements under large stretches. In order to study the behavior of a first-order element, \citet{DHA22a} used a block of $2\times2$~mm subject to a tensile deformation of $100\%$. They employed the material parameters $\mu=10$ and $\kappa=1000$ MPa in their simulations. In their work, the edges subject to tensile deformation were free to move in the transverse direction. \citet{JAH22} studied the problem by considering a block of $1\times1$~mm and the material parameters $\mu=80.192$ and $\kappa=400,933.33$ MPa. Similar to the work of \citet{DHA22a}, the edges subject to tensile deformation were free to move in the transverse direction. He used different meshes of the first-order compatible strain element. The applied deformation was, however, limited to $50\%$ of the height of the block. \citet{SHO18} used the same block as employed in the present study, but in the context of the incompressible elasticity. They considered the pressure as an independent field in the formulation of their element. The tensile deformation that they applied to the block was $300\%$ of its initial dimension.

\begin{figure}
\centering
\includegraphics[width=0.35\textwidth]{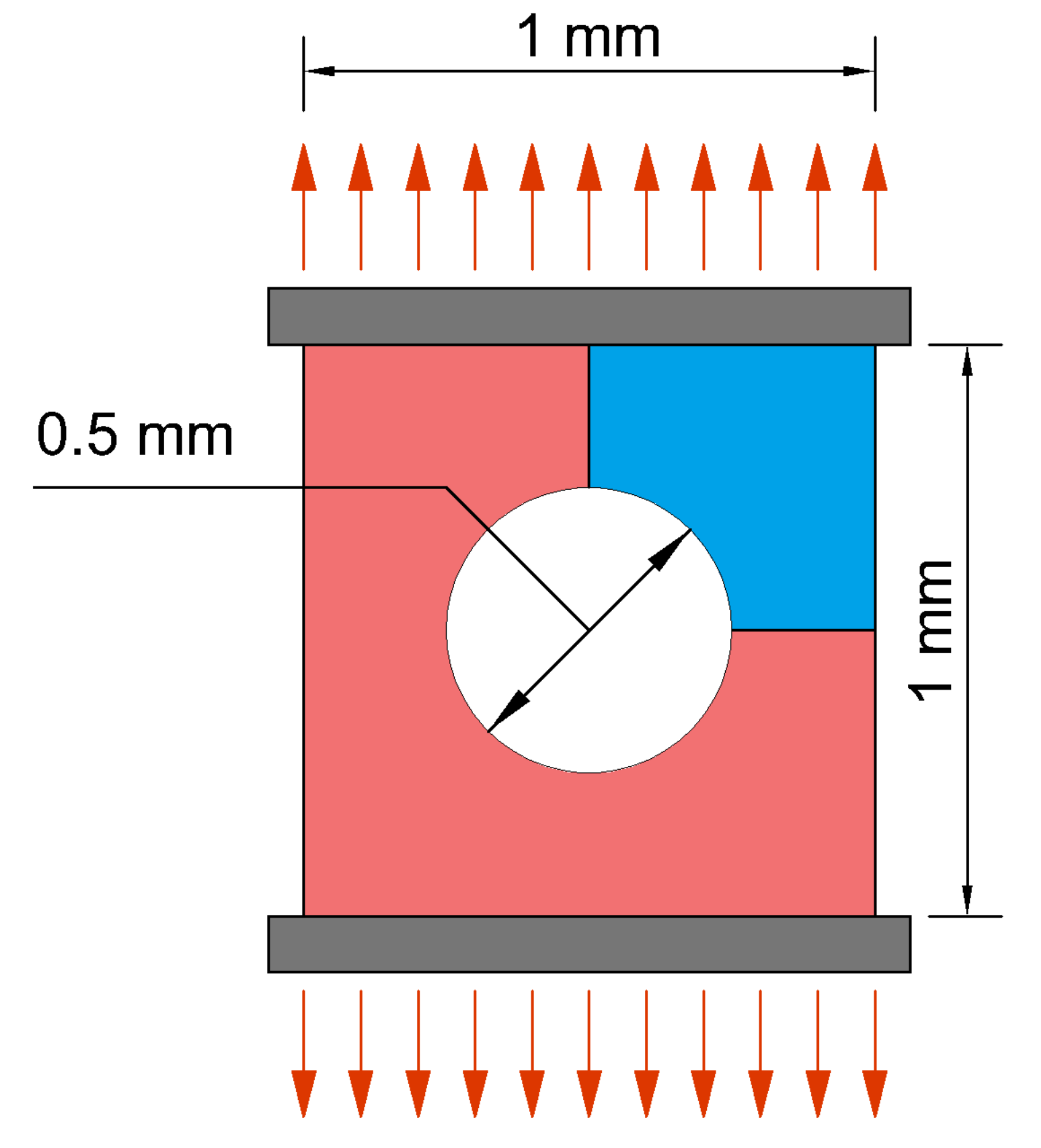}
\vskip 0.1in
\caption{Geometry and boundary conditions for the perforated block under tension: A tensile deformation equal to $300\%$ of the height of the block is applied to the top and bottom edges, which are restrained in the horizontal direction. Due to symmetry in both the horizontal and vertical directions, only the upper right quarter of the block is modeled.}
\label{fig:561}
\end{figure}

In the present work, a block with a square cross section of $1\times1$~mm and a central hole is subjected to a tensile deformation equal to $300\%$ of the height of the block at the top and bottom edges. The specified displacement is applied in $1500$ load steps. The top and bottom edges are restrained in the horizontal direction. The diameter of the hole is $0.50$~mm. The geometry and boundary conditions of the problem are shown in Figure~\ref{fig:561}. The first material model in \eqref{eqn:5.1}, with $\mu = 10$ and $\kappa = 1000$ MPa, is used for simulations. Due to symmetry in both the horizontal and vertical directions, only a quarter of the block is modeled. We consider several meshes in order to study the problem. The number of elements, the number of degrees of freedom and the final value of the load applied to the top edge for these meshes are provided in Table \ref{tab:5.5}. In order to compare the performance of the second-order CSMFEs developed in this work with the first-order CSMFEs developed in \citep{JAH22}, the problem is also solved with the material parameters $\mu=80.192$ and $\kappa=400,933.33$~MPa assigned to the first material model in \eqref{eqn:5.1}. The top and bottom edges for this second model are free to move in the transverse direction. The maximum tensile deformation that can be applied to the model using the first-order CSMFEs is limited to $175$\% of the height of the block. A total of $875$ load steps is used to apply the deformation. Beyond this value, the convergence is not possible using the first-order CSMFEs. The performance of the second-order CSMFEs is compared with that of the first-order CSMFEs through the load-deflection curves and contour plots of the Kirchhoff stress tensor.

\begin{figure}
\centering
\includegraphics[width=0.7\textwidth]{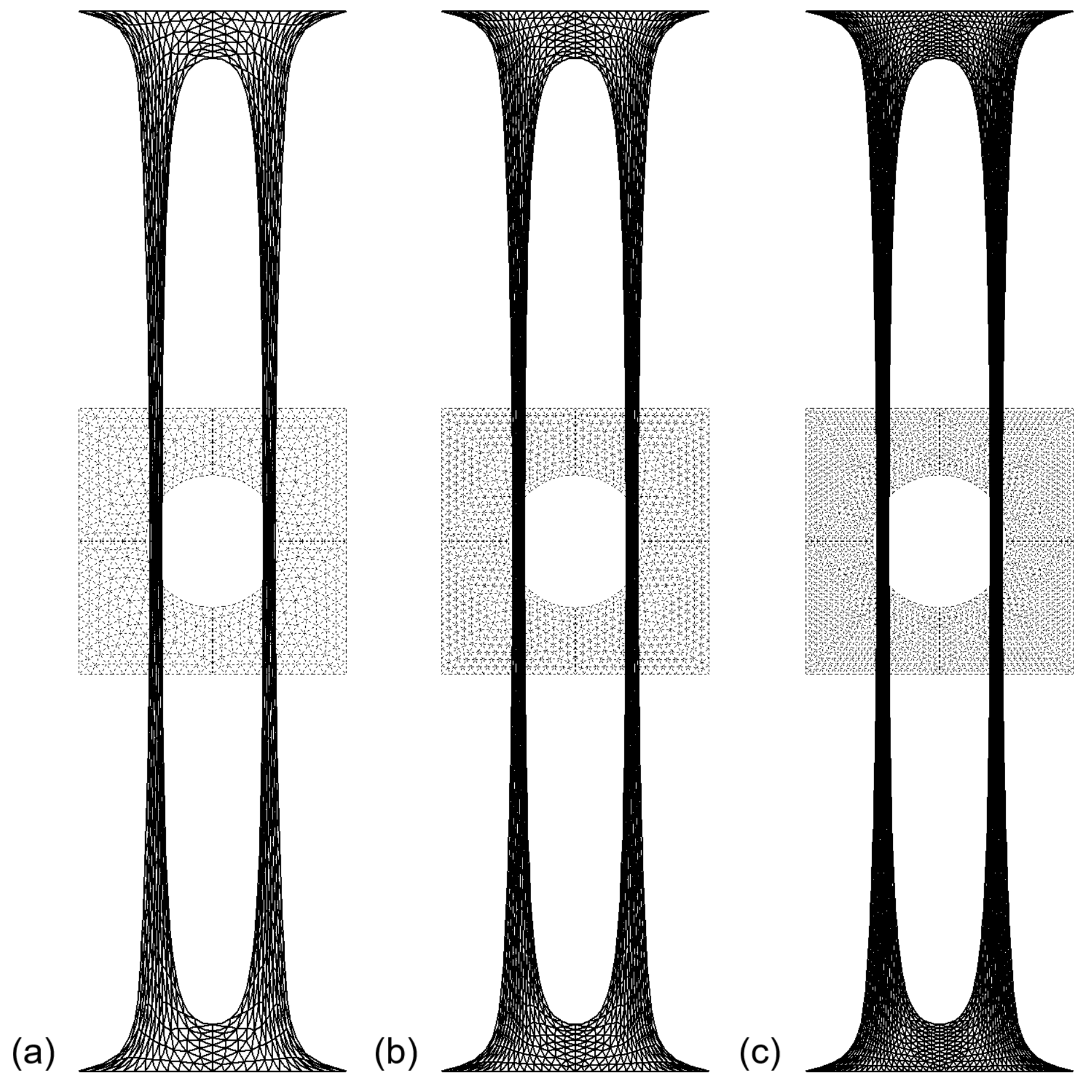}
\vskip 0.1in
\caption{Deformed configurations of the perforated block subject to a tensile deformation equal to $300$\% of the height of the block using meshes consisting of (a) $258$, (b) $457$ and (c) $845$ second-order CSMFEs developed in this work; material parameters are $\mu=10$ and $\kappa=1000$~MPa, and the top and bottom edges are restrained in the horizontal direction.}
\label{fig:562}
\end{figure}

\begin{figure}
\centering
\includegraphics[width=0.5\textwidth]{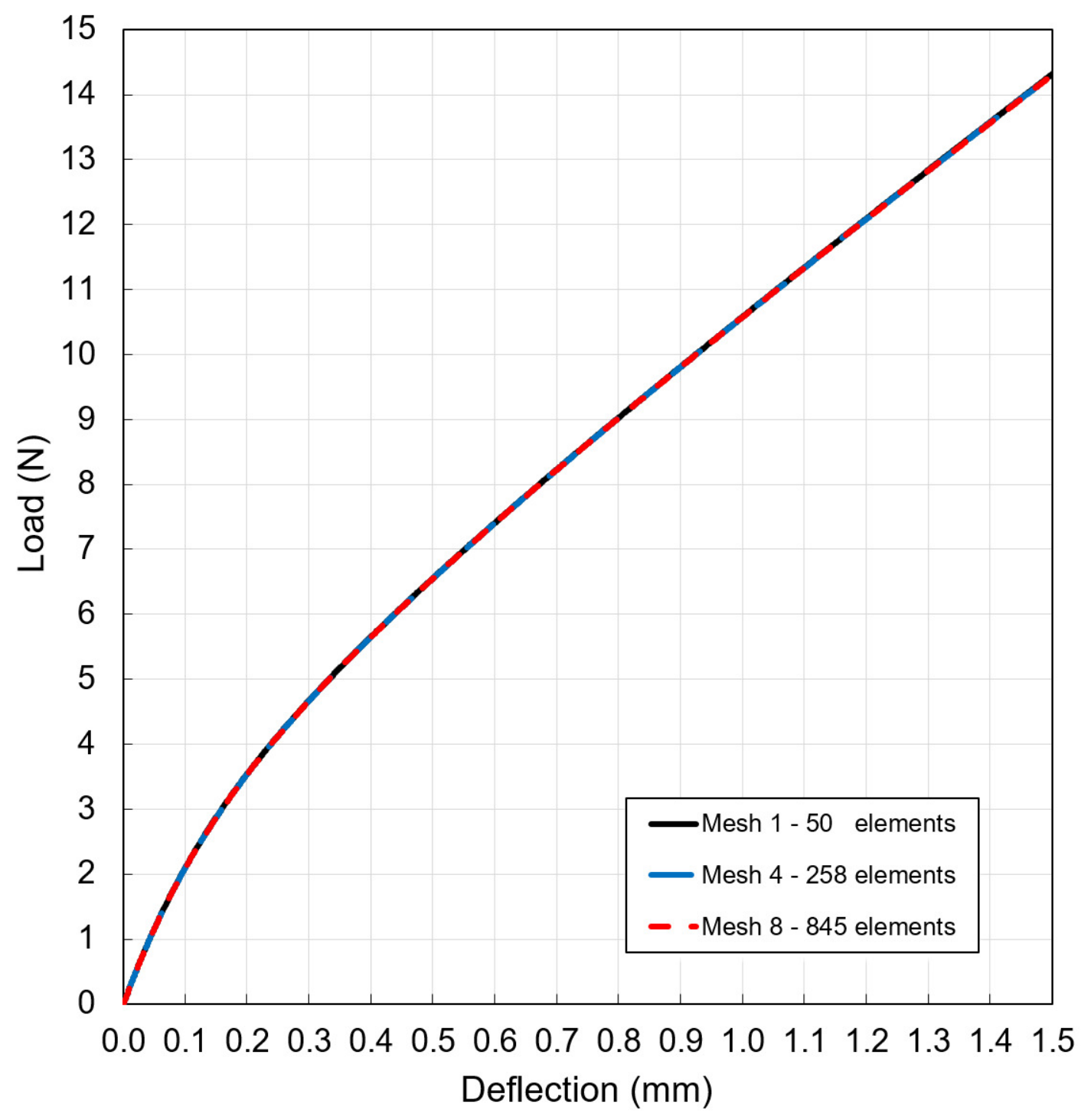}
\vskip 0.1in
\caption{Load-deflection curves for the perforated block subject to a tensile deformation equal to $300$\% of the height of the block using different meshes of second-order CSMFEs developed in this work; material parameters are $\mu=10$ and $\kappa=1000$~MPa, and the top and bottom edges are restrained in the horizontal direction.}
\label{fig:563}
\end{figure}

\begin{figure}
\centering
\includegraphics[width=0.5\textwidth]{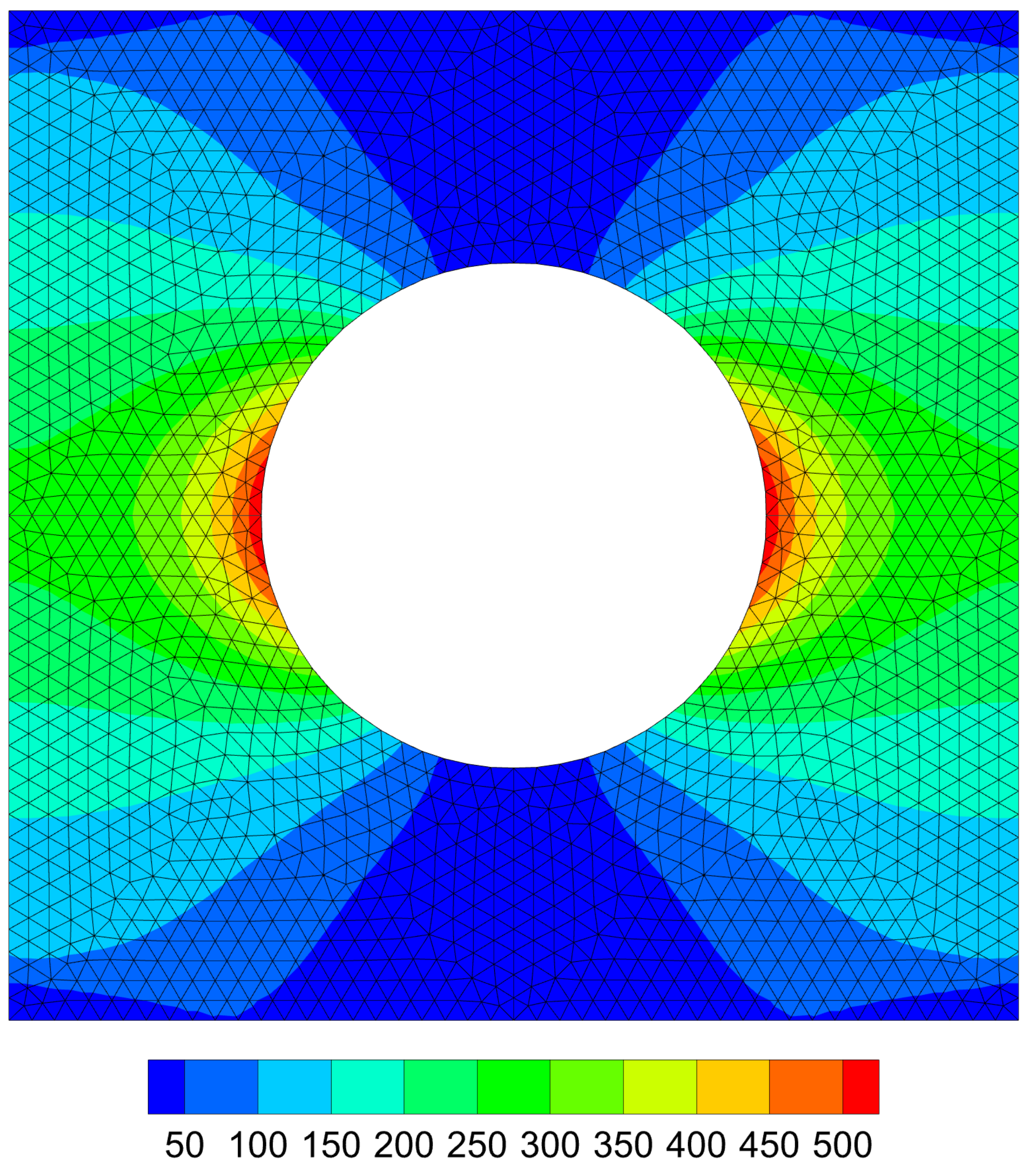}
\vskip 0.1in
\caption{Contour of the Kirchhoff stress component $\tau_{22}$~(MPa) for the perforated block subject to a tensile deformation equal to $300$\% of the height of the block using a mesh consisting of $845$ second-order CSMFEs developed in this work; material parameters are $\mu=10$ and $\kappa=1000$~MPa, and the top and bottom edges are restrained in the horizontal direction.}
\label{fig:564}
\end{figure}

\begin{figure}
\centering
\includegraphics[width=0.7\textwidth]{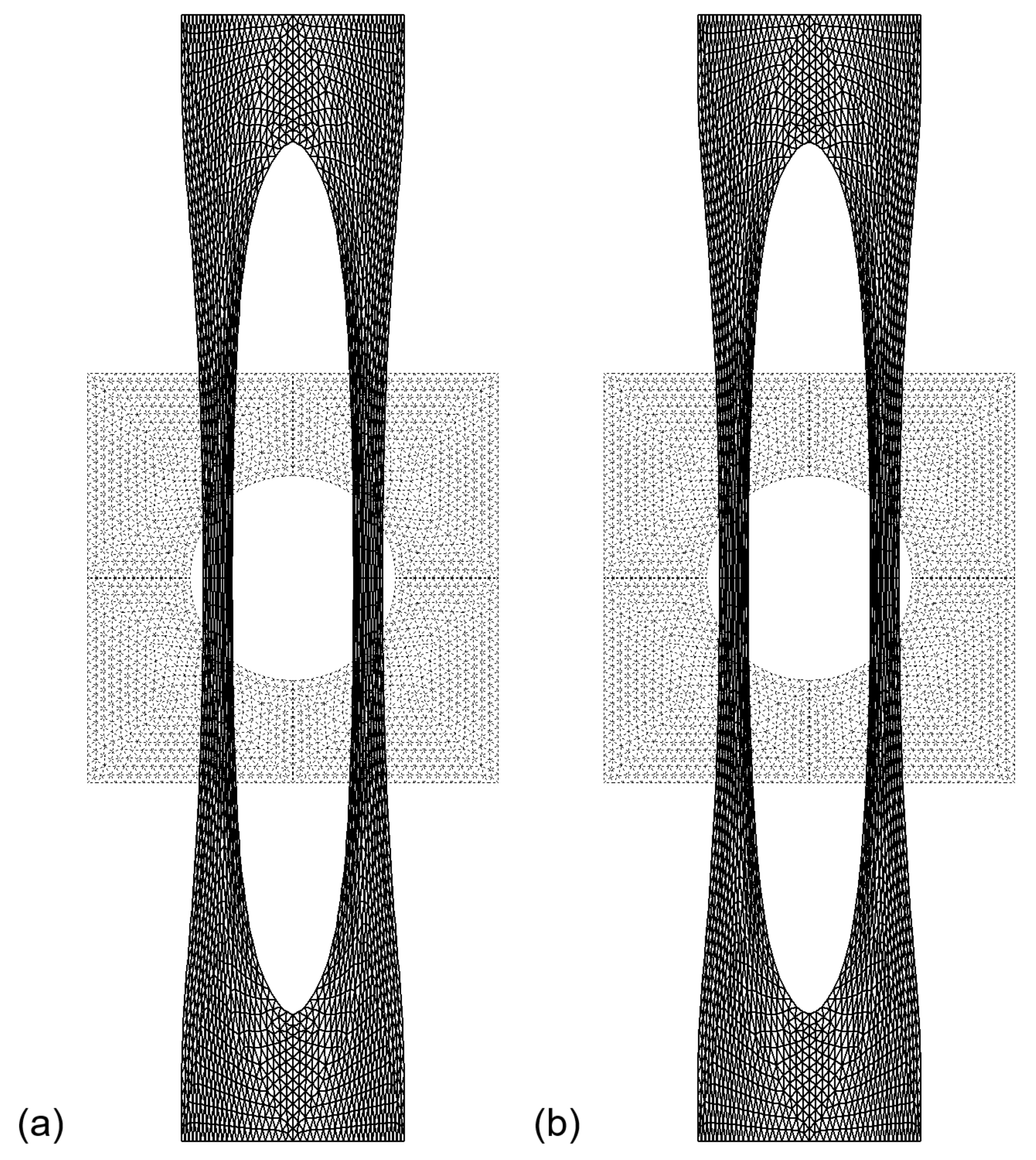}
\vskip 0.1in
\caption{Deformed configurations of the perforated block subject to a tensile deformation equal to $175$\% of the height of the block using meshes consisting of (a) $845$ second-order CSMFEs developed in this work and (b) $845$ first-order CSMFEs developed in \citep{JAH22}; material parameters are $\mu=80.192$ and $\kappa=400,933.33$~MPa, and the top and bottom edges are free to move in the transverse direction.}
\label{fig:565}
\end{figure}

\begin{figure}
\centering
\includegraphics[width=0.5\textwidth]{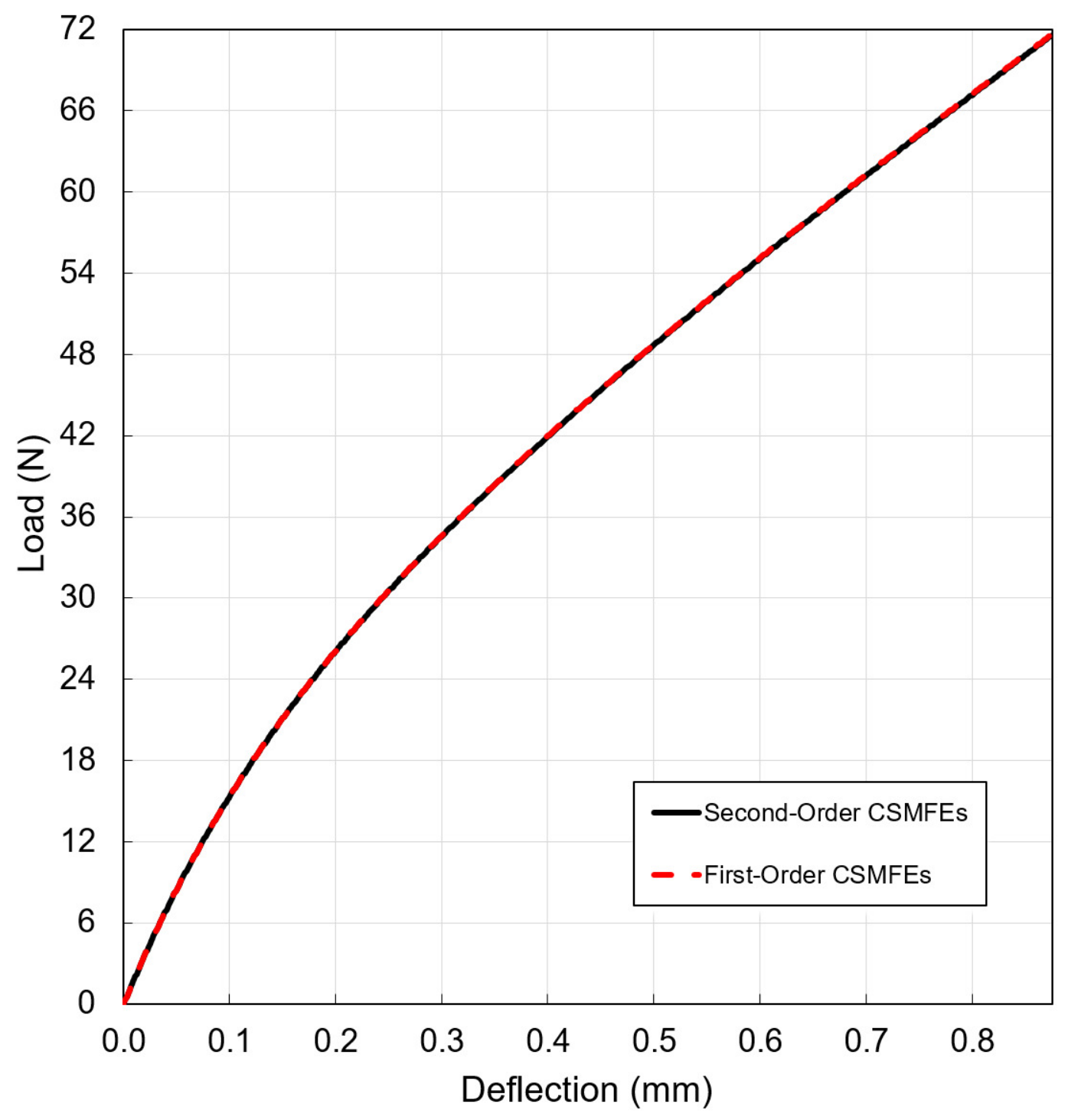}
\vskip 0.1in
\caption{Load-deflection curves for the perforated block subject to a tensile deformation equal to $175$\% of the height of the block using $845$ second-order CSMFEs developed in this work and $845$ first-order CSMFEs developed in \citep{JAH22}; material parameters are $\mu=80.192$ and $\kappa=400,933.33$~MPa, and the top and bottom edges are free to move in the transverse direction.}
\label{fig:566}
\end{figure}

\begin{figure}
\centering
\includegraphics[width=\textwidth]{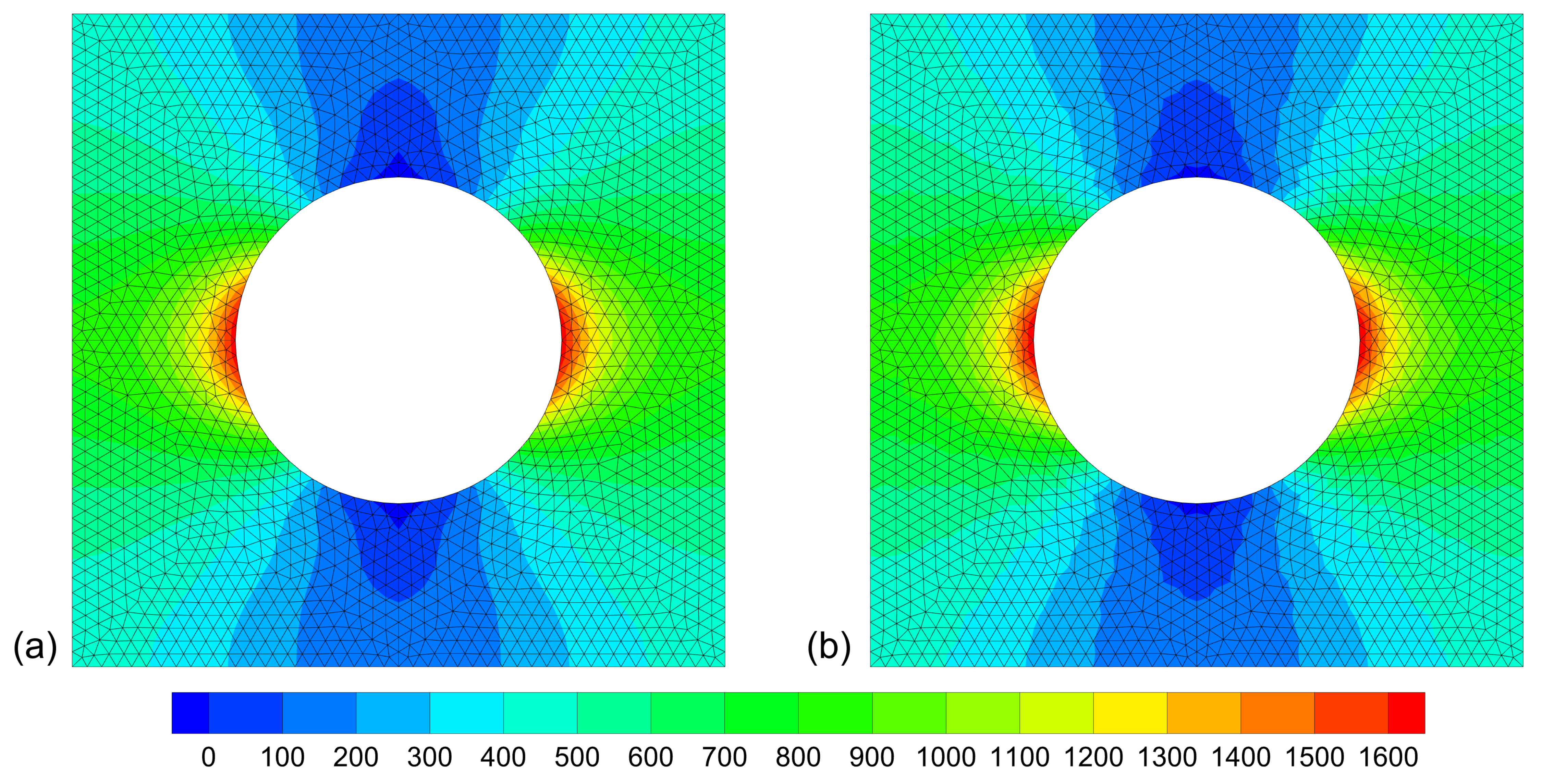}
\vskip 0.1in
\caption{Contours of the Kirchhoff stress component $\tau_{22}$~(MPa) for the perforated block subject to a tensile deformation equal to $175$\% of the height of the block using (a) $845$ second-order CSMFEs developed in this work and (b) $845$ first-order CSMFEs developed in \citep{JAH22}; material parameters are $\mu=80.192$ and $\kappa=400,933.33$~MPa, and the top and bottom edges are restrained in the horizontal direction.}
\label{fig:567}
\end{figure}

The deformed and undeformed configurations of the block subject to a tensile deformation equal to $300$\% of its height are shown in Figure~\ref{fig:562} for the meshes consisting of $258$, $457$ and $845$ second-order CSMFEs. The load-displacement curves resulting from this value of tensile deformation applied to top edge (for a quarter of the block) are plotted in Figure~\ref{fig:563} for the meshes of $50$, $258$ and $845$ second-order CSMFEs. It is observed that the curves corresponding to the meshes with $258$ and $845$ elements coincide with each other (in fact, all meshes in rows $2$ to $8$ of Table \ref{tab:5.5}), indicating excellent convergence. The contour plot of the Kirchhoff stress component $\tau_{22}$ is shown in Figure~\ref{fig:564} for the applied tensile deformation of $300$\%. Figure~\ref{fig:565} presents the deformed and undeformed configurations of the block for an applied tensile deformation equal to $175$\% of the height of the block using meshes of $845$ second-order and first-order CSMFEs. The load-displacement curves for the block subject to this tensile deformation are shown in Figure~\ref{fig:566} using meshes of second-order and first-order CSMFEs. Clearly, the two curves coincide with each other. Figure~\ref{fig:567} presents the contour plots of the Kirchhoff stress component $\tau_{22}$ using $845$ second-order and first-order CSMFEs. It is observed that the minimum and maximum contour values for both meshes are the same. However, the second-order CSMFEs yields a smoother contour plot compared with that of the first-order CSMFEs.
\begin{remark}\label{rem:5.3}
It is important to note that in this work the deformation of $300\%$ is applied to the block in the context of compressible elasticity. This is in contrast to \citep{SHO18}, where the same deformation is applied in the context of incompressible elasticity. This proves the stability of the current element for such large deformations.
\end{remark}

\begin{table}
\centering
\caption{The number of elements, the number of degrees of freedom, and the final value of the load applied to the top edge at a target displacement equal to $300\%$ of the height of the block for the meshes of second-order CSMFEs developed in this work.}
\label{tab:5.5}
\renewcommand{\arraystretch}{1.5}
\renewcommand{\tabcolsep}{0.2cm}
\begin{tabular}{c c c c}
\hline
Mesh & No. Elements & DOFs   & Load (N)  \\
\hline
1    & 50           &  1557  & 14.32     \\
2    & 89           &  2725  & 14.31     \\
3    & 180          &  5425  & 14.31     \\
4    & 258          &  7731  & 14.31     \\
5    & 345          & 10291  & 14.31     \\
6    & 457          & 13583  & 14.31     \\
7    & 566          & 16781  & 14.31     \\
8    & 845          & 24953  & 14.31     \\
\hline
\end{tabular}
\end{table}

\clearpage
\section{Conclusions}
\label{sec:6}

In this paper, second-order compatible-strain mixed finite elements were formulated. Exact shape functions for the displacement gradient and stress tensor were provided in the natural coordinate system. These shape functions were then transformed to the physical space using the covariant and contravariant Piola transformation. A three-field functional of Hu-Washizu type, with displacement, displacement gradient and the first Piola-Kirchhoff stress tensor as independent fields, was used to variationally derive the governing equations of compressible nonlinear elasticity. Our efficient implementation of second-order CSMFEs has led to a finite element code that is reasonably fast and stable, making it well-suited for a broad range of practical problems. The performance of the second-order CSMFEs has been evaluated through several numerical examples. The load-deflection curves and convergence analysis have demonstrated the excellent convergence properties of the second-order CSMFEs, particularly for near-incompressible solids. Our numerical results indicate that the second-order CSMFEs are free from numerical artifacts and do not exhibit the stiff behavior commonly seen in low-order and certain enhanced strain-based elements. It has been observed that even coarse meshes can approximate the correct solutions with sufficient accuracy.
In summary, the good convergence properties and numerical stability of the proposed second-order CSMFEs demonstrate their potential for solving problems prone to numerical instabilities, such as hourglassing and shear locking.
It is important to note that, compared to conventional finite elements, second-order CSMFEs require the evaluation of shape functions for the displacement gradient and stress tensor. However, since the optimal forms of the shape functions are used and their transformation to the physical space is performed with minimal calculations, this additional computation is not costly. The proposed second-order CSMFEs can be extended to incompressible solids by interpolating the pressure field as an independent variable. Similar concepts can be applied to develop first and second-order quadrilateral CSMFEs.

\bibliographystyle{abbrvnat}
\bibliography{ref}

\appendix

\section{Explicit Forms of the Second Order Shape Functions}
\label{app:A}

In this appendix, the local shape functions that interpolate the displacement gradient and the stress tensor are presented in explicit form for the reference element shown in Figure~\ref{fig:421}. In the first step, the integrals in \eqref{eqn:3.5}, \eqref{eqn:3.14}, \eqref{eqn:3.17} and \eqref{eqn:3.23} are transformed from the simplicial element $\mathscr{T}$ in the physical space (see Figure \ref{fig:21}) to the reference element $\widehat{\mathscr{T}}$ in the natural coordinate system. The local shape functions $\widehat{\mathbf{v}}$ are then computed by enforcing the conditions in \eqref{eqn:3.7} and \eqref{eqn:3.13}. The local $\mathcal{P}^{c}_{2}(T\widehat{\mathscr{T}})$ and $\mathcal{P}^{d}_{2}(T\widehat{\mathscr{T}})$ shape functions $\widehat{\bar{\mathbf{v}}}$ are obtained using \eqref{eqn:3.8}, and $\mathcal{P}^{c-}_{2}(T\widehat{\mathscr{T}})$ and $\mathcal{P}^{d-}_{2}(T\widehat{\mathscr{T}})$ shape functions are determined via \eqref{eqn:3.19}. These shape functions can be transformed back to the physical space using the covariant and contravariant Piola transformations. Having computed the local shape functions in the physical space, the global shape functions are obtained through \eqref{eqn:3.10} and \eqref{eqn:3.11}.

The shape functions that belong to the spaces $\mathcal{P}^{c}_{2}\left(T\mathscr{T}\right)$ and $\mathcal{P}^{c-}_{2}\left(T\mathscr{T}\right)$ (shape functions used to interpolate the displacement gradient) are transformed from the natural coordinate system to the physical space using the covariant Piola transformation defined as follows \citep{ROG10}:
\begin{equation}\label{eqn:A1}
\mathbf{v}\left(\mathbf{x}\right)=\mathbf{J}^{-\mathsf{T}}\widehat{\mathbf{v}}\circ\bm{\psi}^{-1}\left(\mathbf{x}\right),
\end{equation}
where $\bm{\psi}\left(\widehat{\mathbf{x}}\right):\widehat{\mathscr{T}}\rightarrow\mathscr{T}$ maps the reference element in Figure \ref{fig:421} to the simplicial element in Figure \ref{fig:21}, $\mathbf{J}=\nabla_{\widehat{x}}\bm{\psi}$ is the Jacobian of the mapping, $\mathbf{v}\left(\mathbf{x}\right)$ is the shape function in the physical space (reference configuration) and $\widehat{\mathbf{v}}\circ\bm{\psi}^{-1}\left(\mathbf{x}\right)$ is its pull-back to the natural coordinate system. On the other hand, the shape functions that belong to the spaces $\mathcal{P}^{d}_{2}\left(T\mathscr{T}\right)$ and $\mathcal{P}^{d-}_{2}\left(T\mathscr{T}\right)$ (shape functions that interpolate the stress tensor) can be obtained in the physical space via the contravariant Piola transformation \citep{ROG10,AZN22}:
\begin{equation}\label{eqn:A2}
\mathbf{v}\left(\mathbf{x}\right)=\frac{1}{\det \mathbf{J}}\,\mathbf{J}\,\widehat{\mathbf{v}}\circ\bm{\psi}^{-1}\left(\mathbf{x}\right).
\end{equation}
It is important to note that if the origins of the coordinate systems in the reference element and the element in the physical space are chosen as indicated in Figures \ref{fig:421}.a and \ref{fig:21}, then the Jacobian $\mathbf{J}$ takes the following simple form:
\begin{equation}\label{eqn:A3}
\mathbf{J}=
\begin{bmatrix}
  x_{2} & x_{3} \\
  y_{2} & y_{3}
\end{bmatrix},
\end{equation}
where $\left(x_{2},y_{2}\right)$ and $\left(x_{3},y_{3}\right)$ are the coordinates of nodes $2$ and $3$ with respect to the origin of the coordinate system in the physical space. It is clear that $\det\mathbf{J}$ is equal to twice the area of the simplicial element in this space.

In the discussion that follows, $e$ denotes one of the three edges of the simplicial element in the physical space (see Figure \ref{fig:21}). The same edge is represented by $\widehat{e}$ in the reference element shown in Figure \ref{fig:421}b. The vectors that extend from the first node to the second node of the edges $e$ and $\widehat{e}$ are denoted by $\mathbf{e}$ and $\widehat{\mathbf{e}}$, respectively. Evidently, these vectors can be expressed in the alternative forms as:
\begin{equation}\label{eqn:A4}
    \mathbf{e}=\|\mathbf{e}\|\mathbf{t}\,,\qquad
    \widehat{\mathbf{e}}=\|\widehat{\mathbf{e}}\|\,\widehat{\mathbf{t}}\,,
\end{equation}
where $\|\mathbf{e}\|$ and $\|\widehat{\mathbf{e}}\|$ are the lengths of the edges $e$ and $\widehat{e}$, and $\mathbf{t}$ and $\widehat{\mathbf{t}}$ are the unit vectors tangent to these edges.

\subsection{$\mathcal{P}^{c}_{2}(T\mathscr{T})$ shape functions}
\label{sec:A.1}

Noting that the following relationship exists between the vectors $\mathbf{e}$ and $\widehat{\mathbf{e}}$ \citep{AZN22}:
\begin{equation}\label{eqn:A5}
\mathbf{e}=\mathbf{J}\,\widehat{\mathbf{e}}\,,
\end{equation}
the unit vector $\mathbf{t}$ can be expressed in terms of $\widehat{\mathbf{t}}$ as:
\begin{equation}\label{eqn:A6}
\mathbf{t}=\frac{\|\widehat{\mathbf{e}}\|}{\|\mathbf{e}\|}\mathbf{J}\,\widehat{\mathbf{t}}\,.
\end{equation}
Observing that the base functions $s\in\mathcal{P}_{1}\left(e\right)$ and  $\widehat{s}\in\mathcal{P}_{1}\left(\widehat{e}\right)$ are related as:
\begin{equation}\label{eqn:A7}
s=\frac{\|\mathbf{e}\|}{\|\widehat{\mathbf{e}}\|}\widehat{s}\,,
\end{equation}
the integrals in \eqref{eqn:3.5} can be transformed to the natural coordinate system using \eqref{eqn:A1}, \eqref{eqn:A6} and \eqref{eqn:A7} as follows:
\begin{align}\label{eqn:A8}
    \phi^{\mathscr{T},e_{i}}_{j}\left(\widehat{\mathbf{v}}\right)&=\left[\frac{\|\mathbf{e}_{i}\|}{\|\widehat{\mathbf{e}}_{i}\|}\right]^{j-1}
    \int_{\widehat{e}_{i}}(\widehat{\mathbf{v}}\cdot\widehat{\mathbf{t}})\,\widehat{s}^{j-1}d\widehat{s}\,,\qquad i,j=1,2,3,\nonumber\\
    \phi^{\mathscr{T}}_{k}\left(\widehat{\mathbf{v}}\right)&=\int_{\widehat{\mathscr{T}}}\left(\widehat{\mathbf{v}}\cdot\widehat{\mathbf{w}}^{\mathscr{T}}_{k}\right)\left(\det\mathbf{J}\right)d\widehat{A}\,,\qquad k=1,2,3,
\end{align}
where the vectors $\widehat{\mathbf{w}}^{\mathscr{T}}_{k},~k=1,2,3$ are defined as:
\begin{equation}\label{eqn:A9}
\widehat{\mathbf{w}}^{\mathscr{T}}_{k}=\mathbf{J}^{-1}\mathbf{w}^{\mathscr{T}}_{k}\,,
\end{equation}
with $\mathbf{w}^{\mathscr{T}}_{k},~k=1,2,3$ given in \eqref{eqn:3.6}. Observing that the Piola transformation in \eqref{eqn:A1} is linear, a form similar to the one in \eqref{eqn:3.4} can be considered for the shape functions $\widehat{\mathbf{v}}^{\widehat{\mathscr{T}},\widehat{e}_{i}}_{j},~i,j=1,2,3$ in the natural coordinate system.\footnote{The coordinate variables $x$ and $y$ are replaced with $r$ and $s$, respectively, and the coefficient $\widehat{a}^{j}_{1},\ldots,\widehat{f}^{j}_{2}$ are used instead of $a^{j}_{1},\ldots,f^{j}_{2}$.} The coefficients of the shape functions are determined by substituting $\widehat{\mathbf{v}}^{\widehat{\mathscr{T}},\widehat{e}_{i}}_{j}$ for $\widehat{\mathbf{v}}$ in \eqref{eqn:A8} and enforcing the conditions in \eqref{eqn:3.7} on the degrees of freedom $\phi^{\mathscr{T},e_{i}}_{j},~i,j=1,2,3$ and $\phi^{\mathscr{T}}_{k},~k=1,2,3$. Substituting the computed shape functions in \eqref{eqn:3.8} leads to the shape functions $\widehat{\bar{\mathbf{v}}}^{\widehat{\mathscr{T}},\widehat{e}_{i}}_{j},~i,j=1,2,3$, the explicit forms of which are provided in Table \ref{tab:A1}.
\begin{table}[h]
\centering
\caption{Explicit forms of the shape functions $\widehat{\bar{\mathbf{v}}}^{\widehat{\mathscr{T}},\widehat{e}_{i}}_{j}\in\mathcal{P}^{c}_{2}\left(T\widehat{\mathscr{T}}\right)$ on the edges $\widehat{e}_{1}$, $\widehat{e}_{2}$ and $\widehat{e}_{3}$ of the reference element $\widehat{\mathscr{T}}$ in the natural coordinate system.}
\label{tab:A1}
\renewcommand{\arraystretch}{1.5}
\renewcommand{\tabcolsep}{0.2cm}
\begin{tabular}{c c c c c}
\hline
$j$                & $\widehat{\bar{\mathbf{v}}}$ & $\widehat{e}_{1}$
                                                  & $\widehat{e}_{2}$
                                                  & $\widehat{e}_{3}$ \\
\hline
\multirow{2}{*}{1} & $\widehat{\bar{v}}_{1}$      & $-s\big[4\left(r+s\right)-3\big]$
                                                  & $s\left(4r-1\right)$
                                                  & $\left(4s-1\right)\left(s-1\right)$ \\
                   & $\widehat{\bar{v}}_{2}$      & $r\big[4\left(r+s\right)-3\big]$
                                                  & $-\left(4r-1\right)\left(r-1\right)$
                                                  & $-r\left(4s-1\right)$ \\
\multirow{2}{*}{2} & $\widehat{\bar{v}}_{1}$      & $s\left(2r+3s-2\right)$
                                                  & $s\left(2r-s\right)$
                                                  & $\left(3s-1\right)\left(2r+s-1\right)$ \\
                   & $\widehat{\bar{v}}_{2}$      & $r\left(3r+2s-2\right)$
                                                  & $\left(3r-1\right)\left(r+2s-1\right)$
                                                  & $-r\left(r-2s\right)$ \\
\multirow{2}{*}{3} & $\widehat{\bar{v}}_{1}$      & $-s\left(2s-1\right)$
                                                  & $-s\left(2s-1\right)$
                                                  & $s\left(2s-3\right)+4r\left(r+s-1\right)+1$ \\
                   & $\widehat{\bar{v}}_{2}$      & $r\left(2r-1\right)$
                                                  & $r\left(3-2r\right)-4s\left(r+s-1\right)-1$
                                                  & $r\left(2r-1\right)$ \\
\hline
\end{tabular}
\end{table}

Regarding the shape functions $\widehat{\mathbf{v}}^{\widehat{\mathscr{T}}}_{k},~k=1,2,3$ on the reference element $\widehat{\mathscr{T}}$, it can be shown that they take the following form:
\begin{equation}\label{eqn:A10}
\widehat{\mathbf{v}}^{\widehat{\mathscr{T}}}_{k}=
\begin{bmatrix}
s\left(\widehat{c}^{k}_{1}+\widehat{e}^{k}_{1}r-\widehat{c}^{k}_{1}s\right)\\
\\
r\left(\widehat{b}^{k}_{2}-\widehat{b}^{k}_{2}r+\widehat{e}^{k}_{2}s\right)
\end{bmatrix}.
\end{equation}
If we substitute $\widehat{\mathbf{v}}^{\widehat{\mathscr{T}}}_{k}$ for $\widehat{\mathbf{v}}$ in \eqref{eqn:A8} and impose the conditions \eqref{eqn:3.13} on the resulting degrees of freedom, we can determine the coefficients as given in Table \ref{tab:A2}.

\begin{table}[h]
\centering
\caption{Coefficients of the shape functions $\widehat{\mathbf{v}}^{\widehat{\mathscr{T}}}_{k}\in\mathcal{P}^{c}_{2}\left(T\widehat{\mathscr{T}}\right)$ for the reference element $\widehat{\mathscr{T}}$ in the natural coordinate system.}
\label{tab:A2}
\renewcommand{\arraystretch}{1.5}
\renewcommand{\tabcolsep}{0.2cm}
\begin{tabular}{c c c c c}
\hline
$k$ & $\widehat{c}^{k}_{1}$
    & $\widehat{e}^{k}_{1}$
    & $\widehat{b}^{k}_{2}$
    & $\widehat{e}^{k}_{2}$ \\
\hline
1   & $\frac{12}{\det\mathbf{J}}\left(3x_{2}+x_{3}\right)$
    & $\frac{-24}{\det\mathbf{J}}\left(2x_{2}+x_{3}\right)$
    & $\frac{12}{\det\mathbf{J}}\left(x_{2}+3x_{3}\right)$
    & $\frac{-24}{\det\mathbf{J}}\left(x_{2}+2x_{3}\right)$ \\
2   & $\frac{12}{\det\mathbf{J}}\left(3y_{2}+y_{3}\right)$
    & $\frac{-24}{\det\mathbf{J}}\left(2y_{2}+y_{3}\right)$
    & $\frac{12}{\det\mathbf{J}}\left(y_{2}+3y_{3}\right)$
    & $\frac{-24}{\det\mathbf{J}}\left(y_{2}+2y_{3}\right)$ \\
3   & $\frac{-60}{\det\mathbf{J}}$
    & $\frac{120}{\det\mathbf{J}}$
    & $\frac{-60}{\det\mathbf{J}}$
    & $\frac{120}{\det\mathbf{J}}$ \\
\hline
\end{tabular}
\end{table}

\subsection{$\mathcal{P}^{d}_{2}(T\mathscr{T})$ shape functions}
\label{sec:A.2}

If the unit vectors normal to the edges $e$ and $\widehat{e}$ are denoted by $\mathbf{n}$ and $\widehat{\mathbf{n}}$, respectively, then we can obtain the following relationship between the inner products $\mathbf{v}\cdot\mathbf{n}$ and $\widehat{\mathbf{v}}\cdot\widehat{\mathbf{n}}$ \citep{AZN22}:
\begin{equation}\label{eqn:A11}
\mathbf{v}\cdot\mathbf{n}=\frac{\|\widehat{\mathbf{e}}\|}{\|\mathbf{e}\|}\widehat{\mathbf{v}}\cdot\widehat{\mathbf{n}}\,.
\end{equation}
Using \eqref{eqn:A2}, \eqref{eqn:A7} and \eqref{eqn:A11}, the integrals in \eqref{eqn:3.14} can be transformed to the natural coordinate system as follows:
\begin{align}\label{eqn:A12}
\phi^{\mathscr{T},e_{i}}_{j}\left(\widehat{\mathbf{v}}\right)&=\left(\frac{\|\mathbf{e}_{i}\|}{\|\widehat{\mathbf{e}}_{i}\|}\right)^{j-1}
\int_{\widehat{e}_{i}}\left(\widehat{\mathbf{v}}\cdot\widehat{\mathbf{n}}\right)\widehat{s}^{j-1}d\widehat{s}\,,\qquad i,j=1,2,3,\nonumber\\
\phi^{\mathscr{T}}_{k}\left(\widehat{\mathbf{v}}\right)&=\int_{\widehat{\mathscr{T}}}\left(\widehat{\mathbf{v}}\cdot\widehat{\mathbf{w}}^{\mathscr{T}}_{k}\right)\left(\det\mathbf{J}\right)d\widehat{A}\,,\qquad k=1,2,3,
\end{align}
where the vectors $\widehat{\mathbf{w}}^{\mathscr{T}}_{k},k=1,2,3$ are defined as:
\begin{equation}\label{eqn:A13}
\widehat{\mathbf{w}}^{\mathscr{T}}_{k}=\frac{1}{\det\mathbf{J}}\mathbf{J}^{\mathsf{T}}\mathbf{w}^{\mathscr{T}}_{k}\,,
\end{equation}
with $\widehat{\mathbf{w}}^{\mathscr{T}}_{k},~k=1,2,3$ given in \eqref{eqn:3.15}. A form similar to \eqref{eqn:3.4} is considered for the shape functions $\widehat{\mathbf{v}}^{\widehat{\mathscr{T}},\widehat{e}_{i}}_{j},~i,j=1,2,3$ in the natural coordinate system. The coefficients of the shape functions can be determined by substituting $\widehat{\mathbf{v}}^{\widehat{\mathscr{T}},\widehat{e}_{i}}_{j}$ into \eqref{eqn:A12} and enforcing the conditions in \eqref{eqn:3.7}. The shape functions are then substituted into \eqref{eqn:3.8} to obtain the shape functions $\widehat{\bar{\mathbf{v}}}^{\widehat{\mathscr{T}},\widehat{e}_{i}}_{j},~i,j=1,2,3$. The explicit forms of $\widehat{\bar{\mathbf{v}}}^{\widehat{\mathscr{T}},\widehat{e}_{i}}_{j}$ are provided in Table \ref{tab:A3}.
\begin{table}[h]
\centering
\caption{Explicit forms of the shape functions $\widehat{\bar{\mathbf{v}}}^{\widehat{\mathscr{T}},\widehat{e}_{i}}_{j}\in\mathcal{P}^{d}_{2}\left(T\widehat{\mathscr{T}}\right)$ on the edges $\widehat{e}_{1}$, $\widehat{e}_{2}$ and $\widehat{e}_{3}$ of the reference element $\widehat{\mathscr{T}}$ in the natural coordinate system.}
\label{tab:A3}
\renewcommand{\arraystretch}{1.5}
\renewcommand{\tabcolsep}{0.2cm}
\begin{tabular}{c c c c c}
\hline
$j$                & $\widehat{\bar{\mathbf{v}}}$ & $\widehat{e}_{1}$
                                                  & $\widehat{e}_{2}$
                                                  & $\widehat{e}_{3}$ \\
\hline
\multirow{2}{*}{1} & $\widehat{\bar{v}}_{1}$      & $r\big[4\left(r+s\right)-3\big]$
                                                  & $-\left(4r-1\right)\left(r-1\right)$
                                                  & $-r\left(4s-1\right)$ \\
                   & $\widehat{\bar{v}}_{2}$      & $s\big[4\left(r+s\right)-3\big]$
                                                  & $-s\left(4r-1\right)$
                                                  & $-\left(4s-1\right)\left(s-1\right)$ \\
\multirow{2}{*}{2} & $\widehat{\bar{v}}_{1}$      & $r\left(3r+2s-2\right)$
                                                  & $\left(3r-1\right)\left(r+2s-1\right)$
                                                  & $-r\left(r-2s\right)$ \\
                   & $\widehat{\bar{v}}_{2}$      & $-s\left(2r+3s-2\right)$
                                                  & $-s\left(2r-s\right)$
                                                  & $-\left(3s-1\right)\left(2r+s-1\right)$ \\
\multirow{2}{*}{3} & $\widehat{\bar{v}}_{1}$      & $r\left(2r-1\right)$
                                                  & $r\left(3-2r\right)-4s\left(r+s-1\right)-1$
                                                  & $r\left(2r-1\right)$  \\
                   & $\widehat{\bar{v}}_{2}$      & $s\left(2s-1\right)$
                                                  & $s\left(2s-1\right)$
                                                  & $s\left(3-2s\right)-4r\left(r+s-1\right)-1$ \\
\hline
\end{tabular}
\end{table}

It is straightforward to verify that the shape functions $\widehat{\mathbf{v}}^{\widehat{\mathscr{T}}}_{k},~k=1,2,3$ on the reference element $\widehat{\mathscr{T}}$ take the following form:
\begin{equation}\label{eqn:A14}
\widehat{\mathbf{v}}^{\widehat{\mathscr{T}}}_{k}=
\begin{bmatrix}
r\left(\widehat{b}^{k}_{1}-\widehat{b}^{k}_{1}r+\widehat{e}^{k}_{1}s\right)\\
\\
s\left(\widehat{c}^{k}_{2}+\widehat{e}^{k}_{2}r-\widehat{c}^{k}_{2}s\right)
\end{bmatrix}.
\end{equation}
Substituting $\widehat{\mathbf{v}}^{\widehat{\mathscr{T}}}_{k}$ for $\widehat{\mathbf{v}}$ in \eqref{eqn:A12} and enforcing the conditions \eqref{eqn:3.13} leads to the coefficients given in Table \ref{tab:A4}.
\begin{table}[h]
\centering
\caption{Coefficients of the shape functions $\widehat{\mathbf{v}}^{\widehat{\mathscr{T}}}_{k}\in\mathcal{P}^{d}_{2}\left(T\widehat{\mathscr{T}}\right)$ for the reference element $\widehat{\mathscr{T}}$ in the natural coordinate system.}
\label{tab:A4}
\renewcommand{\arraystretch}{1.5}
\renewcommand{\tabcolsep}{0.2cm}
\begin{tabular}{c c c c c}
\hline
$k$ & $\widehat{b}^{k}_{1}$
    & $\widehat{e}^{k}_{1}$
    & $\widehat{c}^{k}_{2}$
    & $\widehat{e}^{k}_{2}$ \\
\hline
1   & $\frac{12}{\det\mathbf{J}}\left(y_{2}+3y_{3}\right)$
    & $\frac{-24}{\det\mathbf{J}}\left(y_{2}+2y_{3}\right)$
    & $\frac{-12}{\det\mathbf{J}}\left(3y_{2}+y_{3}\right)$
    & $\frac{24}{\det\mathbf{J}}\left(2y_{2}+y_{3}\right)$  \\
2   & $\frac{-12}{\det\mathbf{J}}\left(x_{2}+3x_{3}\right)$
    & $\frac{24}{\det\mathbf{J}}\left(x_{2}+2x_{3}\right)$
    & $\frac{12}{\det\mathbf{J}}\left(3x_{2}+x_{3}\right)$
    & $\frac{-24}{\det\mathbf{J}}\left(2x_{2}+x_{3}\right)$ \\
3   & $\frac{60}{\det\mathbf{J}}$
    & $\frac{-120}{\det\mathbf{J}}$
    & $\frac{-60}{\det\mathbf{J}}$
    & $\frac{120}{\det\mathbf{J}}$ \\
\hline
\end{tabular}
\end{table}

\subsection{$\mathcal{P}^{c-}_{2}(T\mathscr{T})$ shape functions}
\label{sec:A.3}

Following the same procedure as described in \S\ref{sec:A.1}, we can transform the integrals in \eqref{eqn:3.17} to the natural coordinate system as follows:
\begin{equation}\label{eqn:A15}
\begin{aligned}
    \phi^{\mathscr{T},e_{i}}_{j}\left(\widehat{\mathbf{v}}\right)&=\left[\frac{\|\mathbf{e}_{i}\|}{\|\widehat{\mathbf{e}}_{i}\|}\right]^{j-1}
    \int_{\widehat{e}_{i}}(\widehat{\mathbf{v}}\cdot\widehat{\mathbf{t}})\,\widehat{s}^{j-1}d\widehat{s}\,,&& i=1,2,3\,,~j=1,2,\\
    \phi^{\mathscr{T}}_{k}\left(\widehat{\mathbf{v}}\right)&=\int_{\widehat{\mathscr{T}}}\left(\widehat{\mathbf{v}}\cdot\widehat{\mathbf{w}}^{\mathscr{T}}_{k}\right)\left(\det\mathbf{J}\right)d\widehat{A}\,,&& k=1,2,
\end{aligned}
\end{equation}
where the vectors $\widehat{\mathbf{w}}^{\mathscr{T}}_{k},~k=1,2$ are given in \eqref{eqn:A9}. A form similar to \eqref{eqn:3.16} can be considered for the shape functions $\widehat{\mathbf{v}}^{\widehat{\mathscr{T}},\widehat{e}_{i}}_{j},~i=1,2,3,~j=1,2$ in the natural coordinate system. The coefficients are determined by substituting the shape functions into \eqref{eqn:A15} and imposing the conditions in \eqref{eqn:3.7} on the resulting degrees of freedom. The substitution of the shape functions into \eqref{eqn:3.19} leads to the shape functions $\widehat{\bar{\mathbf{v}}}^{\widehat{\mathscr{T}},\widehat{e}_{i}}_{j},~i=1,2,3,~j=1,2$ with explicit forms given in Table \ref{tab:A5}.
\begin{table}[h]
\centering
\caption{Explicit forms of the shape functions $\widehat{\bar{\mathbf{v}}}^{\widehat{\mathscr{T}},\widehat{e}_{i}}_{j}\in\mathcal{P}^{c-}_{2}\left(T\widehat{\mathscr{T}}\right)$ on the edges $\widehat{e}_{1}$, $\widehat{e}_{2}$ and $\widehat{e}_{3}$ of the reference element $\widehat{\mathscr{T}}$ in the natural coordinate system.}
\label{tab:A5}
\renewcommand{\arraystretch}{1.5}
\renewcommand{\tabcolsep}{0.2cm}
\begin{tabular}{c c c c c}
\hline
$j$                & $\widehat{\bar{\mathbf{v}}}$ & $\widehat{e}_{1}$
                                                  & $\widehat{e}_{2}$
                                                  & $\widehat{e}_{3}$ \\
\hline
\multirow{2}{*}{1} & $\widehat{\bar{v}}_{1}$      & $-s\big[4\left(r+s\right)-3\big]$
                                                  & $s\left(4r-1\right)$
                                                  & $\left(4s-1\right)\left(s-1\right)$ \\
                   & $\widehat{\bar{v}}_{2}$      & $r\big[4\left(r+s\right)-3\big]$
                                                  & $-\left(4r-1\right)\left(r-1\right)$
                                                  & $-r\left(4s-1\right)$ \\
\multirow{2}{*}{2} & $\widehat{\bar{v}}_{1}$      & $-\frac{1}{3}s\big[4\left(r-s\right)+1\big]$
                                                  & $-\frac{1}{3}s\left(4r+8s-5\right)$
                                                  & $\frac{1}{3}\left(4s-3\right)\left(2r+s-1\right)$ \\
                   & $\widehat{\bar{v}}_{2}$      & $\frac{1}{3}r\big[4\left(r-s\right)-1\big]$
                                                  & $\frac{1}{3}\left(4r-3\right)\left(r+2s-1\right)$
                                                  & $-\frac{1}{3}r\left(8r+4s-5\right)$ \\
\hline
\end{tabular}
\end{table}

The shape functions $\widehat{\mathbf{v}}^{\widehat{\mathscr{T}}}_{k},~k=1,2$ on the reference element $\widehat{\mathscr{T}}$ have the following simple form:
\begin{equation}\label{eqn:A16}
\widehat{\mathbf{v}}^{\widehat{\mathscr{T}}}_{k}=
\begin{bmatrix}
s\left( \widehat{c}^{k}_{1}-\widehat{d}^{k}_{1}r-\widehat{c}^{k}_{1}s\right)\\
\\
r\left(-\widehat{d}^{k}_{1}+\widehat{d}^{k}_{1}r+\widehat{c}^{k}_{1}s\right)
\end{bmatrix}.
\end{equation}
Substituting $\widehat{\mathbf{v}}^{\widehat{\mathscr{T}}}_{k}$ for $\widehat{\mathbf{v}}$ in \eqref{eqn:A15} and enforcing the conditions in \eqref{eqn:3.13} leads to the coefficients provided in Table~\ref{tab:A6}.
\begin{table}[h]
\centering
\caption{Coefficients of the shape functions $\widehat{\mathbf{v}}^{\widehat{\mathscr{T}}}_{k}\in\mathcal{P}^{c-}_{2}\left(T\widehat{\mathscr{T}}\right)$ for the reference element $\widehat{\mathscr{T}}$ in the natural coordinate system.}
\label{tab:A6}
\renewcommand{\arraystretch}{1.5}
\renewcommand{\tabcolsep}{0.2cm}
\begin{tabular}{c c c}
\hline
$k$ & $\widehat{c}^{k}_{1}$
    & $\widehat{d}^{k}_{1}$ \\
\hline
1   & $\frac{8}{\det\mathbf{J}}\left(2x_{2}-x_{3}\right)$
    & $\frac{8}{\det\mathbf{J}}\left(x_{2}-2x_{3}\right)$ \\
2   & $\frac{8}{\det\mathbf{J}}\left(2y_{2}-y_{3}\right)$
    & $\frac{8}{\det\mathbf{J}}\left(y_{2}-2y_{3}\right)$ \\
\hline
\end{tabular}
\end{table}

\subsection{$\mathcal{P}^{d-}_{2}(T\mathscr{T})$ shape functions}
\label{sec:A.4}

We can use the same approach as outlined in \S\ref{sec:A.2} to transform the integrals in \eqref{eqn:3.23} from the physical space to the natural coordinate system. This leads to the following degrees of freedom for the reference element~$\widehat{\mathscr{T}}$:
\begin{equation} \label{eqn:A17}
\begin{aligned}
    \phi^{\mathscr{T},e_{i}}_{j}\left(\widehat{\mathbf{v}}\right)&=\left[\frac{\|\mathbf{e}_{i}\|}{\|\widehat{\mathbf{e}}_{i}\|}\right]^{j-1}
    \int_{\widehat{e}_{i}}\left(\widehat{\mathbf{v}}\cdot\widehat{\mathbf{n}}\right)\widehat{s}^{j-1}d\widehat{s}\,,&& i=1,2,3\,,~j=1,2\,,\\
    \phi^{\mathscr{T}}_{k}\left(\widehat{\mathbf{v}}\right)&=\int_{\widehat{\mathscr{T}}}\left(\widehat{\mathbf{v}}\cdot\widehat{\mathbf{w}}^{\mathscr{T}}_{k}\right)\left(\det\mathbf{J}\right)d\widehat{A}\,,&& k=1,2,
\end{aligned}
\end{equation}
where the vectors $\widehat{\mathbf{w}}^{\mathscr{T}}_{k}\,,~k=1,2$ are given in \eqref{eqn:A13}. We can consider a form similar to \eqref{eqn:3.22} for the shape functions $\widehat{\mathbf{v}}^{\widehat{\mathscr{T}},\widehat{e}_{i}}_{j},~i=1,2,3,~j=1,2$ in the natural coordinate system. The coefficients can be computed by substituting the shape functions into \eqref{eqn:A17} and enforcing the conditions in \eqref{eqn:3.7}. These shape functions are then used in \eqref{eqn:3.19} to derive the shape functions $\widehat{\bar{\mathbf{v}}}^{\widehat{\mathscr{T}},\widehat{e}_{i}}_{j},~i=1,2,3,~j=1,2$ with explicit forms given in Table~\ref{tab:A7}.
\begin{table}[h]
\centering
\caption{Explicit forms of the shape functions $\widehat{\bar{\mathbf{v}}}^{\widehat{\mathscr{T}},\widehat{e}_{i}}_{j}\in\mathcal{P}^{d-}_{2}\left(T\widehat{\mathscr{T}}\right)$ on the edges $\widehat{e}_{1}$, $\widehat{e}_{2}$ and $\widehat{e}_{3}$ of the reference element $\widehat{\mathscr{T}}$ in the natural coordinate system.}
\label{tab:A7}
\renewcommand{\arraystretch}{1.5}
\renewcommand{\tabcolsep}{0.2cm}
\begin{tabular}{c c c c c}
\hline
$j$                & $\widehat{\bar{\mathbf{v}}}$ & $\widehat{e}_{1}$
                                                  & $\widehat{e}_{2}$
                                                  & $\widehat{e}_{3}$ \\
\hline
\multirow{2}{*}{1} & $\widehat{\bar{v}}_{1}$      & $r\big[4\left(r+s\right)-3\big]$
                                                  & $-\left(4r-1\right)\left(r-1\right)$
                                                  & $-r\left(4s-1\right)$ \\
                   & $\widehat{\bar{v}}_{2}$      & $s\big[4\left(r+s\right)-3\big]$
                                                  & $-s\left(4r-1\right)$
                                                  & $-\left(4s-1\right)\left(s-1\right)$ \\
\multirow{2}{*}{2} & $\widehat{\bar{v}}_{1}$      & $\frac{1}{3}r\big[4\left(r-s\right)-1\big]$
                                                  & $\frac{1}{3}\left(4r-3\right)\left(r+2s-1\right)$
                                                  & $-\frac{1}{3}r\left(8r+4s-5\right)$ \\
                   & $\widehat{\bar{v}}_{2}$      & $\frac{1}{3}s\big[4\left(r-s\right)+1\big]$
                                                  & $\frac{1}{3}s\left(4r+8s-5\right)$
                                                  & $-\frac{1}{3}\left(4s-3\right)\left(2r+s-1\right)$ \\
\hline
\end{tabular}
\end{table}

The shape functions $\widehat{\mathbf{v}}^{\widehat{\mathscr{T}}}_{k},~k=1,2$ on the reference element $\widehat{\mathscr{T}}$ take the following form:
\begin{equation}\label{eqn:A18}
\widehat{\mathbf{v}}^{\widehat{\mathscr{T}}}_{k}=
\begin{bmatrix}
r\left(\widehat{b}^{k}_{1}-\widehat{b}^{k}_{1}r-\widehat{c}^{k}_{2}s\right)\\
\\
s\left(\widehat{c}^{k}_{2}-\widehat{b}^{k}_{1}r-\widehat{c}^{k}_{2}s\right)
\end{bmatrix}.
\end{equation}
Substituting $\widehat{\mathbf{v}}^{\widehat{\mathscr{T}}}_{k}$ for $\widehat{\mathbf{v}}$ in \eqref{eqn:A17} and enforcing the conditions \eqref{eqn:3.13}, one can obtain the coefficients provided in Table~\ref{tab:A8}.
\begin{table}[h]
\centering
\caption{Coefficients of the shape functions $\widehat{\mathbf{v}}^{\widehat{\mathscr{T}}}_{k}\in\mathcal{P}^{d-}_{2}\left(T\widehat{\mathscr{T}}\right)$ for the reference element $\widehat{\mathscr{T}}$ in the natural coordinate system.}
\label{tab:A8}
\renewcommand{\arraystretch}{1.5}
\renewcommand{\tabcolsep}{0.2cm}
\begin{tabular}{c c c}
\hline
$k$ & $\widehat{b}^{k}_{1}$
    & $\widehat{c}^{k}_{2}$ \\
\hline
1   & $-\frac{8}{\det\mathbf{J}}\left(y_{2}-2y_{3}\right)$
    & $-\frac{8}{\det\mathbf{J}}\left(2y_{2}-y_{3}\right)$ \\
2   & $\frac{8}{\det\mathbf{J}}\left(x_{2}-2x_{3}\right)$
    & $\frac{8}{\det\mathbf{J}}\left(2x_{2}-x_{3}\right)$ \\
\hline
\end{tabular}
\end{table}

\begin{remark}\label{rem:A1}
The shape functions $\widehat{\bar{\mathbf{v}}}^{\widehat{\mathscr{T}},\widehat{e}_{i}}_{j}\in\mathcal{P}^{d}_{2}(T\widehat{\mathscr{T}})$ in Table~\ref{tab:A3} can be obtained from the $90$-degree clockwise rotation of the shape functions $\widehat{\bar{\mathbf{v}}}^{\widehat{\mathscr{T}},\widehat{e}_{i}}_{j}\in\mathcal{P}^{c}_{2}(T\widehat{\mathscr{T}})$ in Table~\ref{tab:A1}. In the same manner, a $90$-degree rotation in the clockwise direction can be applied to the shape functions $\widehat{\bar{\mathbf{v}}}^{\widehat{\mathscr{T}},\widehat{e}_{i}}_{j}\in\mathcal{P}^{c-}_{2}(T\widehat{\mathscr{T}})$ in Table~\ref{tab:A5} to obtain the shape functions $\widehat{\bar{\mathbf{v}}}^{\widehat{\mathscr{T}},\widehat{e}_{i}}_{j}\in\mathcal{P}^{d-}_{2}(T\widehat{\mathscr{T}})$ in Table~\ref{tab:A7}.
\end{remark}

\begin{remark}\label{rem:A2}
It is observed that the origin of the coordinate system for the simplicial element in Figure \ref{fig:21} is attached to the first vertex of the edge $e_{3}$. Therefore, to compute the shape functions $\widehat{\bar{\mathbf{v}}}^{\widehat{\mathscr{T}},\widehat{e}_{i}}_{j}$ on the edge $e_{i}$ we can set the origin of the coordinate system to the first vertex of that edge and use the shape functions corresponding to the edge $e_{3}$ in Tables \ref{tab:A1}, \ref{tab:A3}, \ref{tab:A5} and \ref{tab:A7} to evaluate the desired shape functions. However, it is necessary to apply the cyclic permutations $e_{3}\rightarrow e_{1}\rightarrow e_{2}$ and $v_{2}\rightarrow v_{3}\rightarrow v_{1}$ to the Jacobian in \eqref{eqn:A3} to correctly transform the shape functions from the natural coordinate system to the physical space.
\end{remark}

\end{document}